\documentclass[12pt]{article}
\usepackage{amsfonts}
\usepackage{mathrsfs}
\usepackage{amsmath}
\usepackage{latexsym}
\usepackage{amssymb}
\usepackage{amsthm}
\usepackage{amscd}
\title{\bf Nichols algebras over classical Weyl groups, Fomin-Kirillov algebras and Lyndon basis
}
\author{ \small  Shouchuan Zhang $^{a}$,  Weicai Wu $^{b}$,  Zhengtang Tan $^{a}$    and Yao-Zhong Zhang $^{c,d}$   \\
\small $a$. Department  of Mathematics,   Hunan University,
\small   Changsha  410082,   P.R. China \\
\small $b$. Department  of Mathematics,   Zhejiang  University,
\small   Hangzhou  310007,  P.R. China \\
\small $c$. School of Mathematics and Physics,   The University of Queensland\\
\small Brisbane 4072,   Australia\\
\small $d$. CAS Key Laboratory of Theoretical Physics, Institute of Theoretical Physics\\
\small  Chinese Academy of Sciences, Beijing 100190, P.R. China\\
\small {\tt Emails: z9491@sina.cn (SZ); yzz@maths.uq.edu.au (YZZ)}
 }
\date {}
\begin{document}
\newtheorem{Proposition}{Proposition}[section]
\newtheorem{Theorem}[Proposition]{Theorem}
\newtheorem{Definition}[Proposition]{Definition}
\newtheorem{Corollary}[Proposition]{Corollary}
\newtheorem{Lemma}[Proposition]{Lemma}
\newtheorem{Example}[Proposition]{Example}
\newtheorem{Remark}[Proposition]{Remark}
\newtheorem{Conjecture}[Proposition]{Conjecture}

\maketitle

\begin{abstract} We show that  except in several  cases conjugacy classes of classical Weyl groups $W(B_n)$
and $W(D_n)$ are of type {\rm D}. We prove that except in three  cases Nichols algebras of irreducible
Yetter-Drinfeld ({\rm YD} in short )modules over the classical Weyl groups are infinite dimensional.
We establish the relationship between Fomin-Kirillov algebra $\mathcal E_n$
and Nichols algebra $\mathfrak{B} ({\mathcal O}_{{(1, 2)}} ,  \epsilon \otimes {\rm sgn})$
of transposition over symmetry group by means of quiver Hopf algebras. We generalize {\rm FK } algebra.
The   characteristic of finiteness of Nichols algebras in thirteen  ways  and of  {\rm FK } algebras
${\mathcal E}_n$ in nine   ways is given. All irreducible representations of finite dimensional Nichols algebras
and a complete  set of hard super- letters of Nichols algebras of finite Cartan types are  found.
The sufficient and necessary condition for Nichols algebra $\mathfrak B(M)$ of reducible {\rm YD} module  $M$
over  $A \rtimes \mathbb{S}_n$ with ${\rm supp } (M) \subseteq A$
to be finite dimensional is given.
It is shown that hard braided Lie Lyndon word,  standard  Lyndon word,  Lyndon basis path, hard  Lie Lyndon word
and standard  Lie Lyndon word are the same with respect to $ \mathfrak B(V)$, Cartan matrix $A_c$
and $U(L^+)$, respectively, where  $V$ and $L$ correspond to the same finite Cartan matrix $A_c$.



\end {abstract}


\section{Introduction}\label{introduction}

Nichols algebras play a fundamental role in the classification of finite-dimensional complex
pointed Hopf algebras by means of the lifting method developed by Andruskiewitsch and Schneider
\cite {AS02}. In this context, given a group $ G$, an important step to classify all finite-dimensional complex pointed
Hopf algebras $H$ with group-like $G(H) = G$ is to determine all the pairs $(\mathcal O, \rho )$ such that the
associated Nichols algebra $  \mathfrak B(\mathcal O, \rho )$ is finite-dimensional. Here the pairs $(\mathcal O, \rho )$ are such that
$\mathcal O$ runs over all conjugacy classes of $G$ and $\rho $ runs over all irreducible representations of the
centralizer of $g$ in $ G$, with $g \in \mathcal O$ fixed. In general, this is a difficult task since Nichols algebra
is defined by generators and relations.
In practice, it is often useful to discard those pairs such that $\dim  \mathfrak B (\mathcal O, \rho ) = \infty$. There are
properties of the conjugacy class $\mathcal O$ that imply that $\dim  \mathfrak B (\mathcal O, \rho )  = \infty$ for any $\rho $, one
of which is the property of being of type {\rm D}. This is useful since it reduces the computations to
operations inside the group and avoids hard calculations of generators and relations of the
corresponding Nichols algebra.

 Fomin and Kirillov \cite [Conjecture 2.2]{FK99} point out a conjecture `` ${\mathcal E}_n$ is finite dimensional " to study  the cohomology ring of the flag manifold. In \cite {FK99} it is shown that
${\mathcal E}_3 = 12 $  and ${\mathcal E}_4 = 24 ^2 $ (see \cite [Section 3.4]{AS02}). Many papers ( for example, \cite {   MS00, GHV11, Ba06, Ma10}) refer to this conjecture.

Kharchenko and Heckenberger give hard Lyndon bases of   Nichols algebras of diagonal type  in \cite  {Kh99}  and
  \cite {He05}. All standard Lyndon words  in $\mathfrak B(V)$ were obtained when $V$ is finite Cartan type in   \cite [Section 4.2]{AA08}.

This work contributes to the classification of finite-dimensional Hopf algebras over an algebraically closed field  of characteristic $0$. This problem was
posed by Kaplansky in 1975. The lifting method of
Andruskiewitsch and Schneider  describes
a way to classify finite-dimensional complex Hopf algebras. In \cite {AS10} they
obtained the classification over the finite abelian groups whose order is
relative prime with 210.
Nichols algebras of braided vector spaces  $(\mathbb CX,  cq)$,  where $X$ is a
rack and $q$ is a 2-cocycle in $X$, were studied in \cite {AG03}.  It was shown
\cite {AFGV08,  AFZ09,  AZ07} that Nichols algebras  $\mathfrak B({\mathcal O}_
\sigma, \rho)$ over symmetry groups  are infinite dimensional,  except in a
number of remarkable cases corresponding to ${\mathcal O}_\sigma$. Two of the
present authors \cite{ZZ12} showed  that except in three cases Nichols algebras
of irreducible Yetter-Drinfeld (YD) modules over classical Weyl groups $A
\rtimes \mathbb S_n$ supported by $\mathbb S_n$ are infinite dimensional.
However, the classification has not been completed for Nichols algebras over
general classic Weyl groups $W(B_n)$ and $W(D_n)$. In  \cite [Section 4.2]{AA08} and \cite {LR95},  all standard Lyndon words  in $\mathfrak B(V)$ and $U(L^*)$ were obtained when $V$ is finite Cartan type and $L$ is a simple Lie algebra.

Note that $\mathbb Z_2^n \rtimes \mathbb S_n$
is isomorphic to Weyl groups $W(B_n)$ and $W(C_n)$ of $B_n$ and
$C_n$ for $n>2$. If $A= \{a \in
     \mathbb Z_2^n \mid \ a = ({a_1}, {a_2}, \cdots, {a_n}) $  $\hbox {
with   } a_1 +a_2 + \cdots + a_n \equiv 0  \ ({\rm mod } \ 2) \}$,
then $A \rtimes \mathbb S_n$ is isomorphic to Weyl group
$W(D_n)$ of $D_n$ for $n>3$. Obviously, when $A= \{a \in
     {\mathbb Z}_2^n \mid \ a = ({a_1}, {a_2}, \cdots, {a_n}) $  $\hbox {
with  all } a_i \equiv 0 \ ({\rm mod } \ 2)\}$, $A \rtimes \mathbb S_n$ is isomorphic to
Weyl group of $A_{n-1}$ for $n>1$. Note that $\mathbb{S}_n$ acts on
$A$ as follows: for any $a\in A$ with $a = ({a_1}, {a_2},
\cdots, {a_n})$ and $\sigma\in \mathbb{S}_n$,
$\sigma \cdot a :=( {a_{\sigma^{-1}(1)}}, {a_{\sigma^{-1}(2)}}, \cdots,{a_{\sigma^{-1}(n)}} ).$
It is clear that
 \begin {eqnarray*}
&& (a, \sigma )^{- 1} = (- (a _{\sigma (1)}, a_{\sigma (2)}, \cdots  , a_{\sigma (n)}), \sigma^{ - 1}) = (- \sigma ^{ - 1}(a), \sigma ^{- 1}), \\
&&(b, \tau  )(a, \sigma )(b, \tau  )^{- 1} = (b + \tau  (a) -  \tau  \sigma \tau ^{ - 1}(b), \tau  \sigma \tau ^{ - 1}).
\end  {eqnarray*}
Without specification,  $K_n := \{ a\in \mathbb Z_n \mid a_1 +a_2 +\cdots + a_n =0\}$ and $ K_n \rtimes  \mathbb{S}_n$  is a subgroup of $\mathbb Z_2^n \rtimes \mathbb S_n$. Let $\mathbb W_n$ denote $  K_n \rtimes  \mathbb{S}_n$  or $\mathbb Z_2^n \rtimes \mathbb S_n$ throughout this paper.

In this paper we prove  that except in several  cases conjugacy classes of classical Weyl groups $\mathbb W_n$ are of type
{\rm D}, and  except in three  cases Nichols algebras of irreducible
YD modules over the classical Weyl groups are infinite dimensional.
We also establish the relationship
between the Fomin-Kirillov ({\rm FK}) algebra $\mathcal E_n$ introduced in \cite{FK99} and the Nichols algebra
$\mathfrak{B} ({\mathcal O}_{{(1, 2)}} ,  \epsilon \otimes {\rm sgn})$ of transposition over symmetry group with the help of quiver Hopf algebras. That is, if  $\dim \mathfrak{B} ({\mathcal O}_{{(1, 2)}} , \epsilon\otimes {\rm sgn}) = \infty$,  then so is  $\dim \mathcal E_n$. We give the characteristic of finiteness of Nichols algebras in thirteen ways  and of ({\rm FK}) algebra $\mathcal E_n$ in nine   ways, and found all irreducible representations of finite dimensional Nichols
 algebras and a complete set of hard super- letters of Nichols algebras of finite Cartan type. We give the sufficient and necessary condition for Nichols algebra
 $\mathfrak B(M)$ of reducible {\rm YD} module  $M$ over  $A \rtimes \mathbb{S}_n$
 with ${\rm supp } (M) \subseteq A$  to be finite dimensional.
Some conditions for a braided vector space to become a {\rm YD} module over finite commutative group are also obtained. Finally we show that hard braided Lie Lyndon word,  standard  Lyndon word,  Lyndon basis path, hard  Lie Lyndon word and    standard  Lie Lyndon word    are the same with respect to $ \mathfrak B(V)$, Cartan matrix $A_c$
and $U(L^+)$, respectively, where  $V$ and $L$ correspond to the same finite Cartan matrix $A_c$.



The work is organized as follows. In Section \ref{notation} we provide some preliminaries and set our notations.
In Section \ref{extension}  we determine when the conjugacy classes of  juxtapositions of two elements
are of type {\rm D}. In  Section \ref{rack}  we  prove that except in several  cases conjugacy classes of classical Weyl groups $\mathbb W_n$ are of type {\rm D}.
In  Section \ref{nichols algebra} we classify Nichols algebras of irreducible {\rm YD} modules over the classical Weyl groups.
In  Section \ref{fk conjecture} we give the relationship between $\mathfrak{B}({\mathcal O}_{\sigma}, \rho)$ and ${\mathcal E}_n$,  where $\rho= {\rm sgn} \otimes {\rm sgn}$ or $\rho= \epsilon \otimes {\rm sgn}$.
 In Section \ref{bi-one arrow}  we give an estimate for the dimensions of the {\rm PM} Nichols  algebras and {\rm FK} algebra $\mathcal E_n$. In Section \ref{finiteness} the characteristic of finiteness of Nichols algebras and  {\rm FK} algebras is given in several ways. In Section \ref{irreps section}  all irreducible representations of finite dimensional Nichols algebras are found. In Section \ref{YD modules} we obtain some conditions for a braided vector space to become a {\rm YD} module over finite commutative group, and give a sufficient and necessary condition for Nichols algebra $\mathfrak B(M)$ of reducible
{\rm YD} module  $M$ over  $A \rtimes \mathbb{S}_n$ with ${\rm supp } (M) \subseteq A$
 to be finite dimensional. In Section \ref{Lyndon basis} we show that hard braided Lie Lyndon word,  standard  Lyndon word,  Lyndon basis path, hard  Lie Lyndon word and    standard  Lie Lyndon word    are the same with respect to $ \mathfrak B(V)$, Cartan matrix $A_c$ and $U(L^+)$, respectively.

\section*{Preliminaries and conventions}\label{notation}

A quiver $Q=(Q_0, Q_1, s, t)$ is an oriented graph,  where  $Q_0$ and
$Q_1$ are the sets of vertices and arrows,  respectively; $s$ and $t$
are two maps from  $Q_1$ to $Q_0$. For any arrow $a \in Q_1$,  $s(a)$
and $t(a)$ are called {\it its start vertex and end vertex},  respectively,
and $a$ is called an {\it arrow} from $s(a)$ to $t(a)$. For any $n\geq 0$,
an $n$-path or a path of length $n$ in the quiver $Q$ is an ordered
sequence of arrows $p=a_na_{n-1}\cdots a_1$ with $t(a_i)=s(a_{i+1})$
for all $1\leq i\leq n-1$. Note that a 0-path is exactly a vertex
and a 1-path is exactly an arrow. In this case,  we define
$s(p)=s(a_1)$,  the start vertex of $p$,  and $t(p)=t(a_n)$,  the end
vertex of $p$. For a 0-path $x$,  we have $s(x)=t(x)=x$. Let $Q_n$ be
the set of $n$-paths. Let $^yQ_n^x$ denote the set of all $n$-paths
from $x$ to $y$,  $x,  y\in Q_0$. That is,  $^yQ_n^x=\{p\in Q_n\mid
s(p)=x,  t(p)=y\}$.

Let ${ k}$ be the complex field.
Let $G$ be a group, $\widehat{G}$ denote the set of all isomorphism classes
of irreducible representations of $G$, $G^\sigma$ be the centralizer of
$\sigma$, and ${\mathcal O}_\sigma$ or ${\mathcal O}_\sigma^G$ be the
conjugacy class of $\sigma$ in $G$. We use the notation in \cite {Ka95, Sw69a,Mo93}.

Let ${\mathcal K}(G)$ denote the set of conjugacy classes in group $G$. A
formal sum $r=\sum_{C\in {\mathcal K}(G)}r_CC$  of conjugacy classes
of $G$  with cardinal number coefficients is called a {\it
ramification} (or {\it ramification data} ) of $G$,  i.e.  for any
$C\in{\mathcal K}(G)$,  \  $r_C$ is a cardinal number. In particular,
a formal sum $r=\sum_{C\in {\mathcal K}(G)}r_CC$  of conjugacy
classes of $G$ with non-negative integer coefficients is a
ramification of $G$.

 For any ramification $r$ and  $C \in {\mathcal K}(G)$ we can choose a set $I_C(r)$ such that its cardinal number is $r_C$
without loss of generality.
 Let ${\mathcal K}_r(G):=\{C\in{\mathcal
K}(G)\mid r_C\not=0\}=\{C\in{\mathcal K}(G)\mid
I_C(r)\not=\emptyset\}$.  If there exists a ramification $r$ of $G$
such that the cardinal number of $^yQ_1^x$ is equal to $r_C$ for any
$x,  y\in G$ with $x^{-1}y \in C\in {\mathcal K}(G)$,  then $Q$ is
called a {\it Hopf quiver with respect to the ramification data
$r$}. In this case,  there is a bijection from $I_C(r)$ to $^yQ_1^x$,
and hence we write  ${\ }^yQ_1^x=\{a_{y, x}^{(i)}\mid i\in I_C(r)\}$
for any $x,  y\in G$ with $x^{-1}y \in C\in {\mathcal K}(G)$. If $r_C =1$ for any  $C \in {\mathcal
K}_r(G)$, then  the arrow  from $x$ to $y$ is denoted by $a_{y, x}$ in short. Let $\underline {a_{i_m, i_{m-1}}a_{i_{m-1}, i_{m-2}} \cdots a_{i_2, i_{1}}}$ denote a path from $i_1$ to $i_m.$

$(G,  r,  \overrightarrow \rho,  u)$ is called a {\it ramification system
with irreducible representations}  (or {\rm RSR } in short),  if $r$
is a ramification of $G$, $u$ is a map from ${\mathcal K}(G)$ to $G$
with $u(C)\in C$ for any $C\in {\mathcal K}(G)$;  $I_{C} (r,  u)$ is a set  and
$\overrightarrow \rho=\{\rho_C^{(i)} \}_ { i\in I_C(r, u ),
C\in{\mathcal K}_r(G)} \ \in \prod _ { C\in{\mathcal K}_r(G)}
(\widehat{    {G^{u(C)}}}) ^{\mid I_{C} (r,  u) \mid }$ with
$\rho_C^{(i)} \in \widehat{ {G^{u(C)}}} $ for any $i$ in a set
 $I_C(r,  u),  C\in {\mathcal K}_r(G)$.
In this paper we assume that $I_C(r,  u)$ is a finite set for any $C\in
{\mathcal K}_r (G).$
 Furthermore,  if $\rho _C^{(i)}$ is a one dimensional representation for any $C\in
{\mathcal K}_r(G)$,  then  $(G,  r,  \overrightarrow \rho,  u)$ is called a {\it ramification system
with characters }  (or {\rm RSC }$(G,  r,  \overrightarrow \rho,  u)$ in short) (see \cite [Definition 1.8]{ZZC07}).

For ${\rm RSR}(G,  r,  \overrightarrow \rho,  u)$,  let $\chi
_{C}^{(i)}$ denote the character of $\rho _C^{(i)}$ for any $i \in
I_C(r,  u)$,  $C \in {\mathcal K}_r(G)$. If $r = r_CC$
and $I_C(r, u) = \{ i  \}$ then we say that ${\rm RSR}(G,  r,  \overrightarrow \rho,  u)$ is   bi-one,  written  as
${\rm RSR } (G,  {\mathcal O}_{s},  \rho )$ with $s= u(C)$ and $\rho = \rho _C^{(i)}$ in short,
since $r$ only has one conjugacy class $C$ and $\mid\!I_C(r, u)\!\mid
=1$. Quiver Hopf algebras,   Nichols algebras and
Yetter-Drinfeld modules,  corresponding to a bi-one ${\rm RSR}(G,  r,
\overrightarrow \rho,  u)$,  are said to be    bi-one.
Therefore  we also say that
$\mathfrak{B}({\mathcal O}_s, \rho)$ is  a bi-one Nichols Hopf
algebra.

\begin {Remark} \label {1.2'} The representation $\rho$ in   $\mathfrak{B}({\mathcal O}_s, \rho)$   introduced in \cite {Gr00, AZ07} and
$\rho _C^{(i)}$ in {\rm RSR} are different.  $\rho (g)$ acts on its representation
space from  the left and $\rho _C^{(i)} (g)$
acts on its representation space from  the right.
\end {Remark}

Otherwise, when $\rho = \chi$ is a one dimensional representation, then $(kQ_1^1, ad( G, r, \overrightarrow{\rho}, u)) $ is PM (see \cite [Def. 1.1] {ZZC07}). Thus the formulae  are available in   \cite [Lemma 1.9] {ZZC07}. That is, $g \cdot a _{t} = a _{g t_i, g}$, $  a _{ t_i} \cdot g = \chi (\zeta _{i}(g)) a_{t_ig, g}$.

Let $S_{m}\in {\rm End}_{k}(T(V)^{m})$
and $S_{1, j}\in {\rm End}_{k}(T(V)^{j+1})$ denote the maps $S_{m}=\prod \limits_{j=1}^{m-1}({\rm id}^{\bigotimes m-j-1} \\ \bigotimes S_{1, j})$ ,
 $S_{1, j}={\rm id}+C_{12}^{-1}+C_{12}^{-1}C_{23}^{-1}+\cdots+C_{12}^{-1}C_{23}^{-1}\cdots C_{j, j+1}^{-1}$
(in leg notation) for $m\geq 2$ and $j\in \mathbb N$. Then the subspace
$S=\bigoplus \limits _{m=2}^{\infty}  {\rm ker} S_{m}$ of the tensor $T(V)=\bigoplus \limits _{m=0}^{\infty} T(V)^{\bigotimes m}$
is a two-sided ideal,  and algebra $\mathfrak B(V)=T(V)/S$ is termed the Nichols algebra associated to $(V, C)$.

For  $s\in G$ and  $(\rho,  V) \in  \widehat {G^s}$,  here is a
precise description of the {\rm YD} module $M({\mathcal O}_s,
\rho)$,  introduced in \cite {Gr00}. Let $t_1 = s, t_2,   \cdots,
t_{m}$ be a numeration of ${\mathcal O}_s$,  which is a conjugacy
class containing $s$,   and let $g_i\in G$ such that $g_i \rhd s :=
g_i s g_i^{-1} = t_i$ for all $1\le i \le m$. Then $M({\mathcal
O}_s,  \rho) = \oplus_{1\le i \le m}g_i\otimes V$. Let $g_iv :=
g_i\otimes v \in M({\mathcal O}_s, \rho)$,  $1\le i \le m$,  $v\in V$.
If $v\in V$ and $1\le i \le m $,  then the action of $h\in G$ and the
coaction are given by
\begin {eqnarray} \label {e0.11}
\delta(g_iv) = t_i\otimes g_iv,  \qquad h\cdot (g_iv) =
g_j(\nu _i(h)\cdot v),  \end {eqnarray}
 where $hg_i = g_j\nu _i(h)$,  for
unique  $1\le j \le m$ and $\nu _i(h)\in G^s$. The explicit formula for
the braiding is then given by
\begin{equation} \label{yd-braiding}
C(g_iv\otimes g_jw) = t_i\cdot(g_jw)\otimes g_iv =
g_{j'}(\nu _j (t_i)\cdot w)\otimes g_iv\end{equation} for any $1\le i, j\le
m$,  $v, w\in V$,  where $t_ig_j = g_{j'}\nu _j(t_i)$ for unique $j'$,  $1\le
j' \le m$ and $\nu _j(t_i) \in G^s$. Let $\mathfrak{B} ({\mathcal O}_s,
\rho )$ denote $\mathfrak{B} (M ({\mathcal O}_s,  \rho ))$.
$M({\mathcal O}_s,  \rho )$ is a simple {\rm YD} module (see \cite
 {DPR91, Ci97,    AZ07}).

We briefly recall the definition and main properties of racks;
see \cite {AG03} for details, more information and bibliographical references.
A rack is a pair $(X, \rhd )$,  where $ X$
is a non-empty set and  $\rhd :  X \times X \rightarrow  X$ is an operation such that
  $ x\rhd x = x$, $x\rhd (y \rhd z) = (x\rhd y)\rhd (x\rhd z)$  and $\phi _x$
is invertible for any $x, y, z
\in  X$, where $\phi _x$ is a map from $X$ to $X$ sending $y$ to $x\rhd y$  for any $x, y\in X.$

For example, $({\mathcal O}_s^G, \rhd )$ is a rack  with $x\rhd y:= x y x^{-1}.$

 If $R$ and $S$ are two subracks of $X$  with $R\cup S = X$, $R\cap S =\emptyset$,  $x\rhd y \in S$, $y\rhd x \in R$,  for any $x\in R, y\in S,$ then $R\cup S$ is called a decomposition of subracks of $X$. Furthermore, if there exists $a\in R$, $b\in S$ such that ${\rm sq} (a, b) := a\rhd (b \rhd (a \rhd b)) \not= b$, then $X$ is called to be of type {\rm D}. Notice that if a rack $Y$ contains a subrack  $X$ of type ${\rm D}$, then $Y$  is also of  type ${\rm D}$ ( see \cite {AFGV08}).

Now we keep on the work in \cite[Page 295-299 ]{Su78}.
 Let $Y_i $ be the set of all letters which
belong to those cycles of length $i$ in the independent  cycle
decomposition of $\sigma$. Clearly,  $Y_{i}\bigcap Y_{j}=\emptyset$
for $i\neq j$ and $\bigcup\limits_iY_i=\{1, 2, \cdots, n \}$.

By \cite {ZWW08} and \cite[Page 295]{Su78}, we have

  {(i)} $\rho \in (\mathbb S_n)^\sigma $ if and only
if  $Y_{i}$ is $\rho$-invariant,  namely $\rho(Y_i)\subset Y_i$, and
 the restriction $\rho_{i}$ of $\rho$ on $Y_{i}$ commutes with the
restriction $\sigma_{i}$ of $\sigma$ on $Y_{i}$ for  $i = 1,  2,
\cdots,  n$;

 {\rm (ii)}  $(\mathbb S_n)^\sigma= \prod\limits_{i} ^n(\mathbb S_{Y_i})^{\sigma_i};$.

 {\rm (iii)}  $\widehat {  (\mathbb S_n)^\sigma }= \otimes _{i=1} ^n   \widehat {   (\mathbb S_{Y_i})^{\sigma_i}};$.

 Let $G= \mathbb W_n$. $(a, \sigma)\in G$ is called a
sign cycle if $\sigma = (i_1, i_2, \cdots, i_r)$  is cycle and $a =
(g_2 ^{a_1}, \cdots, g_2 ^{a_n})$ with $a_i=0$ for $i \notin \{i_1,
i_2, \cdots, i_r\}$.  A sign cycle $(a, \sigma)$ is called positive
( or negative ) if $\sum _{i=1}^n a_i $ is even (or odd).

$B = \oplus _{i=0} ^\infty  B_i$ is  called  a ${\mathbb N}_0 $-graded vector space  or a graded vector space if $B = \oplus _{i=0} ^\infty  B_i$ is a direct sum as vector spaces. Let $B_{>0}$ or $B^+ $  denote $ \oplus _{i=1} ^\infty  B_i.$
If $B$ is an algebra and $B = \oplus _{i=0} ^\infty  B_i$ is a graded vector space  with unit element $1_A \in B_0$ and $B_iB_j \subseteq B_{i+j}$ for any $i, j \ge 0$, then $B$ is called a ${\mathbb N}_0 $-graded algebra, or graded algebra in short.

 If $M$ is a module over a graded algebra $B = \oplus _{i=0} ^\infty  B_i$ with $B_0=k$ and $ h \cdot x = h_0 x $ for any $h = \sum \limits_{i=0} ^\infty h_i \in B$ with $h_i \in B_i$ $ (0\le i ) $ and $x\in M$, then $M$ is called a $B$-module with trivial action.

 If   $\{x_1,   \cdots,   x_n\}$ is  a basis of  vector space  $V$ and
$C(x_i\otimes  x_j) = q_{ij} x_j\otimes x_i$ with $q_{ij} \in F$,
then $V$  is called a braided vector space of diagonal type, $\{x_1,   \cdots,   x_n\}$
is called canonical basis and $(q_{ij})_{n\times n}$  is called braided matrix.

 We use  notations in  \cite {Sw69a},  \cite {Mo93}, \cite {Ma95} and \cite {Ka95}.

\section{Conjugacy classes of juxtapositions}\label {extension}
In this section we determine when the conjugacy classes of  juxtapositions of two elements are of type {\rm D}.

Let $\alpha := (1, 1, \cdots, 1)\in \mathbb Z_2^n.$

\begin {Lemma}\label {1.1} Let $G= \mathbb W_n$ and $(a, \tau), (b, \mu)\in G$. Let $(c, \lambda)$ denote $(a, \tau) \rhd ((b, \mu) \rhd((a, \tau)\rhd (b, \mu)))$. Then

{\rm (i) }   
$c =(a+\tau \cdot [ b+\mu  \cdot(a+\tau \cdot b+(\tau \rhd \mu )\cdot a)+ (\mu \rhd (\tau \rhd \mu)) \cdot b ] 
+(\tau \rhd (\mu \rhd (\tau \rhd \mu))) \cdot a.$

 {\rm (ii)} If  $\tau $ and $\mu$ are commutative,  then 
 $c=a+\tau \mu \cdot a+\tau \mu^2 \cdot a+\mu \cdot a+\tau \cdot b+\tau ^2 \mu \cdot b+\tau \mu \cdot b.$ 

{\rm (iii)} If $\tau $ and $\mu$ are commutative,  then ${\rm sq }(a\tau,  b\mu) = b\mu $ if and only if
\begin {eqnarray}\label {e2.1.3}
a+\tau \mu \cdot a+\tau \mu^2 \cdot a+\mu\cdot a=b+\tau \cdot b+\tau ^2 \mu \cdot b+\tau \mu \cdot b.
\end {eqnarray}

{\rm (iv)} If $\tau $ and $\mu$ are commutative with  $(b, \mu)=\xi\rhd (a, \tau )$ and $\xi \cdot a=a$,  then $c=a+\tau \mu^2 \cdot a+\mu \cdot a+\tau \cdot a+\tau ^2 \mu \cdot a$.

{\rm (v)} If  $\tau $ and $\xi$ are commutative with $\tau ^2=1$ and $(b, \mu)=\xi \rhd (a, \tau )$ and $\xi\cdot a=a$,  then $c=a$.
\end {Lemma}
\noindent {\it Proof.}  It is clear. \hfill $\Box$

\vskip.1in
If the lengths of independent sign cycles of $(a, \pi)$ and $(b, \tau)$ are different,  then they are called mutually orthogonal,  written as $a\pi \bot b\tau. $ Obviously, $a\pi \bot b\tau $ if and only if $\pi \bot \tau. $ If ${\rm ord (\pi)}$ and ${\rm ord (\tau)}$ are coprime  and one of the two elements does not have any fixed point, then  $\pi \bot \tau. $

For any $a\pi\in \mathbb W_n $ and $b\tau\in \mathbb W_m$,  define $a\pi \# b\tau \in \mathbb W_{m+n}$ as follows:

$(a\# b)_{i} :=\left
\{\begin
{array} {ll}  a_i &\hbox {when } i\le n\\
b_{i-n} &\hbox {when } i>n
\end {array} \right.,  $  $  (\pi \# \tau ) (i):=\left
\{\begin
{array} {ll} \pi (i) &\hbox {when } i\le n\\
\tau (i-n)+n &\hbox {when } i>n
\end {array} \right. $
\hbox { and } $ a\pi\#b \tau : = (a\# b,   \pi \# \tau ) $.
Obviously $a\pi\#b \tau \in \mathbb W_{m+n}$ and it is called a juxtaposition of $a\pi$ and $b \tau$  (see \cite {AZ07} ) ($\#$ in this place is not smash product). Let $\overrightarrow{\nu_{n,  m}}$ be a map from $\mathbb W_n$ to $\mathbb W_{m+n}$ by sending $a\pi$ to $\overrightarrow{\nu_{n,  m}}(a\pi) :=a\pi \# 1_{\mathbb W_m}$; let $\overleftarrow {\nu_{n,  m}}$ be a map from $\mathbb W_m$ to $\mathbb W_{m+n}$ by sending $b\tau$ to $\overleftarrow{\nu_{n,  m}}(b\tau) :=1_{\mathbb W_n} \# b \tau$.

\begin {Lemma}\label {2.4} Assume  $a\pi \bot b\tau$ with $a\pi,  a'\pi'\in \mathbb W_n $ and $b\tau,  b'\tau'\in \mathbb W_m$. Then

 {\rm (i)} $(a\pi\#b\tau)(a'\pi'\#b'\tau')=(a\pi a'\pi'\# b\tau b'\tau')$.

 {\rm (ii) } $a\pi\#b\tau=\overrightarrow{\nu_{n, m}}(a\pi)\overleftarrow{\nu_{n, m}}(b\tau)=\overleftarrow{\nu_{n,  m}}(b\tau)
\overrightarrow{\nu_{n, m}}(a\pi)$.

{\rm (iii)} $\mathbb W_{m+n}^{a\pi\#b\tau}=\mathbb W_n^{a\pi}\#\mathbb W_m^{b\tau}=\overrightarrow{\nu_{n, m}}(\mathbb W_{n} ^{a\pi})\overleftarrow{\nu_{n, m}}(\mathbb W_m^{b\tau})$ as directed products.

{\rm (iv)} For any $\rho\in\widehat{\mathbb W_{m+n}^{a\pi\# b\tau}}$,  there exist $\mu\in\widehat{\mathbb W_{n}^{a\pi}}$,  $\lambda   \in\widehat{\mathbb W_{m}^{b\tau}}$  such that $\rho=\mu\otimes\lambda $.

{\rm (v)} $(a\pi\# b\tau)\rhd(a'\pi'\# b'\tau')=(a\pi \rhd a'\pi')\#(b\tau\rhd b'\tau')$.

{\rm (vi)} $\mathcal O_{a\pi\# b\tau}^{\mathbb W_{m+n}}=\mathcal O_{a\pi}^{\mathbb W_n}\#\mathcal O_{b\tau}^{\mathbb W_m}$.
\end {Lemma}
\noindent {\it Proof.} {\rm (i)},  {\rm (ii)} and {\rm (v)} are clear.

{\rm (iii)} By \cite [Section 2.2] {AZ07},  $\mathbb S_{m+n}^{\pi\#\tau}=\mathbb S_n^{\pi}\#\mathbb S_m^{\tau}.$
Obviously,  $\mathbb W_n^{a\pi}\#\mathbb W_m^{b\tau}\subseteq \mathbb W_{m+n} ^{a\pi \# b\tau} $. For any
$c\xi\in\mathbb W_{m+n}^{a\pi\# b\tau}, $ then $\xi\in\mathbb S_{m+n}^{\pi \#\tau}$ and there exist $\mu\in\mathbb S_n^\pi$ and $\lambda\in\mathbb S_m ^\tau $ such that $\xi=\mu \# \lambda$. Consequently,  $c\xi=d\mu\#f\lambda.$ Considering $c\xi(a\pi \# b \tau)=(a\pi \# b\tau)c\xi$,  we have $d\mu \in \mathbb W_n^{a\pi}$ and
 $f\lambda  \in \mathbb W_m^{b\tau}$. This completes the proof.

{\rm (iv)} It follows from {\rm (iii)}.

(vi) By {\rm (v)},  $  \mathcal O_{a\pi}^{\mathbb W_{n}}\#\mathcal O_{ b\tau}^{\mathbb W_{m}}\subseteq\mathcal O_{a\pi\# b
\tau}^{\mathbb W_{m+n}}. $  Consequently {\rm (vi)} follows from {\rm (iii)}. \hfill $\Box$

\begin{Remark} {\rm (i),  (ii)} and {\rm (v)} above still hold when $a\pi$ and $b\tau$ are not mutually orthogonal. \end{Remark}

\begin {Theorem}\label {2.5} If $\mathcal O_{a\tau}$ is of type {\rm D},  then $\mathcal O_{a\tau \# b\mu}$ is also of type {\rm D}.
\end {Theorem}
\noindent {\it Proof.} Let $X= R\cup S$ be a subrack decomposition of $\mathcal O_{a\tau}$ and of type {\rm D}. It is clear that
$X\# b\mu=R\# b\mu\cup S \# b\mu$ is a subrack decomposition of $ \mathcal O_{a\tau\# b\mu}$ and of type {\rm D}. \hfill $\Box$

\begin {Lemma}\label {2.7} Assume $a\pi \bot b\tau$. Let $\rho=\mu\otimes\lambda \in \widehat {\mathbb W_{m+n}^{a\pi \# b\tau}}$,  $ \mu \in \widehat{\mathbb W_{n}^{a\pi }} $,  $\lambda \in \widehat {\mathbb W_{m}^{b\tau}}$,  $a\pi\in \mathbb W_n$ and $b\tau\in \mathbb W_m$ with  $q_{a\pi, a\pi} {\rm id} = \mu (a\pi)$ and $q_{b\tau, b\tau } {\rm id}= \lambda (b\tau)$. If $\dim \mathfrak{B} ({\mathcal O}_{a \pi\# b \tau }^ {\mathbb W_{n+m} }, \rho)<\infty $,  then

{\rm (i)} $M( {\mathcal O}_{a \pi }^{\mathbb W_n},  \mu  )$ is isomorphic to a {\rm YD} submodule of  $M ({\mathcal O}_{a\pi\# b\tau }^  {\mathbb W_{n+m} },  \rho)$ over $\mathbb W_n$ when $q_{b\tau, b\tau}=1$; hence $\dim \mathfrak{B} ({\mathcal O}_{a \pi}^  {\mathbb W_{n} },  \mu) <\infty$.

{\rm (ii)} $q_{a\pi, a \pi} q_{b\tau, b\tau}=-1$.

{\rm (iii)} $q_{b\tau, b\tau}=1$ and $q_{a\pi, a\pi}=-1$ when ${\rm ord}(b\tau)\le 2$ and ${\rm ord}(q_{a\pi, a\pi})\not=1$.

{\rm (iv)} $q_{b\tau, b\tau} =1$ and $q_{a\pi, a\pi}=-1$ when ${\rm ord}(a\pi)$ and ${\rm ord }(b\tau)$ are coprime and ${\rm ord}(b\tau)$ is odd.
\end {Lemma}
\noindent {\it Proof.} {\rm (i)} The proof is similar to the one in \cite [Section 2.2] {AZ07}. Indeed,  Let $t_1 = \pi, t_2,   \cdots,
t_{m}$ be a numeration of ${\mathcal O}_\pi$,    and let $g_i\in G$ such that $g_i \rhd \pi :=
g_i \pi g_i^{-1} = t_i$ for all $1\le i \le m$. Let $s_1 = \tau, s_2,   \cdots,
s_{n}$ be a numeration of ${\mathcal O}_\tau$,    and let $h_i\in G$ such that $h_i \rhd \tau :=
h_i \tau  h_i^{-1} = s_i$ for all $1\le i \le n$. Thus  $\{ t_i \# s_j \mid 1 \le i \le m, 1\le  j \le n  \} = {\mathcal O}_{\pi \#\tau}$.  Let  $V$ and $W$ be  representation spaces of $\mu$ and $\rho$, respectively. Let  $0\not= w_0\in W$. Define a map $\psi $
from  $M({\mathcal O}_\pi^{{\mathbb W}_n },
\mu)$ to $M({\mathcal O}_ {\pi \# \tau}^{{\mathbb W}_m  \# {\mathbb W}_n},
\mu \otimes \rho)$  by sending $ g_i v$ to $(g_i \# h_1) v\otimes w_0$ for any $v\in V, $ $1\le i \le m.$  It is clear that $\psi$ is injective. Now we show that $\psi$ is a homomorphism of braided vector spaces. For any $v, v'\in V,$ we have that
\begin {eqnarray*}  (\psi \otimes \psi) ( C (g_iv  \otimes g_jv')) &=&
(\psi \otimes \psi )(  g_{i\rhd j} \mu (\nu _j(t_i))(v') \otimes g_i v  ) \\
&=& (g_{i\rhd j} \# h_1\mu (\nu  _j(t_i)(v') \otimes w_0) \otimes (g_i \# h_1 v \otimes w_0) \ \ \ \hbox {and }
\end {eqnarray*}
\begin {eqnarray*} C (\psi (g_iv) \otimes \psi (g_jv') &=&
C( (g_i \# h_1 v\otimes w_0) \otimes  (g_j \# h_1 v'\otimes w_0)   )\\
&=&  ( g_{i\rhd j} \# h_1( \mu \otimes \rho ) (  \nu _{j, 1} (t_i\# s_1) ) (v'\otimes w_0) ) \otimes (g_i \#h_1 v \otimes w_0).
\end {eqnarray*}
It is clear that $\nu _  {j, 1} (t_i \# h_1) = \nu _{ j} (t_i) \# s_1$. Therefore $\psi$
is a homomorphism of braided vector spaces.

{\rm (ii)} and  {\rm (iii)} are clear.

{\rm (iv)} Obviously,  ${\rm ord}(q_{a\pi, a\pi}) \mid {\rm ord}(a\pi)$  and ${\rm ord}(q_{b\tau, b\tau }) \mid {\rm ord}(b\tau)$. Therefore,  ${\rm ord}(q_{b\tau, b\tau})$ is odd.  The least  common factor $({\rm ord} (q_{a\pi, a\pi}),  {\rm ord}(q_{b\tau, b\tau}))=1$ since $({\rm ord}{(a\pi}), {\rm ord}({b\tau}))=1$. By Part {\rm (ii)},   ${\rm ord}(q_{a\pi, a\pi}){\rm ord}(q_{b\tau, b\tau})=2$. Consequently,  ${\rm ord}(q_{a\pi, a\pi})=2$ and ${\rm ord}(q_{b\tau, b\tau})=1$.

\hfill $\Box$

\section{Conjugacy classes of $\mathbb W_n$}\label {rack}
In this section we  prove that except in several  cases conjugacy classes of classical Weyl groups $\mathbb W_n$ are of type {\rm D}.

\begin {Lemma}\label {1.3} Let $p$ be odd with $p\ge 5$ and $a\tau \in \mathbb W_n$ with $\tau = (1\ 2\cdots \ p)$. Then $\mathcal{O}_{a\tau}^{\mathbb W_n}$ is of type {\rm D}.
\end {Lemma}
\noindent {\it Proof.} {\rm (i)} Assume that $a\tau$ is a negative cycle (defined in \cite [Appedix ]{ZZ12}) and $b=(1, 0, \cdots, 0, a_{p+1}, a_{p+2}, \cdots, a_n)$. Thus   $a \tau $ and $b \tau$ are conjugate. We assume $a =  (1, 1, \cdots, 1, $ $ a_{p+1}, a_{p+2}, \cdots, a_n) $ without lost generality. Obviously,  the right hand side of (\ref {e2.1.3}) is non-vanishing for $\mu=\tau ^2$,  i.e. ${\rm sq}(a\tau, b\tau^2) \not= b\tau^2.$ Let $R:=\mathbb Z_2^n \rtimes\tau \cap  \mathcal O _{a \tau}
^{\mathbb W _n} $   and $S:=\mathbb Z_2^n \rtimes \tau^2 \cap\mathcal O _{a \tau} ^{\mathbb W _n}$. It is clear that $R \cup
S$ is a subrack decomposition and is of type {\rm D}; notice that $\tau$ and $\tau ^2$ are conjugate in $\mathbb S_n.$

{\rm (ii)} Assume that $a\tau$ is a positive cycle. Let $b=(1, 0, 0, 1, 0, a_{p+1}, a_{p+2}, \cdots, a_n)$ when $p=5;$ $b=(1, 1, 0, \cdots, 0, a_{p+1}, a_{p+2}, \cdots, a_n)$ when $p>5.$ Thus $0\tau$,  $a \tau $ and $b \tau$ are conjugate. We assume $a=0$ without lost generality. It is clear that the right hand side of (\ref
{e2.1.3}) is not equal to $0$. Consequently,  $R \cup S$ is a subrack decomposition and is of type {\rm D} as Part {\rm (i)}.

\hfill $\Box$

\begin {Lemma}\label {1.3'} If $\sigma$  is of type ($3^2$),  then $\mathcal O_{a\sigma} ^{\mathbb W_6}$ is of type ${\rm D}$ for all $a \in \mathbb Z_2^6.$
\end {Lemma}
\noindent {\it Proof.} Let $\pi = (1 \ 2\ 3)$, $\xi = (4\ 5\ 6)$, $\tau =(1\ 2\ 3)(4\ 5\ 6)$ and  $\mu = (1\ 2\ 3)^2(4\ 5\ 6)$. We have that
the $4$-th,  $5$-th and $6$-th components of (\ref {e2.1.3}) are
\begin {eqnarray}\label {e1.3'}  (a_6+a_5,  a_4+a_6,  a_5+a_4) =  (b_6+b_5,  b_4+b_6,  b_5+b_4). \end {eqnarray}
By (\ref {e2.1.3}),
{\rm (i)} If $a=\alpha := (1, 1, \cdots, 1)$ and $b=(1, 0, 0, 1, 0, 0)$,  then (\ref {e1.3'}) does not hold.

{\rm (ii)} If $a=0$ and $b=(0, 0, 0, 1, 1, 0)$,  then (\ref {e1.3'}) does not hold.

{\rm (iii)} If $a=(1, 0, 0, 0, 0, 0)$ and $b=(1, 0, 0, 1, 1, 0)$,  then (\ref {e1.3'}) does not hold.

{\rm (iv) } If $a=(0, 0, 0, 1, 0, 0)$ and $b=(0, 0, 0, 0, 1, 0)$,  then (\ref {e1.3'}) does not hold.
Let $R := \mathbb Z_2^6 \rtimes \tau \cap \mathcal O_{a\sigma}^{\mathbb W_6}$ and $S := \mathbb Z_2^6\rtimes \mu\cap \mathcal O_{a\sigma}^{\mathbb W_6}$. It is clear that $R \cup S$ is a subrack decomposition of  $\mathcal O_{a\sigma}^{\mathbb W_6}$. Consequently $\mathcal O_{a\sigma }^{\mathbb W_6}$ is of type {\rm D}. \hfill $\Box$

\begin {Lemma}\label {1.3''} If $\sigma$ is of type ($2^2, 3^1$),  then $\mathcal O_{a\sigma}^{\mathbb W_7}$ is of type ${\rm D}$ for all $\mathbb Z_2^7.$
\end {Lemma}
\noindent {\it Proof.} Let $\pi=(5\ 6\ 7)$,  $\xi=(1\ 2)(3\ 4)$,   $\lambda=(1\ 3)(2\ 4)$,  $\tau = \pi \xi$ and $\mu=\pi \lambda$.
If $a=(a_1, a_2, a_3, a_4, 0, 0, 0)$,  let $b=(b_1, b_2, b_3, b_4, 1, 1, 0)$. If $a=(a_1, a_2, a_3, a_4, 1, 1, 1)$,  let $b=(b_1, b_2, b_3, b_4, 1, 0, 0)$. Then
the $5$-th,  $6$-th and $7$-th components of (\ref {e2.1.3}) are $(a_6+a_7, a_7+a_5, a_5+a_6)=(b_6+b_7, b_7+b_5, b_5+b_6)$, respectively. Consequently,  (\ref {e2.1.3}) does not hold.

Let $R := \mathbb Z_2^7 \rtimes \tau \cap \mathcal O_{a\sigma}^{\mathbb W_7}$ and $S := \mathbb Z_2^7 \rtimes \mu \cap \mathcal O_{a\sigma}^{\mathbb W_7}$. It is clear that $R \cup S$ is a subrack decomposition of  $\mathcal O_{a\sigma }^{\mathbb W_7}$. Consequently $\mathcal O_{a\sigma }^{\mathbb W_7}$ is of type $D.$ \hfill $\Box$

\begin {Example}\label {10.3}   If $\tau$ and $\mu$ are  of type $( 2^2)$, then $a\tau$ and $b\mu$ are  square commutative when they are conjugate to each other or $(-1) ^ {\sum _i^4 a_i} = (-1) ^ {\sum _i^4 b_i} $.

\end {Example}
\noindent {\it Proof.}    Let $\tau = (1\ 2) (3\ 4)$. Then  $\mu = (1\ 3) (2\ 4)$  or $\mu = \tau.$ It is clear that equation (\ref {e2.1.3}) becomes $a_1 + a_2 + a_3 +a_4
= b_1 +b_2 + b_3 + b_4.$  \hfill $\Box$

 \begin {Lemma}\label {2.6} {\rm (i)} Assume that $\tau$ and $\mu$ are conjugate  with ${\rm sq}(\tau, \mu) \not=\mu$ in $\mathbb S_n$ and $\tau(n) = \mu(n) =n$. If  $a \in \mathbb Z_2^n$ and there exists $i$ such that $a_i\not=a_n$ with $\tau(i)=\mu (i)= i$, then $\mathcal  O_{a\tau}$ is of type $D.$

{\rm (ii)} Assume that $n>4$ and $a\tau\in\mathbb W_n$ with type $(1^{n-2}, 2)$ of $\tau.$  If there exist $i, j$ such that $\tau(i)=i$ and $\tau(j)=j$ with $a_i\not=a_j$,  then $\mathcal O_{a\tau}$ is of type {\rm D}.

{\rm (iii)} Assume that $n>5$ and $a\tau\in\mathbb W_n$ with type $(1^{n-3}, 3)$ of $\tau.$ If there exist $i, j$ such that $\tau(i)=i$ and $\tau(j)=j$ with $a _i\not=a_j$,  then $\mathcal O_{a\tau}$ is of type {\rm D}.

\end {Lemma}
 \noindent {\it Proof.} {\rm  (i)}
 Let $R:= \{d\xi \in \mathcal O _{a \tau}^{\mathbb W _n}\mid\xi\in \mathbb S_{n-1}; d_n=0\}$ and $S:= \{d\xi\in\mathcal O_{a\tau} ^{\mathbb W_n}\mid \xi\in\mathbb S_{n-1};d_n=1\}$. Obviously,  $R\cup S$ is a subrack decomposition.
 Let $\xi \in \mathbb S_n$ with $\xi (i) = i$ and $\xi (n) = n$ such that $\xi \rhd \tau  = \mu.$ By simple computation we have $(i, n) \xi \rhd (a \tau ) = b \mu$ with $b_n = a_i.$ Consequently, $a\tau$ and $b\mu$ are in the same set of  $R$ and $S$, which implies that $R\cup S$ is of type ${\rm D}$.

{\rm (ii)} It is clear that ${\rm sq}(\tau, \mu)\not=\mu$ with  $\tau :=(1, 2)$ and $\mu :=(2, 3)$. Applying Part (i) we complete the proof.

{\rm (iii)} It is clear that ${\rm sq}(\tau, \mu)\not=\mu.$ with $\tau=(1, 2, 3)$ and $\mu=(2, 4, 3)$. Applying Part (i) we complete the proof.
 \hfill $\Box$

\begin {Lemma}\label {1.2} Let $\sigma=(a, \tau)\in \mathbb W_n$.  If  $\mathcal{O}_{\tau}^{\mathbb S_n}$ is of type {\rm D},  then so is $\mathcal{O}_{(a, \tau)}^{\mathbb W_n}$.
\end {Lemma}
\noindent {\it Proof.} Let $X=S\cup T$ be a subrack decomposition of $\mathcal{O}_{\tau}^{\mathbb S_n}$ and $s\in S, t\in T$ such that
\begin {eqnarray} \label {e1.1}s\rhd (t \rhd (s \rhd t)) \not= t.\end {eqnarray}
Let $h \rhd \tau=s$ and $g\rhd \tau=t$ with $h, g\in {\mathbb S_n}$. It is clear ${\rm sq }((h\cdot a,  s), (g\cdot a,  t))\not=(g \cdot a,  t) $
since (\ref{e1.1}); $(h\cdot a, s)=h\rhd(a, \tau )$ and  $(g\cdot a, t)=g\rhd(a, \tau)$.

$(<{\mathbb S_n}\cdot a>,  X)$ is a subrack,  where $<{\mathbb S_n}\cdot a>$  is the subgroup generated by subset ${\mathbb S_n}\cdot a$ of $\mathbb Z_2^n$. In fact,  for any $h,  g\in {\mathbb S_n},  \xi,  \mu \in X, $  we have
\begin {eqnarray*}(h\cdot a, \xi)\rhd (g\cdot a, \mu)=((h +\xi g+\xi\mu\xi^{-1}h) \cdot a, \xi\rhd\mu )\in (<{\mathbb S_n}\cdot a>,  X).\end {eqnarray*}
Thus $(<{\mathbb S_n}\cdot a>,  X)$ is a subrack. Consequently  $(<{\mathbb S_n}\cdot a>, X) \cap \mathcal{O}_{(a,  \tau)}^{\mathbb W_n}=(<{\mathbb S_n}\cdot a>, S)
\cap \mathcal{O}_{(a,  \tau)}^{\mathbb W_n}\cup (<{\mathbb S_n}\cdot a>,  T) \cap\mathcal{O}_{(a,  \tau)}^{\mathbb W_n}$ is a subrack decomposition of $\mathcal{O}_{(a,  \tau)}^{\mathbb W_n}$ and is of type {\rm D}. \hfill $\Box$

\begin{Theorem}\label{1.4}
Let $G = \mathbb W_n$ with $n>4$. Let $\tau\in\mathbb S_n$ be of type $(1^{\lambda_1}, 2^{\lambda_2},  \dots,
n^{\lambda_n})$ and $a\in \mathbb Z_2^n$ with $\sigma=(a, \tau)\in G$ and $\tau\not=1$. If $\mathcal{O}_{\sigma}^G$ is not of
type {\rm D},  then the type of $\tau$ belongs to one in the following list.
\renewcommand{\theenumi}{\roman{enumi}}   \renewcommand{\labelenumi}{(\theenumi)}
\begin{enumerate}
\item $(2,  3)$; $(2^3);$
\item   $(2^4);$ $(1,  2^2), (1^2, 3), (1^2, 2^2) ;$
\item     $(1^{n-2},  2)$ and  $(1^{n-3},  3)$ $(n >5)$ with $a_i = a_j$ when $\tau (i) = i $ and $\tau (j)=j$.
\end{enumerate}
\end{Theorem}
\noindent {\it Proof.} It follows from Lemma \ref {1.3''}, Lemma \ref {2.6}, Lemma \ref {1.2}, Lemma \ref {1.3}, Lemma \ref {1.3'} and \cite [Theorem 4.1] {AFGV08}. \hfill $\Box$

\section{Nichols algebras  of irreducible {\rm YD} module over $\mathbb W_n$}\label {nichols algebra}
In this section we show that  except in three  cases Nichols algebras of irreducible {\rm YD} modules over classical Weyl groups
$\mathbb W_n$ are infinite dimensional.

 Let ${\rm supp} M := \{g \in G \mid M_g \not= 0\}$ for $G$-comodule $(M, \delta )$, where $M_g := \{x \in M \mid \delta (x) = g \otimes x \}$.

We shall use the following facts:

\begin{Theorem} \label {}(\cite [Cor. 8.4] {HS08}) Let $n\in \mathbb N $, $n\ge 3$, and assume that $G=\mathbb{S}_n$ is the
  symmetric group.
  Let $U$ be a {\rm YD} module  over $G$. If $\mathfrak B (U)$ is  finite dimensional, then $U$ is an irreducible
  {\rm YD} module over $G$.
\end {Theorem}

\begin{Theorem}\label{th:sm-intro} ( \cite [Th. 1.1] {AFGV08})
    Let $m\ge 5$.  Let $\sigma\in \mathbb S_m$ be of type $(1^{n_1},2^{n_2},\dots,m^{n_m})$,
    let $  \mathcal O_\sigma $ be the conjugacy class of $\sigma$ and let $\rho=(\rho,V) \in
    \widehat{{S_m } ^{\sigma}}
    $. If $\dim \mathfrak B ( \mathcal O_\sigma, \rho) < \infty$, then
    the type of $\sigma$ and $\rho$ are in the following list:
    \renewcommand{\theenumi}{\roman{enumi}}\renewcommand{\labelenumi}{(\theenumi)}
    \begin{enumerate}
        \item\label{it:caso2}
            $(1^{n_1}, 2)$, $\rho_1 = {\rm sgn}$ or $\epsilon$, $\rho_2 ={\rm sgn}$.
        \item
            $(2, 3)$ in $\mathbb S_5$, $\rho_2 ={\rm sgn}$, $\rho_3= \overrightarrow{\chi_{0}}$.
        \item\label{it:caso222}
            $(2^3)$ in $\mathbb S_6$, $\rho_2=\overrightarrow{\chi_{1}}\otimes\epsilon$ or
            $\overrightarrow{\chi_{1}}\otimes {\rm sgn}$.
    \end{enumerate}
\end{Theorem}

\begin {Theorem}\label {1.5}

If $G$ is a finite group and  $\mathcal{O}_{\sigma}^G$ is of type {\rm D},  then {\rm dim} $\mathfrak B(\mathcal{O}_{\sigma}^G,  \rho) = \infty $ for any $\rho \in \widehat { G^\sigma}$.
\end {Theorem}
\noindent {\it Proof.} It follows from \cite {AFGV08} (or see the Appendix).
\hfill $\Box$

\begin{Theorem}\label{1.6} Assume $n>4$. Let $\tau\in\mathbb S_n$ be of type $(1^{\lambda_1},  2^{\lambda_2}, \dots, n^{\lambda_n})$ and $a\in\mathbb Z_2^n$ with $\sigma=(a, \tau)\in \mathbb W_n$ and $\tau\not=1$. If {\rm dim} $\mathfrak B(\mathcal{O}_{\sigma}^{\mathbb W_n}, \rho)<\infty $,  then the type of $\tau$ belongs to one in the following list:
\renewcommand{\theenumi}{\roman{enumi}}   \renewcommand{\labelenumi}{(\theenumi)}
\begin{enumerate}
\item $(2,  3)$; $(2^3);$
\item   $(2^4);$ $(1,  2^2), (1^2, 3), (1^2, 2^2) ;$
\item     $(1^{n-2},  2)$ and  $(1^{n-3},  3)$ $(n >5)$ with $a_i = a_j$ when $\tau (i) = i $ and $\tau (j)=j$.
\end{enumerate}
\end{Theorem}
\noindent {\it Proof.} It follows from Theorem \ref {1.4} and Theorem \ref {1.5}. \hfill $\Box$

\vskip.1in
 Let $a\mu =c\tau \# d\xi \in\mathbb W_n$ with $c\tau \in \mathbb W_m$,  $d\xi \in \mathbb
B_{n-m}$ and   $c\tau \bot d\xi$. Let
  $\rho=\rho _1\otimes \rho_2 \in  \widehat {
\mathbb W_n ^{a\mu}}=\widehat {\mathbb W_m ^{c\tau}} \times\widehat {\mathbb W_{n-m} ^{d\xi}}$.
If  $c=(1, 1, \cdots, 1)$, then    $\rho_1=(\chi _1\otimes\mu _1)\uparrow _{
G^{a\mu}_{\chi_1}} ^{G^{a\mu}}$ with $\chi_1 \in \widehat {(Z_2^ {m})^\tau} $,  $\mu _1\in \widehat{(\mathbb S_m^\tau)_{\chi_1}}$ (see \cite [Section 2.5] {ZZ12} and  \cite {Se77}).
Case $a=0$ and $a=(1, 1, \cdots, 1)$ were studied in
paper \cite [Theorem 1.1 and Table 1]{ZZ12}. Other cases are listed  as follows:

\begin {Corollary}\label {2.9} Under notation above assume  $\dim \mathfrak B( {\mathcal O} _{a\mu} ^{\mathbb W_n},  \rho)< \infty$ with  $d_1=d_2=\cdots=d_{n-m}$ and   $\xi ={\rm id}$.

{\rm (i)} Then $\rho _1 (c\tau)=-id $ when $\rho _2(d\xi)={\rm id} $ and $\rho _1(c\tau)={\rm id}$ when $\rho _2(d\xi)=-{\rm id} $.

{\rm (ii)} Case $\tau=(1\ 2)$,  $c=0$ and $d=(1, 1, \cdots, 1)$. Then $\rho_1(c\tau)=\pm {\rm id}$,  $\rho_2(d\xi)=\mp{\rm id} $,  $\chi_1(c)=1$ and $\mu _1 (\tau)=\pm 1.$

{\rm (iii)} Case $\tau=(1\ 2)$,  $c=(1, 1)$ and $d =0$. Then  $\rho_1 (c\tau) =-{\rm id}$,  $\rho _2 (d\xi)={\rm id}$.

{\rm (iv) } Case $\tau=(1\ 2\ 3)$,  $c=0$ and  $d =(1, 1, \cdots, 1)$. Then $\rho_1(c\tau)={\rm id}$,  $\rho_2(d\xi)=-id$,  $\chi_1(c)=1$ and $\mu _1 (\tau)= 1.$

{\rm (v)} Case $\tau=(1\ 2\ 3)$,  $c=(1, 1, 1)$ and $d =0$. Then $\rho_1(c\tau)=-id$,  $\rho_2(d\xi)={\rm id}$,  $\chi_1(c)=-1$ and $\mu _1 (\tau)=1.$

{\rm (vi)} Case $\tau=(1\ 2)( 3\ 4)$,  $c=0$ and $d=(1, 1)$. Then $\rho_1(c\tau)={\rm id}$,  $\rho_2(d\xi)=-id$.

{\rm (vii)} Case $\tau=(1\ 2)( 3\ 4)$,  $c=(1, 0, 1, 0)$ and $d=(1, 1)$. Then $\rho_1(c\tau)=\pm {\rm id}$,  $\rho_2(d\xi)=\mp {\rm id}$.

{\rm (viii)} Case $\tau=(1\ 2)( 3\ 4)$,  $c=(1, 0, 1, 0)$ and $d=0$. Then $\rho_1(c\tau)=-id$,  $\rho _2(d\xi)={\rm id}$.

{\rm (ix)} Case $\tau=(1\ 2)( 3\ 4)$,  $c=(1, 0, 0, 0)$ and $d=(1, 1)$. Then $\rho_1(c\tau)=\pm {\rm id} $,  $\rho_2(d\xi)=\mp {\rm id}$.

{\rm (x)} Case $\tau =(1\ 2)( 3\ 4)$,  $c=(1, 0, 0, 0)$ and $d=(0, 0)$. Then $\rho_1(c\tau)=-id$,  $\rho_2(d\xi)={\rm id}$.
\end {Corollary}
\noindent {\it Proof.} {\rm (iv)} If $\rho_2(d\xi)={\rm id}$,  then $\dim \mathfrak{B}({\mathcal O}_{c\tau}^{\mathbb W_{3}}, \rho_1)<\infty$ by Lemma \ref {2.7},  which constracts to \cite [Theorem 1.1] {ZZ12}.

{\rm (v)} If $\chi _1(c) =1$,  then there exists a contradiction by \cite [Proposition 2.4,  Theorem 1.1] {ZZ12}.

The others  follow from Lemma \ref {2.7} and \cite [Theorem 1.1]{ZZ12}. \hfill $\Box$

\section {Relationship between Nichols algebras and FK algebras} \label {fk conjecture}
In this section we give the relationship between Nichols algebra $\mathfrak{B}({\mathcal O}_{\sigma}, \rho)$ and {\rm FK} algebra ${\mathcal E}_n$ defined in \cite [Definition 2.1]{FK99},  where $\sigma$ is a transposition in $\mathbb S_n$, and  $\rho= {\rm sgn} \otimes {\rm sgn}$ or $\rho= \epsilon \otimes {\rm sgn}$.
We generalize {\rm FK } algebra.

$\mathfrak{B}({\mathcal O}_{\sigma}, \rho)$ is finite dimensional when $n\le 5$ according to \cite {MS00,  FK99,  AZ07}. However,  it has been an open problem whether or not $\mathfrak{B}({\mathcal O}_{\sigma}, \rho)$ is finite dimensional when $n>5.$

Let $\sigma=(12)\in S_{n}$,  $\mathcal O_{\sigma}=\{(ij)| 1 \le i,  j \le n\}$,
$G^{\sigma}=\{g\in G\mid g\sigma=\sigma g\}= \mathbb S_{ \{3, 4, \cdots, n\}} \times \mathbb S_{\{1, 2\}}$£¬
 \begin {eqnarray}\label {pe4.5.1}
G=\bigcup\limits_{1\le i < j\le n}G^{\sigma}g_{ij}. \end  {eqnarray}
Notice that $(ij) = (ji)$ since $(ij)$ is a transposition.
Let
$g_{k j}  := \left \{\begin
{array} {lll} {\rm id} \ \ \ \ \  \ \ \ \ \ \ \ \ \ \ \ \ k=1,  j=2\\
(2j) \ \  \ \ \ \ \ \ \ \ \ \ \ \  k=1,  j>2\\
(1j)\ \  \ \ \ \ \ \ \ \ \ \ \ \  k=2,  j>2 \\
(1k)(2j)\ \ \ \ \ \ \ \ \ k>2,  j>k \\
\end {array} \right.$ and $t_{ij} = (ij).$

Let $a_{ij}$ denote the arrow $a_{t_{ij}, 1}$ from $1$ to $t_{ij}$.
By \cite [Lemma 1.1]{ZZWCY09}  or appendix,
$\{a_{ij}\mid i\not=j, 1\le i, j\le n\}$ generates an algebra $\mathfrak{B}({\mathcal O}_{(1 2)}, {\rm ad } (\chi))$, which is isomorphic to Nichols algebra $\mathfrak{B}({\mathcal O}_\sigma, \rho),$ in copath Hopf algebra $kQ^c$.

\begin {Lemma} \label {6.2} In bi-one Nichols algebra  $\mathfrak{B}({\mathcal O}_{(1 2)}, {\rm ad } (\chi))$ with $\chi={\rm sgn} \otimes {\rm sgn}$ or $\chi= \epsilon \otimes {\rm sgn}$,

{\rm (i)} If $i, j $ and $k$ are different,  then there exist $\alpha_{i, j, k}$,  $\beta_{i, j, k}\in\{1, -1\}$ such that
\begin {eqnarray}\label {e6.1.1} a_{ij}a_{jk}+\alpha_{ijk}a_{jk}a_{ki}+\beta_{ijk}a_{ki}a_{ij}=0. \end  {eqnarray}

{\rm (ii)} \begin {eqnarray}\label {e6.1.3}  \hbox { Left hand side of (\ref {e6.1.1})} &=&
(\chi(\zeta_{ij}(t_{jk}))\underline  { a_{t_{ij}t_{jk}, t_{jk}}a_{t_{jk}, 1}} +\underline {a_{t_{ij}t_{jk}, t_{ij}}a_{t_{ij}, 1}}) \nonumber \\
  & &+ \alpha _{ijk}(\chi (\zeta_{jk}(t_{ik}))\underline {a_{t_{jk}t_{ik}, t_{ik}}a_{t_{ik}, 1}} +\underline {a_{t_{jk}t_{ik}, t_{jk}}a_{t_{jk}, 1}}) \nonumber \\
  && + \beta_{ijk}(\chi(\zeta_{ik}(t_{ij}))\underline {a_{t_{ik}t_{ij}, t_{ij}}a_{t_{ij}, 1}} +\underline {a_{t_{ik}t_{ij}, t_{ik}}a_{t_{ik}, 1} } ) \end {eqnarray}

{\rm (iii)} If $\chi(\zeta_{ij}(t_{jk}))\chi(\zeta_{jk}(t_{ik}))\chi(\zeta_{ik}(t_{ij}))=-1, $ then Part (i) holds.

{\rm (iv)} If $i, j$ and $k$ are different,  then  Part (i) holds if and only if $$\chi(\zeta_{ij}(t_{jk}))\chi(\zeta_{jk}(t_{ik}))\chi(\zeta_{ik}(t_{ij}))=-1.$$

{\rm (v)} If $i, j, k$ and $l$ are different,  then there exist $\lambda_{i j kl} \in\{1, -1\}$
 such that $a_{ij}a_{kl}=\lambda_{ijkl}a_{kl}a_{ij}.$
\end {Lemma}
\noindent {\it Proof.} Let $\chi' :={\rm sgn}\otimes {\rm sgn}$ and $\chi '':=\epsilon \otimes {\rm sgn}.$ It is clear that
$M({\mathcal O}_{(12)}, {\rm ad }(\chi))$ is a PM $\mathbb S_n$-{\rm YD} module (see \cite [Definition 1.1] {ZZC07}) and $g\cdot a_{ij}=a_{g t_{ij},  g}$,  $a_{{ij}} \cdot g=\chi(\zeta_{ij}(g))a_{t_{ij}g, g}$ (see \cite [Lemma 1.9] {ZZC07}). By \cite {CR02},
\begin {eqnarray}\label {e6.1.3'}a_{{ij}}a_{{kl}}=\chi(\zeta_{ij}({t_{kl}}))\underline {a_{{t_{ij}}{t_{kl}}, {t_{kl}}}a_{{t_{kl}}, 1}}+ \underline {a_{{t_{ij}}{t_{kl}}, {t_{ij}}}a_{{t_{ij}}, 1}}.
\end {eqnarray}

{\rm (ii)} It  follows from (\ref {e6.1.3'}). In fact, $\alpha _{i jk} = - \chi (\zeta _{ij} (t_{jk})) $ and $\beta _{i jk} = - \frac {1}{ \chi (\zeta _{ik} (t_{ij})) }.$

{\rm (iii)} and {\rm (iv)}  follow from Part {\rm (ii)}.

{\rm (i)} Let $a, b$ and $c$ stand for $\zeta_{ij}(t_{jk})$,  $\zeta_{jk}(t_{ik})$ and  $\zeta_{ik}(t_{ij})$, respectively, in Table $1$ below.
\vskip.1in
\begin{tabular}{|l|l|l|l|l|l|l|l|l|l|}
\hline { \bf case } & $a$ & $b$& $c$ &$\chi' (a) $ & $\chi' (b)$ &$\chi' (c) $ &$\chi ''(a) $ & $\chi ''(b)$ &$\chi'' (c) $  \\
\hline $ 2< i<j<k $ &$ (kj)$  &$ (12)(ijk)$ &$  (ij)$&$  -1$ &$  -1 $&$  -1$ &$ 1$ &$  -1 $&$  1$ \\
\hline $ i=1,  j=2 <k $ &$ (1)$  &$ (1)$ &$  (12)$&$  1$ &$  1 $&$  -1$ &$ 1$ &$ 1 $&$  -1$ \\
\hline $ i=1,  2< j<k $ &$ (kj)$  &$ (12)(kj)$ &$  (1)$&$  -1$ &$  1 $&$  1$ &$1$ &$  -1 $&$  1$ \\
\hline $ i=2<j<k $ &$ (kj)$  &$ (1)$ &$ (12) (jk)$&$  -1$ &$  1 $&$  1$&$ 1$ &$  1 $&$  -1$ \\
\hline $ 2< i<k<j $ &$ (kj)$  &$ (ik)$ &$(12)(ijk)$&$  -1$ &$  -1 $&$  -1$ &$ 1$ &$  1 $&$  -1$ \\
\hline $ i=1,  k=2<j $ &$ (1)$  &$ (12)$ &$  (1)$&$  1$ &$  -1 $&$  1$ &$ 1$ &$  -1 $&$  1$ \\
\hline $ i=1,  2< k<j $ &$ (kj)$  &$ (1)$ &$  (12) (jk)$&$  -1$ &$  1 $&$  1$ &$ 1$ &$  1 $&$ - 1$ \\
\hline $ i=2< k<j $ &$ (kj)$  &$ (kj) (12)$ &$  (1)$&$ - 1$ &$  1 $&$  1$ &$ 1$ &$  -1 $&$  1$ \\ \hline
\end{tabular} $$\hbox {Table } 1$$
By Table $1$,  $\chi (\zeta _{ij} (t_{jk}))\chi (\zeta _{jk} (t_{ik}))\chi (\zeta _{ik} (t_{ij})) =-1$.
Consequently, Part (i) holds by Part {\rm (ii-iv)}.

(v) By (\ref {e6.1.3'}),  $a_{ij} a_{kl}=\lambda_{ijkl}a_{kl}a_{ij} $ if and only if $(\chi(\zeta_{ij}(t_{kl}))-\lambda_{ijkl})=(1- \chi(\zeta_{kl}(t_{ij}))\lambda_{ijkl})=0$.  Since $(ij)(kl)=(kl)(ij)$,  we have $g_{ij}(kl)=g_{ij}(kl)g_{ij}g_{ij}$ with $g_{ij} (kl)g_{ij}\in\mathbb S_n^{(12)}$ and $\zeta_{ij}((kl))=g_{ij}(kl)g_{ij}.$  See $\zeta_{ij}(t_{kl})=g_{ij}t_{kl}g_{ij}=$\\
$\left \{\begin
{array} {lll l} (kl) & \hbox  { if } k,  l >2\\
(kl)  & \hbox { if }  i=1,  j = 2\\
(2j) (kl)(2j) = (lj) \hbox  { or } (kj) \hbox  { or } (kl)  & \hbox  { if }   i=1,  j>2,  k=2 \hbox  { or } i=1,  j>2,   l=2  \\ & \hbox  { or }  i=1,  j>2,  k\not=2,  l\not= 2\\
(1j) (kl)(1j) = (lj) \hbox  { or } (kj) \hbox  { or } (kl)  & \hbox  { if }   i=2,  j>2,  k=1 \hbox  { or } i=1,  j>2,   l=1  \\ & \hbox  { or }  i=1,  j>2,  k\not=1,  l\not= 1\\
 \hbox { a transposition of } \\  \hbox { two numbers greater than } $2$ & \hbox  { if }   i,  j>2
\end {array} \right. .$

\noindent Thus,  $\zeta_{ij}(t_{kl})$ is a transposition of two numbers greater than $2$ and
$\chi'(\zeta _{ij}(t_{kl}))=-1$ and $ \chi'' (\zeta _{ij}(t_{kl}))=1$. Similarly,  $\chi'(\zeta_{kl}(t_{ij}))=-1$ and $\chi''(\zeta _{kl}(t_{ij}))=1$. Consequently,  it is enough to set $\lambda_{ijkl} =: \chi(\zeta_{ij}(t_{kl}))$. \hfill $\Box$

\vskip.1in
Obviously,  $\mathcal E_n={\mathcal E}_n(1, 1, -1, 1), $ i.e. $\mathcal E_n= {\mathcal E}_n(\alpha, \beta, \gamma, \lambda )$ with $\alpha _{i,  j,  k} = 1, \beta _{i,  j,  k} =1, \gamma _{ij} =-1, \lambda _{i,  j,  k,  l} =1$ for any distinct $i, j, k$ and $l.$ ${\mathcal E}_n(\alpha, \beta, \gamma, \lambda )$ is called a generalized {\rm FK} algebra.

\begin {Definition} \label {4.4'}  (\cite [Definition 2.1]{FK99}) {\rm FK} algebra ${\mathcal E}_n$ is generated by $\{x_{ij} \mid 1\le i< j \le n\}$ with defining relations:

{\rm (i)} $x_{ij}^2=0$ for $i<j;$

{\rm (ii)}  $ x_{ij}x_{jk}=x_{jk}x_{i k}+x_{i k}x_{ij}$ and $x_{jk}x_{ij}=x_{i k}x_{jk}+x_{ij}x_{i k}, $ for $i<j<k;$

{\rm (iii)} $x_{ij}x_{kl}=x_{kl}x_{ij}$
  for any distinct $i, j, k$ and $l, $ $i<j, k<l.$

\noindent
Equivalently,  {\rm FK} algebra ${\mathcal E}_n$ is generated by $\{x_{ij}\mid i\not=j, 1\le i, j \le n\}$ with defining relations:

{\rm (i)} $x_{ij}^2=0$,  $x_{i j}=-x_{ji}$;

{\rm (ii)} $ x _{ij}x_{jk}+x_{jk}x_{ki}+x_{ki}x_{ij}=0;$

{\rm (iii)} $x_{ij}x_{kl}=x_{kl}x_{ij}$ for any distinct $i, j, k$ and $l.$
\end {Definition}

By \cite [Theorem  7.1]{FK99},  a subring of $\mathcal E_n$ is isomorphic to the cohomology ring of the flag manifold. It has been conjectured by
Fomin and Kirillov \cite [Conjecture 2.2]{FK99} that $ \dim {\mathcal E}_n<\infty$.
Consequently,    we have

\begin {Conjecture} \label {4.4} Let $\alpha_{ijk}, \beta_{ijk}, \gamma_{ij}, \lambda_{ijkl}\in\{1, -1\}$ for any distinct $i, j, k$ and $l$ with $1\le i,  j, k,  l \le n$. Assume that algebra ${\mathcal E}_n(\alpha, \beta, \gamma, \lambda )$ is generated by $\{x_{ij}\mid i\not=j, 1\le i, j\le n\}$ with defining relations:

  {\rm  (i)} $x_{ij}^2=0$,  $x_{i j}=\gamma_{ij}x_{ji}$;

 {\rm  (ii)} $x_{ij}x_{jk}+\alpha_{ijk}x_{jk}x_{ki}+\beta_{ijk}x_{ki}x_{ij}=0;$

  {\rm (iii)} $x_{ij}x_{kl}=\lambda_{ijkl}x_{kl}x_{ij}$ for any distinct $i, j, k$ and $l.$
\noindent
Then ${\mathcal E}_n(\alpha, \beta, \gamma, \lambda )$ is finite dimensional when $n>5$.

\end {Conjecture}

\begin {Theorem} \label {5.8} Assume  $n>3$.

{\rm (i)} $\mathfrak{B} ({\mathcal O}_{{(1, 2)}} ,  \epsilon \otimes {\rm sgn} )$ is an image of $\mathcal E_n$.

{\rm (ii)} If  $\dim \mathfrak{B} ({\mathcal O}_{{(1, 2)}} ,  \epsilon\otimes {\rm sgn}) = \infty$,  then so is  $\dim \mathcal E_n$.
\end {Theorem}
\noindent {\it Proof.}
Let $b_{ij} = -a_{ij} $ when $ i= 2$  and $j>2$;  $b_{ij} =  a_{ij}$ otherwise.  It follows from  (\ref {e6.1.3}) that  $$
  \begin{tabular}{|l|l|l|l|l|l|l|l|l|l|}
  \hline
  { \bf case } & $\alpha _{ijk}$ & $\beta _{ijk}$ \\
\hline
$ 2< i<j<k $ &$-1$ & $-1$
 \\
\hline
$ i=1,  j=2 <k $ &$-1$ & $1$
 \\
\hline
$ i=1,  2< j<k $ &$-1$ & $-1$
 \\
\hline
$ i=2<j<k $ &$-1$ & $1$
 \\
\hline
$ 2< i<k<j $ &$-1$ & $1$
 \\
\hline
$ i=1,  k=2<j $ &$-1$ & $-1$
\\
\hline
$ i=1,  2< k<j $ &$-1$ & $1$
\\
\hline
$ i=2< k<j $ &$-1$ & $-1$
\\
\hline
\end{tabular} \ \ \ . $$
 $$\hbox {Table } 2$$
\noindent Using Table 2 we can show that  $b_{ij}^2=0$ for $i<j$,
$ b_{ij} b_{jk} = b _{jk} b_{i k} + b _{i k} b_{ij}  $ and  $b_{jk} b _{ij} = b_{i k} b _{jk} + b_{ij}b _{i k} ,  $ for $i<j<k.$ By  Lemma  \ref {6.2}(iv),  $b_{ij} b_{kl} =b_{kl} b_{ij}$   for any distinct $i,  j,  k$ and $l, $  $i<j,  k<l.$
Consequently,   Part (i) holds since $\{b_{ij} \mid  1\le i <j \le n\}$ generates $\mathfrak{B} ({\mathcal O}_{{(1, 2)}},  \epsilon\otimes {\rm sgn}) $.  \hfill $\Box$

\begin{Remark}    ${\mathcal E}_5$ is finite dimensional  (see \cite {FK99} and
 \cite [Section 3.4]{AS02}).
Consequently,  $\dim \mathfrak{B} ({\mathcal O}_{{(1, 2)}} ,  \epsilon\otimes {\rm sgn}) < \infty$ when $n=5.$
\end{Remark}

Let $\chi ' := {\rm sgn} \otimes {\rm sgn}$ and $\chi '' := \epsilon \otimes {\rm {\rm sgn}}$. Let $\phi _1 $ and $\phi_2$ be maps from
$\mathbb S_n\times T$ to ${\bf k}\setminus 0$ such that  $\phi_1 (g,  t) :=  \left \{ \begin {array} {ll} 1 &  \hbox {if } g(i) < g(j) \\
 -1 &  \hbox {if } g(j) < g(i) \\ \end {array} \right. $ and $\phi_2 (g,  t) := (-1)^{l(g)}$ where $t = (i, j)$,   $T := \{(u,  v) \mid 1\le u,  v \le n,  u \not= v \}\subseteq \mathbb S_n$  and  $l(g)$ is the length of $g$,  that is,  $l(g)$  is the minimal number $q$ such that $g = g_1g_2\cdots g_q$ with $g_i \in T$ ,  $1\le i\le q.$  Let $\{x_{ij} \mid (i,  j) \in T\}$ be a basis of $M(\mathbb S_n,  T,  \phi  )$ with $\phi = \phi_1$ or $\phi =\phi_2.$  Define module and comodule operations as follows: $g \rhd x_{ij} = \phi (g,  (i,  j)) x_{g \rhd (i,  j)}$,  $\delta ^- (x_{ij}) = (i,  j) \otimes  x_{ij}, $ for any $g\in \mathbb S_n,  (i,  j) \in T.$ By \cite [Definition 5.1 and Example 5.3]{MS00},  $M(\mathbb S_n,  T,  \phi  )$ is a ${\rm YD}$ module over $\mathbb S_n.$ Its Nichols algebra is written as $\mathfrak{B} (\mathbb S_n,  T,  \phi)$. By \cite [Definition 2.1] {FK99} and \cite [Example 6.2 ] {MS00},  $\mathfrak{B} (\mathbb S_n,  T,  \phi)$ is an image of $\mathcal E_n$.

 \begin {Lemma} \label {5.9} If  $n>4$,  then $M(\mathbb S_n,  T,  \phi_1  )$  is not isomorphic to $M ({\mathcal O}_{{(1, 2)}} ,  \chi' )$  as {\rm YD} modules over $\mathbb S_n.$    $M(\mathbb S_m,  T,  \phi_2  )$ is not isomorphic to  $M ({\mathcal O}_{{(1, 2)}} ,  \chi'' )$ as {\rm YD} modules over $\mathbb S_n$.

 \end {Lemma}
\noindent {\it Proof.}
 If there exists isomorphism $\psi : M(\mathbb S_n,  T,  \phi  ) \rightarrow M ({\mathcal O}_{{(1, 2)}} ,  \chi )$ as  {\rm YD} mosules over $\mathbb S_n$,  where $\phi =\phi_1$ or $\phi =\phi_2$,  $\chi =\chi'$ or $\chi =\chi''$,  then  there exists $k_{ij} \in {\bf  k}$ such that $\psi (x _{ij}) = k_{ij} a_{ij}$ for any $(i,  j) \in T$,  where $a_{ij}$ denotes arrow $a_{t_{ij},  1}$ in short,  since $\psi$ is a comodule isomorphism.  By
$\psi (g \rhd x_{ij}) = g \rhd \psi (x_{ij})$,  we have
 \begin {eqnarray}\label {e4.5.1} k_{g\rhd (ij) } \phi (g,  (i, j)) a _{ g\rhd (ij)}  = k_{ij} \chi (\zeta _{ij} (g^{-1})) a _{g\rhd (ij)},  \end {eqnarray} for any $g\in \mathbb S_n$,  $(i,  j)\in T.$   For convenience,  set $k_{ij} =k_{ji}.$

  Let $ 2<i<j $,  $g = (1,  2)$. It is clear $g_{ij} g = (1i)(2j) (12) = (ij)(1i)(2j)$. We have $g = g^{-1}$ and  $\zeta _{ij} (g) =(ij)$. Thus $\chi ' (\zeta _{ij} (g)) =-1$ and  $\chi '' (\zeta _{ij} (g)) =1$; $\phi _1 (g,  (ij)) = 1$ and $\phi _2 (g,  (ij)) = -1$. Considering (\ref {e4.5.1}) we have
 $M(\mathbb S_n,  T,  \phi_1  )$  is not isomorphic to $M ({\mathcal O}_{{(1, 2)}} ,  \chi' )$  as {\rm YD} modules over $\mathbb S_n.$    $M(\mathbb S_m,  T,  \phi_2  )$ is not isomorphic to  $M ({\mathcal O}_{{(1, 2)}} ,  \chi'' )$ as {\rm YD} modules over $\mathbb S_n$  since $k_{g\rhd (ij) } \phi_1 (g,  (i, j))   \not= k_{ij} \chi '(\zeta _{ij} (g)) $ and $k_{g\rhd (ij) } \phi _2(g,  (i, j))  \not= k_{ij} \chi'' (\zeta _{ij} (g))$. \hfill $\Box$

\begin {Proposition} \label {5.10}  If there exists a natural  number $n_0>5$ such that
$M(\mathbb S_{n_0},  T,  \phi_1  )$  is not isomorphic to $M ({\mathcal O}_{{(1, 2)}} ,  \chi'' )$  as {\rm YD} modules over $\mathbb S_{n_0}$,  then
 $ \dim \mathfrak{B} (\mathbb S_n,  T,  \phi_1) = \infty$  and $\dim \mathcal E_n = \infty$ for any $n\ge n_0.$.
\end {Proposition}
\noindent {\it Proof.}    $M(\mathbb S_6,  T,  \phi_1 )$  is not isomorphic to $ M ({\mathcal O}_{{(1, 2)(3,  4)( 5,  6)}} ,   \rho)$ as  {\rm YD} modules over $\mathbb S_6$ since they  are not isomorphic as comodules over $\mathbb S_6$,  when $\rho$ is one dimensional representation.

 If  $\dim \mathfrak{B} (\mathbb S_n,  T,  \phi_1) < \infty $,  then $\dim \mathfrak{B} (\mathbb S_{n_0},  T,  \phi_1) < \infty $ since $\mathbb S_{n_0}$  is a subgroup of $\mathbb S_n.$
   $M(\mathbb S_{n_0},  T,  \phi_1 )$  is a reducible   {\rm YD} modules over $\mathbb S_{n_0}$ by Lemma \ref {5.9} and \cite [Theorem 1.1] {AFGV08}. However,  every reducible     {\rm YD} modules over $\mathbb S_{n_0}$ is infinite dimensional by \cite [Corollary 8.4] {HS08}. This is a contradiction.  Consequently,  $\dim \mathfrak{B} (\mathbb S_n,  T,  \phi_1) =\infty .$    By \cite [Definition 2.1] {FK99} and \cite [Example 6.2 ] {MS00},  $\mathfrak{B} (\mathbb S_n,  T,  \phi_1)$ is an image of $\mathcal E_n$. Therefore,   $\dim \mathcal E_n = \infty$ for any $n\ge n_0.$ \hfill $\Box$

\vskip.1in
Recall \cite  [Section 4.2] {Ba06}.
Let  $W$ be   Weyl group of a simple Lie algebra and $\Phi$ be the root system of $W$.
Let $\alpha$ be a root, and $s_\alpha\in W$ be the corresponding reflection.

Let $V_{W}$ be the linear space spanned by symbols
$\widetilde{\alpha}$ where $\alpha$ is a root of $W$, subject to the relation
$\widetilde{-\alpha}=- \widetilde{\alpha}$.
The dimension of $V_{W}$ is thus $|\Phi^+|$.

The $W$-action on $V_{W}$ is given by
$w \widetilde{\alpha}= \widetilde{ w\alpha}$, and the $W$-grading
is given by assigning the degree $s_\alpha$ to the basis element $\widetilde {\alpha }$.
The action and the grading are compatible, so
that $V_{W}$ is a Yetter Drinfeld module over $W$.

By \cite [Section 7] {Ba06}, $\mathcal E_n = T(V_{\mathbb S_n}) /I$, where $I$ is an ideal generated by ${\rm ker} ({\rm id } + C) = {\rm ker} ({\rm id} + C^{-1}) = {\rm ker} S_2$.  Obviously, $I \subseteq S := \bigoplus \limits _{m=2}^{\infty} {\rm ker} S_{m}$. Consequently,
the Nichols algebra
$$\mathcal B _{\mathbb S_n} := \mathfrak B(V _{\mathbb S_n}) = T(V _{\mathbb S_n}) / S \cong  (T(V _{\mathbb S_n})/I) / (S/I)$$ is a  quotient of  $\mathcal E_n$ as algebras.
By \cite [Th. 6.1]{Ma10},  $\mathcal E_n$ is a braided Hopf algebra in  category of {\rm YD} modules over $\mathbb S_n.$

We have the following as Proposition \ref {5.10}.
\begin {Proposition} \label {5.10¡®}  If there exists a natural  number $n_0>5$ such that
$  \mathfrak B(V _{\mathbb S_{n_0}}) $  is not isomorphic to $M ({\mathcal O}_{{(1, 2)}} ,  \chi'' )$  and $M ({\mathcal O}_{{(1, 2)}} ,  \chi' )$ as {\rm YD} modules over $\mathbb S_{n_0}$,  then
 $ \dim   \mathfrak B(V _{\mathbb S_n})   = \infty$  and $\dim \mathcal E_n = \infty$ for any $n \ge  n_0.$.
\end {Proposition}

\section {PM Nichols algebras and {\rm FK} algebra $\mathcal E_n$}  \label {bi-one arrow}

In this section we give an estimate for the dimensions of the {\rm PM} Nichols  algebras and {\rm FK} algebra $\mathcal E_n$.

 \begin {Definition}\label {ppppp1}   If $(X, \rhd)$ is a rack and  $\{v_\alpha \mid \alpha \in X\}$ is a basis of braided vector space $(V, C)$ such that  $C(v_\alpha \otimes v_\beta) = q_{\alpha, \beta} v_{\alpha \rhd \beta} \otimes v_\alpha$ and $0\not= q_{\alpha, \beta } \in \mathbb C$ for any $\alpha, \beta \in X,$ then $(V, C)$ is called a braided vector space of  rack-diagonal type; $(q_{\alpha, \beta})$ is called the braiding matrix and   $\{v_\alpha \mid \alpha \in X\}$ is  called a  canonical basis.
  \end {Definition}

\begin {Proposition}\label {ppppp2}
 $(V, C)$ is  a braided vector space of  rack-diagonal type with braiding matrix $(q_{\alpha, \beta})_{\alpha, \beta \in X}$ if and only if
\begin {eqnarray}\label {pppppe1} q_{a, b\rhd c}q_{b,  c}q_{a, b} =q_{a, b}q_{a,  c}q_{ a \rhd b, a\rhd c}
\end {eqnarray} for any $a, b, c\in X.$
  \end {Proposition}

  \noindent {\it Proof.}  It is clear that {\rm YBE}:
  $(C \otimes {\rm id }) ({\rm id} \otimes C)(C \otimes {\rm id }) = ({\rm id} \otimes C)(C \otimes {\rm id }) ({\rm id} \otimes C)$ holds if and only if (\ref {pppppe1}) holds.

Necessity is clear. Sufficiency. $C^{-1}(v_\alpha \otimes v_\beta) = q_{ \beta, \beta \rhd ^{-1} \alpha} ^{-1} v_{ \beta} \otimes v_{\beta \rhd ^{-1}\alpha}$. Here $\phi _\alpha (\beta) := \alpha \rhd \beta$ and $\alpha \rhd  ^{-1}\beta:=\phi _\alpha ^{-1}(\beta)$. \hfill $\Box$

\begin {Proposition}\label {ppppp3}
Assume $\mid Q_1^1\mid < \infty$ for  ${\rm RSC }(G,  r,  \overrightarrow \chi,  u)$. Then
$(kQ_1^{1}, ad( G, r, \overrightarrow{\chi }, u)) $  is a braided vector space with
of  rack-diagonal type with rack $X = \cup _{C\in {\mathcal K}_r(G)} C$  and   canonical basis $Q_1^1 := \{ a_{\alpha, 1}  \mid  a_{\alpha, 1} \hbox  { is an arrow from } 1  \hbox  { to } \alpha  \hbox { with } \alpha \in X\}$.  \end {Proposition}

  \noindent {\it Proof.}  For any $\alpha, \beta \in X$, by \cite [ Pro. 1.9] {ZZC07},
  $\alpha  \rhd a _{\beta, 1} := q_{\alpha, \beta}a_{\alpha \rhd \beta, 1}$ and   $ C(a _{\alpha , 1}  \otimes a _{\beta, 1} ) = q_{\alpha, \beta} a _{\alpha \rhd \beta, 1}  \otimes a _{\alpha, 1}$.
     \hfill $\Box$

Let $\mathfrak{B} (G,
r, \overrightarrow{\chi }, u)$ denote $\mathfrak{B}(kQ_1^{1}, ad( G, r, \overrightarrow{\chi }, u))$. Set $X = \{ \alpha _1, \cdots, \alpha_n\}$. Assume $\{y_1,  y_2,  \cdots,  y_n\}$ be the dual basis of $\{ a_{\alpha_1,1},  a_{\alpha_2, 1} ,  \cdots, a_{\alpha_n, 1}\}.$ Let $x_i := a_{\alpha_i, 1}$  and  $q_{i, j}:= q_{\alpha _i, \alpha _j} $ for convenience

\begin {Lemma}\label {ppppp4} In $\mathfrak{B} (G,
r, \overrightarrow{\chi }, u)$, $<y_i,  x_i ^s> = (s) _{q_{ii}^{-1}} x_i ^{s-1}$ for any $s\in \mathbb N.$ \end {Lemma}

\noindent {\it Proof.} It can be proved by  \cite [Proposition 2.4]{MS00} and  induction on $s.$ \hfill $\Box$

\begin {Lemma}\label {ppppp5} In $\mathfrak{B} (G,
r, \overrightarrow{\chi }, u)$, if $i_1,  i_2,  \cdots,  i_t$ are different and $0 \le e_j < {\rm ord } q_{i_j}$ for  $ 1\le j \le t $,  then $x_{i_1}^{e_1}x_{i_2}^{e_2} \cdots x_{i_t}^{e_t} \not=0.$ \end {Lemma}

\noindent {\it Proof.} We show this by induction on $m:= e_1 + e_2 + \cdots +e_t.$ It is clear when  $m=1$. Assume that $m >1.$ Considering Lemma \ref {ppppp4},  we can assume $t>1.$ By \cite [Proposition 2.4]{MS00},
\begin {eqnarray*}<y_{i_t}, x_{i_1}^{e_1}x_{i_2}^{e_2} \cdots x_{i_t}^{e_t}  > &=& x_{i_1}^{e_1}x_{i_2}^{e_2} \cdots x_{i_{t-1}}^{e_{t-1}}<y_{i_t},  x_{i_t}^{e_t}> \\
&=& (e_r) _{q_{i_r,  i_r}^{-1}} x_{i_1}^{e_1}x_{i_2}^{e_2} \cdots x_{i_{t-1}}^{e_{t-1}} x_{i_t}^{e_{t}-1} \\
 &\not=& 0  \ \ ( \hbox {by inductive assumption}).
 \end {eqnarray*}
 \hfill $\Box$

 \begin {Proposition}\label {ppppp6}

  {\rm  (i)}
$\{  x_1^{e_1} x_2^{e_2} \cdots  x_n^{e_n} \mid
 0 \le e_j < {\rm ord} (q_{jj}),  1\le j \le n \}$ is linearly independent in $\mathfrak{B} (G,
r, \overrightarrow{\chi }, u)$, Furthermore,  $\dim (\mathfrak{B} (G,
r, \overrightarrow{\chi }, u))  \ge {\rm ord } (q_{11}){\rm ord } (q_{22})\cdots {\rm ord } (q_{nn})$.

 {(\rm (ii))} $\dim (\mathfrak{B} ({\mathcal O}_{{(1, 2)}} ,  \epsilon \otimes {\rm sgn} )) \ge 2 ^{\frac {n(n-1) }{2}}$ and $\dim (\mathcal E_n) \ge 2 ^{\frac {n(n-1) }{2}}$, where $n>3$.

  \end {Proposition}

\noindent {\it Proof.} {\rm  (i)} It follows from Lemma \ref {ppppp5} that $\dim (\mathfrak{B} (G,
r, \overrightarrow{\chi }, u))  \ge {\rm ord } (q_{11}){\rm ord } (q_{22})\cdots {\rm ord } (q_{nn})$.

{\rm  (ii)} It is clear that $q_{ii} =-1$ for  $1\le i \le n$. Consequently, $\dim (\mathfrak{B} ({\mathcal O}_{{(1, 2)}} ,  \epsilon \otimes {\rm sgn} )) \ge 2 ^{\frac {n(n-1) }{2}}$ by Part {\rm (i)}.
By Theorem \ref {5.8}, $  \dim (\mathcal E_n ) \ge 2 ^{\frac {n(n-1) }{2}}$. \hfill $\Box$

Consequently,  $\dim (\mathcal E_4 ) \ge 2^6$ and $\dim (\mathcal E_5 ) \ge 2^{10}.$

\section {Finiteness conditions}\label{finiteness}

In this section we give the characteristic of finiteness of Nichols algebras in ten  ways and of {\rm FK } algebras ${\mathcal E}_n$ in six    ways.

\subsection{ Finiteness of      Nichols (braided) Lie algebras}

The fixed parameters     in \cite [Table A.1, A.2]{He05},  \cite [Table B, C]{He06a} are called  quantum numbers of  generalied Dynkin diagrams. For example,     quantum number  of Row 2 in \cite [Table A.1]{He05} is $q$; quantum numbers   of Row 9 in \cite [Table A.2]{He05} are  $q$,  $r$ and $s$.

\begin {Lemma} \label{ppp6} Assume that $R$  is an associative algebra generated by subset $L$  of $R$ as  associative algebras. If $u_1, u_2, \cdots, u_n$ is a basis of $L$ and $L$ is a Lie algebra with $[a, b] ^- = ab-ba$ for any $a, b \in L$, then $R$ is spanned  by  $\{ u_1 ^{e_1} u_2 ^{e_2} \cdots u_m^{e_m} \mid   e_1, \cdots, e_m \in \mathbb N_0 \}$ as vector spaces.
\end {Lemma}

\noindent {\it Proof.} The  universal  enveloping algebra $U(L)$ has a PBW basis $\{ u_1 ^{e_1} u_2 ^{e_2} \cdots u_m^{e_m} \mid   e_1, \cdots, e_m \in \mathbb N_0 \}$. It is clear that $R$ is a homomorphic  image of  universal  enveloping algebra $U(L)$. Consequently, $R$ is spanned  by  $\{ u_1 ^{e_1} u_2 ^{e_2} \cdots u_m^{e_m} \mid   e_1, \cdots, e_m \in \mathbb N_0 \}$ as vector spaces. \hfill $\Box$

\vskip.1in
$u$ is called an algebraic element over $k$ if there exists a polynomial  $f(x) \in k[x]$  such that $f(u)=0.$  Obviously, every nilpotent element is an algebraic element over $k$.
\begin {Theorem} \label {pppp41} Assume that $V$ is a finite dimensional braided vector space. Then $\mathfrak B(V)$ is
finite-dimensional if and only if $\mathfrak L^-(V)$ is finite-dimensional and there exists a basis, which  consists of nilpotent elements or algebraic elements over $k$,  of  $\mathfrak L^-(V)$.
\end {Theorem}

\noindent {\it Proof.} The necessity is clear. The sufficiency. Let $u_1, u_2, \cdots, u_m$
be a basis, which  consists of nilpotent elements or algebraic elements over $k$, of  $\mathfrak L^-(V)$. Thus  $\mathfrak B(V)$ is  spanned  by $\{ u_1 ^{e_1} u_2 ^{e_2} \cdots u_m^{e_m} \mid   e_1, \cdots, e_m \in \mathbb N_0 \}$ by Lemma \ref {ppp6}. Consequently, $\mathfrak B(V)$ is  finite dimensional. \hfill $\Box$

\begin {Proposition}\label {ppp41} If  $V$ is a finite dimensional  braided vector space, then

{\rm (i) }  $\mathfrak L^-(V)$  is finite dimensional if and only if $\mathfrak L^-(V)$  is nilpotent as Lie algebras;

{\rm (ii) }  $\mathfrak L(V)$  is finite dimensional if and only if $\mathfrak L(V)$  is nilpotent as braided Lie algebras;

\end {Proposition}

\noindent {\it Proof.}
{\rm (ii) }   The sufficiency. There exists  $m \in \mathbb N$ such that  $m$-degree   $\mathfrak L(V)_m$ of $\mathfrak L(V)$ is zero since $\mathfrak L(V)$  is nilpotent as braided Lie algebras. Therefore, $\mathfrak L(V)$  is finite dimensional. The necessity is clear.

Similarly we can prove {\rm (i) }. \hfill $\Box$

\vskip.1in
By Lie theorem in Lie theory, we have
\begin {Corollary}\label {ppp42}  Assume that $V$ is a finite dimensional  braided vector space.
 If  $\mathfrak L^-(V)$ is finite dimensional  or  $\mathfrak B(V)$  is finite dimensional, then ${\rm ad }\mathfrak L^-(V) \subseteq {\rm  End} \mathfrak L^-(V)$ consists of strict low triangular matrices.

\end {Corollary}

Let $ {\rm BLie} (V)$ and $ {\rm Lie} (V)$ denote the braided Lie algebra and Lie algebra generated by $V$, respectively, in $T(V)$.
\begin{Proposition} Assume that $\mathfrak B(V) = T(V)/ I$,  $J = I \cap {\rm BLie} (V) $ and $J ^-= I \cap {\rm Lie} (V).$ Then

{\rm (i)} $\mathfrak L(V)= ({\rm BLie } (V) +I) / I  \cong {\rm BLie} (V) / J$.

{\rm (ii)} $\mathfrak L^-(V) = ({\rm Lie } (V) +I) / I \cong {\rm Lie} (V) / J^-$.

\end {Proposition}

\noindent {\it Proof.} {\rm (i)} Obviously, $({\rm BLie }  (V) +I) / I$ is the  braided Lie algebra generated by $V$ in $T(V)/I$. Let $\psi$   be a map  from ${\rm BLie} (V) / J $ to  $\mathfrak L (V)\subseteq  \mathfrak B(V) := T(V)/I$ by sending  $u +J $ to $ u + I$ for any $u +J \in {\rm BLie }(V) / J$.
It is clear that  $\psi$ is an isomorphism as  braided Lie algebras.

{\rm (ii)} It is similar to Part {\rm (i)}.
\hfill $\Box$

\subsection{ Finiteness of  Nichols  algebras}
Let $r_b, r_l, r_k, r_j, r_{bm}$ denote the Baer radical, locally nilpotent radical, nil radical, Jacobson radical and Brown MacCoy radical of associative algebras (defined in \cite {Sz82}).  An algebra $A$ is called an elementary algebra if $A/ r_j(A) \cong \oplus _{i=1} ^m B_i$ with $B_i \cong k$ as algebras for $1\le i \le m$ (see \cite {ARS95}).

\begin {Theorem}\label {pp3.8''} Assume that $V$ is a finite dimensional  braided vector space. Then the following
conditions are equivalent:

{\rm (F1)} $\dim \mathfrak B(V)< \infty$;

{\rm (F2)} $\mathfrak B(V)$ is an elementary algebra with nilpotent Jacobson radical;

{\rm (F3)} $\mathfrak L^-(V)$  is  nilpotent as  Lie algebras and there exists a basis which consists of nilpotent elements or algebraic elements over $k$.

 {\rm (F4)} $\dim \mathfrak L^-(V)<\infty$  and there exists a   nilpotent basis   of  $\mathfrak L^-(V)$;

 {\rm (F5)} $\mathfrak B(V)^+$ is  nilpotent;

{\rm  (F6)}  $\dim \mathfrak L^-(V)<\infty$   with   $ r_k(\mathfrak B(V)^+ ) =  \mathfrak B(V)^+ $, or   $ r_l(\mathfrak B(V)^+ ) =  \mathfrak B(V)^+ $, or   $ r_b(\mathfrak B(V)^+ ) =  \mathfrak B(V)^+ $, or $ r_k(\mathfrak B(V) ) =  \mathfrak B(V)^+ $, or   $ r_l(\mathfrak B(V) ) =  \mathfrak B(V)^+ $, or   $ r_b(\mathfrak B(V) ) =  \mathfrak B(V)^+;$


Furthermore, if $V$ is a connected diagonal braided vector space and there is not any $m$-infinite elements $($ defined in \cite {WZZ15}$)$ with  $\dim V >1$, then the conditions above  and below  are equivalent each other:

 {\rm (F7)}  $\Delta ( \mathfrak B(V)  )$ is  an arithmetic root system and every hard super-letter $($ defined in \cite {Kh99, He05,  WZZ15}$)$ is nilpotent.

{\rm  (F8) }  $\mathfrak L(V)$ is finite dimensional.

 {\rm (F9)}  $\mathfrak L(V)$ is nilpotent as braided Lie algebras.

Moreover, if $V$ is a connected  Cartan type with $\dim V>1$,  then the conditions above  and below  are equivalent each other:

 {\rm (F10)} $V$ is a finite   Cartan type with   $1< {\rm ord } (q) < \infty$, where $q$ is the  quantum number.
\end {Theorem}

\noindent {\it Proof.} It is clear that $r (\mathfrak B(V)) = \mathfrak B(V)^+$ if and only if $r (\mathfrak B(V)^+) = \mathfrak B(V)^+$ for $r= r_b, r_k, r_{l}$ and $r_j$ by \cite {Sz82}. If $\dim \mathfrak B(V)< \infty$, then $ \mathfrak B(V)$ is nilpotent. Consequently,  $r (\mathfrak B(V)^+) = \mathfrak B(V)^+$ for $r= r_b, r_l= r_k$. If
$r_b (\mathfrak B(V)^+) = \mathfrak B(V)^+$ or $r_l (\mathfrak B(V)^+) = \mathfrak B(V)^+$, then $r_k (\mathfrak B(V)^+) = \mathfrak B(V)^+$ since $r_b (\mathfrak B(V)^+) \subseteq r_l (\mathfrak B(V)^+) \subseteq r_k (\mathfrak B(V)^+)$. However, every element in $r_k (\mathfrak B(V)^+)$ is nilpotent.

It is clear that  {\rm (F1)} and {\rm (F5)} are equivalent.
 It is clear that  {\rm (F1)} and {\rm (F6)} are equivalent by Proposition \ref {ppp41} and  Theorem  \ref {pppp41}.

{\rm (F5)}  $\Rightarrow $ {\rm (F2)} is clear. {\rm (F2)}  $\Rightarrow $ {\rm (F5)}. $ \mathfrak B(V) $ is finite dimensional since $ \mathfrak B(V) /r_j( \mathfrak B(V) )$ and $r_j( \mathfrak B(V) )$ are finite dimensional.
{\rm (F1)} and  {\rm (F3)} are equivalent by Theorem \ref {pppp41}. {\rm (F4)}  $\Rightarrow $ {\rm (F3)}. It follows from Proposition \ref {ppp41}. {\rm (F1)}  $\Rightarrow $ {\rm (F4)} is clear.

{\rm (F1)} and  {\rm (F7)} are equivalent by \cite {He05}.
{\rm (F8)}  $\Leftrightarrow $ {\rm (F1)} follows from \cite [Th.4.11]{WZZ15}. {\rm (F8)} and {\rm (F9)} are equivalent by Proposition \ref {ppp41}.
It follows from \cite [Th. 4]{He06} that {\rm (F1)} and {\rm (F10)} are equivalent.
\hfill $\Box$

\subsection{ Finiteness of   (braided)   Lie algebras in {\rm FK} algebra $\mathcal E_n$}

 We omit some of the proofs in following two subsections since they are similar to those  in the two subsections above.

\begin {Lemma}\label {pp1.1.7''} {\rm (i)} Assume that $B$ is an associative algebra and $I $ is an ideal of $B $. If $\varphi $ is a linear map from $B \otimes B$ to $B \otimes B$  and $\varphi ((I \otimes B) \cup (B \otimes I) )\subseteq B \otimes I+ I \otimes B,$
then there exists a linear map $\bar \varphi $ from  $B/I \otimes B/ I$ to $B/I \otimes B/ I$ such that $\bar \varphi (x \otimes  y) = \sum\limits _i^m (u_i +I) \otimes (v_i +I)$ for any $x, y \in B$, where $\varphi (x \otimes  y) = \sum \limits_i^m u_i  \otimes v_i$. Furthermore, if $\varphi $ is an  inverse map then so is $ \bar \varphi $.

{\rm (ii)} If $(V, C)$ is a braided vector space and $I$ is an ideal of $T(V)$ with $C(I \otimes T(V)) + C (T(V)\otimes I) \subseteq  T(V) \otimes I + I \otimes T(V)$, then $C$
can became a braiding of $T(V)/ I$ as {\rm (i)}.

\end {Lemma}

\noindent {\it Proof.} {\rm (i)}
Let $$\xi \left \{ \begin {array} {ll}  (B/I , B/I) \rightarrow B/I \otimes B/I\\
(x+I, y+I)  \mapsto  \sum \limits_{i=1} ^r (u_i +I) \otimes (v_i + I), \hbox { where } \varphi (x \otimes y) = \sum \limits_{i=1} ^r u_i  \otimes v_i.
\end {array} \right.$$ Obviously, $\xi$ is a bilinear map. Thus there exists   $$ \bar \varphi  \left \{ \begin {array} {ll}  B/I \otimes  B/I \rightarrow B/I \otimes B/I\\
(x+I) \otimes  (y+I)  \mapsto  \sum \limits_{i=1} ^r (u_i +I) \otimes (v_i + I)
\end {array} \right..$$

{\rm (ii)}  $(T(V), C)$ is a braided vector space since $C$ is a braiding of $V$.
By {\rm (i)} we complete the proof. \hfill $\Box$

\vskip.1in
If $(V, C)$ is a braided vector space and there exists a braiding $\bar C$ of $T(V)/ I$ such that  $\bar C((u +I)\otimes (v+I)) = \sum _{i=1} ^m ((u_i +I) \otimes (v_i + I))$,
where $ C(u \otimes v) = \sum _{i=1} ^m ((u_i +I) \otimes (v_i + I))$, then $\bar C$ is called  the lift of $C$ on $T(V)/ I. $  In this case we say that $C$ can be lifted  to $T(V)/I.$

If $(V, C)$ is a braided vector space with a graded ideal $I$ of $T(V)$  and $I \not= T(V)$ (i.e $I= \oplus _{i=1} ^\infty I_i$ with $I_i \subseteq V^{\otimes i}$) and $C$ can be lifted to  $T(V)/I$, then $(T(V)/I, C)$ is called a braided graded algebra. Obviously, both Nichols algebra $(\mathfrak B(V), C)$ and $\mathcal E_n$ are  braided graded algebras.
From now on, $(T(V)/I, C)$ is a braided graded algebra if without special announcement.

Let  $\mathfrak L(T(V)/I)$  denote the  braided Lie algebras
 generated by $V$ in $T(V)/I$ under braided Lie operations $[x,  y] = m ( {\rm id} - C)(y \otimes x)$
 for any homogeneous elements
$x,  y \in T(V)/I$, where $m$ is the multiplication of  $T(V)/I$. Let $\mathfrak L^-(T(V)/I)$ denote the   Lie algebras
 generated by $V$ in $(T(V)/I)^-$ under Lie operations $[x,  y]^-  = xy -yx$
 for any
$x,  y \in T(V)/I$.

\begin {Theorem} \label{ppp4000}  Assume that $(T(V)/I, C)$ is a braided graded algebra with  $ \dim V <\infty$. Then $T(V)/I$ is
finite-dimensional if and only if $\mathfrak L^-(T(V)/I)$ is finite-dimensional and there exists a basis, which  consists of nilpotent elements or algebraic elements over $k$ ,  of  $\mathfrak L^-(T(V)/I)$.
\end {Theorem}

\begin {Proposition}\label {ppp41000} Assume that $(T(V)/I, C)$ is a braided graded algebra with  $ \dim V <\infty$.

{\rm (i) }  $\mathfrak L^-(T(V)/I)$  is finite dimensional if and only if $\mathfrak L^-(T(V)/I)$  is nilpotent as Lie algebras;

{\rm (ii) }  $\mathfrak L(T(V)/I)$  is finite dimensional if and only if $\mathfrak L(T(V)/I)$  is nilpotent as braided Lie algebras;

\end {Proposition}

\begin{Proposition} \label {p7.10} Assume that $(T(V)/I, C)$ is a braided graded algebra with  $ \dim V <\infty$,  $J = I \cap {\rm BLie} (V) $ and $J ^-= I \cap {\rm Lie} (V).$ Then

{\rm (i)} $\mathfrak L( T(V)/I)= ({\rm BLie } (V) +I) / I  \cong {\rm BLie} (V) / J$.

{\rm (ii)} $\mathfrak L^-(T(V)/I) = ({\rm Lie } (V) +I) / I \cong {\rm Lie} (V) / J^-$.

\end {Proposition}

\subsection{ Finiteness of {\rm FK} algebra $\mathcal E_n$}

\begin {Theorem}\label {pp3.8''000} Assume that $(T(V)/I, C)$ is a braided graded algebra with  $ \dim V <\infty$. Then the following
conditions are equivalent:

{\rm (F1)} $\dim T(V)/I< \infty$;

{\rm (F2)} $T(V)/I$ is an elementary algebra with nilpotent Jacobson radical;

{\rm (F3)} $\mathfrak L^-(T(V)/I)$  is  nilpotent as  Lie algebras and there exists a basis which consists of nilpotent elements, or algebraic elements over $k$.

 {\rm (F4)} $\dim \mathfrak L^-(T(V)/I)<\infty$  and there exists a   nilpotent basis   of  $\mathfrak L^-(T(V)/I)$;

 {\rm (F5)} $(T(V)/I)^+$ is  nilpotent;

{\rm  (F6)}  $\dim \mathfrak L^-(T(V)/I)<\infty$   with   $ r_k((T(V)/I)^+ ) =  (T(V)/I)^+ $, or   $ r_l((T(V)/I)^+ ) =  (T(V)/I)^+ $, or   $ r_b((T(V)/I)^+ ) =  (T(V)/I)^+ $, or $ r_k(T(V)/I ) =  (T(V)/I)^+ $, or   $ r_l(T(V)/I ) = ( T(V)/I)^+ $, or   $ r_b(T(V)/I ) = ( T(V)/I)^+;$


\end {Theorem}

\section {Irreducible representations }\label{irreps section}

In this section we give all irreducible representations of finite dimensional Nichols algebras and of  finite dimensional {\rm FK } algebras ${\mathcal E}_n$.

\begin {Lemma}\label {ppp21}  Assume that  $R$ is an associative algebra. Then $V$ is an irreducible $R$-module if and only if $V$ is  an irreducible $R/r_j(R)$-module.

\end {Lemma}

\noindent {\it Proof.}  If $V$ is an irreducible $R$-module, then $r_j(R) \subseteq (0:V)_R := \{ x\in R \mid x V =0\}. $ Thus $V$ become an $R/r_j(R)$-module under natural module operation. If $0\not= N $ is an $R/ r_j(R)$- submodule of $V$, then $ N $ is an $R$- submodule of $V$ and $N= V$, which implies that  $V$ is  an irreducible $R/r_j(R)$-module. Conversely,  if $V$ is  an irreducible $R/r_j(R)$-module, then  $V$ is  an  $R$-module. Let $0\not= N $ be  an $R$- submodule of $V$. It is clear $r_j(R) \subseteq  (0:V)_R \subseteq (0:N)_R $, which implies that $N$ is an $R/ r_j(R)$ -module and $N = V$.

\hfill $\Box$

\begin {Lemma}\label {ppp71} Assume that  $A = \oplus _{i=0} ^\infty A_i$ is a graded algebra with $A_0=k.$ If $r_j (A) = A_{>0} := \oplus _{i=1} ^\infty  A_i$, then  $A$ has  only  one  irreducible module up to isomorphisms and  the dimension of this module is one with trivial action.
\end {Lemma}

\noindent {\it Proof.}  It is clear that $A/ r_j (A) \cong k$, $k$ has  only  one  irreducible module up to isomorphisms and  the dimension of this module is one with trivial action. Considering  Lemma \ref {ppp21} we complete the proof. \hfill $\Box$

\begin {Theorem}\label {ppp31}  Assume that $(T(V)/I, C)$ is a braided graded algebra with  $ \dim V <\infty$. If $(T(V)/I, C)$ is finite dimensional, then $(T(V)/I, C)$ has  only  one  irreducible module up to isomorphisms and  the dimension of this module is one with trivial action.

For example, $(T(V)/I, C)= \mathfrak B(V)$ or   $\mathcal E_n$.

\end {Theorem}

\noindent {\it Proof.}  It follows from  Lemma \ref {ppp71} and Theorem \ref {pp3.8''000} {\rm (F6)}. \hfill $\Box$

\section { {\rm YD} module of finite commutative groups}\label{YD modules}

In this section we give some conditions for a braided vector space to become a {\rm YD} module over finite commutative group. We obtain the sufficient and necessary condition for Nichols algebra $\mathfrak B(M)$ of  {\rm YD} module  $M$ over  $\mathbb W_n$  with ${\rm supp } (M) \subseteq A$
 to be finite dimensional.

Let $G=(a_1)\times (a_2)\times \cdots \times (a_t)$ be  a finite abelian group  with $ N_i= {\rm ord} (a_i)$  for $1\le i \le t$.  If $N_1 = N_i$ for $1\le i \le t,$ then $G$ is called a weak elementary group. If $G$ is a  weak elementary group and  $N_1$ is a prime number, then $G$ is called an  elementary group.

$A_n, D_n, E_6, E_7, E_8$ is called  finite laced Cartan type since there exists  one edge between  two vertices at most.

\subsection{Bicharacters}

Let $G$ and $H$ be two groups. If $r$ is a map from $G\times H$ to $k^*$ $(k^* := \{x \in k \mid x\not=0\})$ with
$r(ab,c)=r(a,c)r(b,c)$ and $r(a,cd)=r(a,c)r(a,d)$ for $a, b \in G, c, d \in H$, then $r$
is called a bicharacter over pair $(G, H)$. If $G=H$, then $r$  is called a bicharacter over group $G.$

\begin{Lemma} \label {pp2.1.1} If $G=(a)$ and $H=(b)$, then $r$ is a bicharacter over $(G,H)$ if and only if $r$ is a map from $G \times H$ to $k^*$ such that
$r(a^{m},b^{n})=r(a,b)^{mn}$ with $r(a,b)^{N}=1$ when $N = {\rm ord } (a)$ or $N = {\rm ord } (b)$ for any  $m, n \in \mathbb N$.
\end {Lemma}

\begin{Proposition} \label {pp2.1.3} Let $G=\oplus _{i\in I} G_i$ be the direct sum of groups $G_i$ for $i\in I$.

{\rm (i)} If $r$ is a bicharacter of group $G$, then  restricted map $r|_{(G_i,G_j )}$ is a bicharacter over pair $(G_i,G_j)$ for $i,j\in I$.

{\rm (ii)} Conversely, if $r_{i,j}$ is a bicharacter over pair $(G_i,G_j)$ for any $i,j\in I$ and define $r(x,y)=\prod\limits_{i,j\in I}r_{i,j}(x_i,x_j)$ for any $x=\{x_i\}$,$y=\{y_i\}\in G$,  then $r$ is a bicharacter over $G$.
\end {Proposition}

\noindent {\it Proof.} It is clear. \hfill $\Box$

\begin{Proposition} \label {pp2.1.4} Let $G=\oplus _{i\in I} G_i$ be the direct sum of cyclic group $G_i=(a_i)$  for any $i \in I$. Set ${\mathcal R}:=\{r \mid r$  is a bicharacter over $G$ $\}$  and ${\mathcal Q}:=\{\{b_{i,j}\}_{i, j \in I} \mid   b_{ij }
\in k^* ; b_{i,j}^{{\rm ord } (a_j)}=b_{i,j}^{{\rm ord } (a_i)}=1,  \hbox { for any } i,j\in I\}$. Define $\phi :$
 $\begin {cases} {\mathcal R} \longrightarrow {\mathcal Q },\\
r \longmapsto \{r(a_i,a_j)\}_{i, j \in I}.
\end {cases}$ Then $\phi$ is a bijection.
\end {Proposition}

\noindent {\it Proof.} We show this by following several steps.

{\rm (i) } $\phi$ is well defined.

{\rm (ii) }  $\phi$ is surjective. In fact,  for any $\{b_{i,j}\}\in \mathcal Q$, define a bicharacter $r$ over $G$ as follows:

$r(a_i,a_j)=b_{i,j}$ for any $i,j\in I$ and $r(x,y)=\prod \limits_{i,j\in I}(b_{i,j})^{m_in_j}$, where $x=\{a_i^{m_i}\}$,$y=\{a_i^{n_i}\}\in G$.
 In fact, for any $x = \{x_i\}= \{a_i^{m_i}\},$ $ x' = \{x_i'\}= \{a_i^{m_i'}\},$ $ y = \{y_i\}=\{a_i^{n_i}\}, $ $y' = \{y_i'\}=\{a_i^{n_i'}\} \in G$ with $x=x'$ and $y=y'$,
we have that
  $$\frac   {r(x, y)} {r(x', y')}= \prod \limits_{i,j \in I}  r(a_i, a_j) ^{m_i(n_j-n_j') + (m_i-m_i') n_j'}
  =1,$$
  which implies $r(x, y)= r(x', y')$. That is, $r$ is well defined.

{\rm (iii) } It is clear that $\phi$ is injective. \hfill $\Box$

\vskip.1in
We have the  similar result  for $I= \{
1, 2, \cdots, t\}$, i.e. $G = (a_1) \times (a_2) \times \cdots \times (a_t)$.

\subsection{{\rm YD} module of finite commutative groups}

If $H$ is a subgroup of $G$ and there exists a subgroup $H'$ such that $G = H\times H'$ is a direct sum of $H$ and $H'$, then $H$ is called a direct summand of $G$.

 Let $\overrightarrow{g} = (g_1, g_2, \cdots, g_\theta) \in G^\theta $, $
\overrightarrow{\chi} = (\chi _1, \cdots, \chi_\theta) \in ( \widehat G, \widehat G, \cdots, \widehat G    ).$   We call $(G,\overrightarrow{g} , \overrightarrow{\chi} )$ an element system with characters, written as   ${\rm ESC } (G,\overrightarrow{g} , \overrightarrow{\chi} )$ (see \cite [Def. 2.1]{ZZC07}). It is clear that  $V $  is a $G$-${\rm YD}$  if and only if there exists  ${\rm ESC } (G, \overrightarrow{g} , \overrightarrow{\chi} )$
 and a basis $x_1, x_2, \cdots, x_\theta$ of $V$ such that module operation  $g \cdot x_i =\chi _i (g) x_i, $  and comodule operation $\delta (x_i) = g_i \otimes x_i,$  $1\le i \le \theta.$

Let $(q_{ij})_{\theta \times \theta }$ be a matrix. If  $q_{ij} q_{ji}  \left \{  \begin{array}{ll}
 \not= 1,  & \mbox {when }  \mid j-i \mid = 1\\
  =1,   & \mbox {when}  \mid j-i \mid \not= 1    \\
\end{array}\right. $ for any $1\le i\not= j \le \theta,$ then  $(q_{ij})_{\theta \times \theta }$ is called a  chain or labelled  chain.
If $(q_{ij})_{\theta \times \theta }$ is  a  chain and \begin {eqnarray} \label {ppe1}(q_{11} q_{1, 2} q_{2, 1} -1)(q_{11} +1)=0;
 (q_{\theta, \theta } q_{\theta , \theta-1} q_{\theta-1, \theta} -1)(q_{\theta,\theta} +1)=0;\end {eqnarray}
  \begin {eqnarray} \label {ppe2}
   q_{ii} +1= q_{i, i-1} q_{i-1, i}q_{i, i+1} q_{i+1, i}-1=0 \mbox { or  }
  q_{ii}q_{i, i-1} q_{i-1, i}=q_{ii} q_{i, i+1} q_{i+1, i}=1,\end  {eqnarray} ¶Ô $1<i < \theta$,
   then  $(q_{ij})_{\theta \times \theta }$ is called a simple chain (see \cite [Def.1]{He06a}). Let
$q:= q_{\theta, \theta} ^2 q_{\theta -1, \theta}q_{\theta , \theta-1}$.
If $(q_{ij})_{\theta \times \theta }$ is a simple chain and  $q_{i, i-1} q_{i-1, i} =q^{-1}$, $1< i \le \theta$, then  $(q_{ij})_{\theta \times \theta }$ is called a pure simple chain, otherwise, it is called a mixed simple chain.
.

\begin{Lemma} \label{ppp1} If $(q_{ij})_{\theta \times \theta }$ is a simple chain with  $q_{ii} \not= -1$ $(1\le i \le \theta)$,  then   $(q_{ij})_{\theta \times \theta }$ is pure simple chain and  $q_{jj} =q$ for $1\le j\le \theta$.
\end{Lemma}
\noindent {\it Proof.} By  (\ref {ppe1}), $q_{\theta, \theta} =q = ( q_{\theta , \theta-1} q_{\theta-1, \theta} )^{-1}$. It follows from  (\ref {ppe2}) that   $q_{ii} = q = ( q_{i, i+1} q_{i+1, i})^{-1}= ( q_{i, i-1} q_{i-1, i})^{-1} $ for  $1< i < \theta$. Consequently, $q_{11} = q$ by  (\ref {ppe1}).
\hfill $\Box$

\begin{Lemma} \label{ppp2} Let $G$ be a commutative group with odd $\mid G \mid$. If connected $(q_{ij})_{\theta \times \theta }$ is a mixed simple  chain, or there exist $i, j$ such that ${\rm ord } (q_{ij})$  is even ( e.g. $q_{ij} =-1$ ), then  $(q_{ij})_{\theta \times \theta }$ is not a $G$- {\rm YD} module.
\end {Lemma}
\noindent {\it Proof.} It follows from Lemma \ref  {ppp1}.
\hfill $\Box$

\vskip.2in
 Let  $H_2$ denote the following two cases:   $\begin{picture}(100,     15)
\put(27,     1){\makebox(0,    0)[t]{$\bullet$}}
\put(60,     1){\makebox(0,     0)[t]{$\bullet$}}
\put(28,     -1){\line(1,     0){33}}
\put(22,    10){$\zeta$}
\put(58,     10){$q$}
\put(40,     5){$q^{-1}$}
\end{picture}$  and
$\begin{picture}(100,     15)
\put(27,     1){\makebox(0,    0)[t]{$\bullet$}}
\put(80,     1){\makebox(0,     0)[t]{$\bullet$}}
\put(28,     -1){\line(1,     0){53}}
\put(22,    10){$\zeta$}
\put(78,     10){$q^{-1}\zeta $}
\put(40,     5){$\zeta ^{-1}q $}
\end{picture},$ where $\zeta \in R_3, q \in k^* \setminus \{1, \zeta, \zeta ^2\}$. The two cases are in Row 6   in \cite [TableA.1] {He05}.\\

 Let  $H_3$ denote the following two cases:
  $\begin{picture}(100,     15)
\put(27,     1){\makebox(0,    0)[t]{$\bullet$}}
\put(60,     1){\makebox(0,     0)[t]{$\bullet$}}
\put(93,    1){\makebox(0,    0)[t]{$\bullet$}}

\put(28,     -1){\line(1,     0){33}}
\put(61,     -1){\line(1,     0){30}}

\put(22,    10){$\zeta$}
\put(58,     10){$\zeta $}
\put(91,     10){$\zeta ^{-3}$}

\put(40,     5){$\zeta^{-1}$}
\put(73,     5){$\zeta^{-1}$}

\end{picture}$  and    $\begin{picture}(100,     15)
\put(27,     1){\makebox(0,    0)[t]{$\bullet$}}
\put(60,     1){\makebox(0,     0)[t]{$\bullet$}}
\put(93,    1){\makebox(0,    0)[t]{$\bullet$}}

\put(28,     -1){\line(1,     0){33}}
\put(61,     -1){\line(1,     0){30}}

\put(22,    10){$\zeta$}
\put(58,     10){$\zeta ^{-4}$}
\put(91,     10){$\zeta ^{-3}$}

\put(40,     5){$\zeta^{-1}$}
\put(73,     5){$\zeta^{4}$}

\end{picture}$,
  where $\zeta \in R_9 $. The two cases are in Row 18   in \cite [TableA.2] {He05}.

\begin{Proposition} \label{ppp3}
 Assume that  connected $(q_{ij})_{\theta \times \theta }$  is of an arithmetic root system  $\theta >1$. If  one of the following four  cases holds,  then $(q_{ij})_{\theta \times \theta }$ is a connected finite Cartan type  or $H_2$ or $H_3.$

${\rm (i)}$   $(q_{ij})_{\theta \times \theta }$ is a $G$-{\rm YD} module and  $G$ is a commutative group with odd $\mid G \mid $.

${\rm (ii)}$  $(q_{ij})_{\theta \times \theta}$ is not the first diagram of  Row 9  in  \cite [Table A.1] {He05} and  ${\rm ord } (q_{ii})$ is odd  for any $1 \le i \le \theta $.

${\rm (iii)}$  $(q_{ij})_{\theta \times \theta}$ is not   Row 12, the last diagram of   Row 16,   Row 18  in  \cite [Table A.2] {He05},
 the first diagram in  Row 15, the first diagram of   Row 17 in  \cite [Table B] {He06a} and Row 5 in \cite [Table C] {He06a} with  $\theta \ge 3$ and  $q_{ii} \not=-1$  for  $1 \le i \le \theta $.

${\rm (iv)}$  $(q_{ij})_{\theta \times \theta}$ is not Row 5 in  \cite [Table C] {He06a} with  $\theta >4$ and  $q_{ii} \not=-1$  for  $1 \le i \le \theta $.

Furthermore, if  ${\rm ord } (q_{ii})>3$ is odd  for $ 1\le i \le \theta$, then $(q_{ij})_{\theta \times \theta }$ is a  finite Cartan type

\end {Proposition}
\noindent {\it Proof.} It follows from  Lemma \ref {ppp1}, Lemma \ref {ppp2}, \cite [Table A.1, A.2] {He05} and \cite [Table B, C] {He06a}. \hfill $\Box$

\vskip.1in
If $(q_{ij}) _{\theta\times \theta}$ is a matrix and there exists an ${\rm ESC } (G, \overrightarrow{g} , \overrightarrow{\chi} )$ such that  \begin {eqnarray} \label {ppe2.2.2}q_{ij} = \chi _j (g_i), \end {eqnarray} then  we say that
$(q_{ij}) _{\theta\times \theta}$  is determined by ${\rm ESC } (G, \overrightarrow{g} , \overrightarrow{\chi} )$. In this case $(q_{ij})_{\theta \times \theta}$ is said to be  a $G$-{\rm YD} module.

\begin{Proposition} \label{ppp4}  Assume that connected $(q_{ij})_{\theta \times \theta }$  is  a $G$- {\rm YD} module and is of an arithmetic root system with $\theta >1.$

{ \rm (i)} If $G$ is an elementary $p$-group with  $p\not=2$, then
$(q_{ij})_{\theta \times \theta }$ is a connected finite Cartan type with quantum number $q$ and ${\rm ord} (q) =p$.

{ \rm (ii)} If $G$ is an elementary $p$-group with  $p=2$, then
$(q_{ij})_{\theta \times \theta }$ is a connected finite laced Cartan type with quantum number $q=-1$.

\end {Proposition}
\noindent {\it Proof.} { \rm (i)}   follows from Proposition  \ref {ppp3}.

{ \rm (ii)}    It is clear.  \hfill $\Box$

Obviously, if  $r$ is a bicharacter over finite commutative $G$ with  $g \in G$, then  there exists $\chi \in \hat G$, such that  $\chi (h) = r(g, h)$, $h \in G.$

\begin{Proposition} \label{pp2.2.3}

Assume that  $H= H '\times H''$.
If $ \rm ESC (H', \overrightarrow{g'}, \overrightarrow{\chi'}  )$  and  $ \rm ESC (H'', \overrightarrow{g''}, \overrightarrow{\chi''}  )$ are {\rm ESC} over $H'$ and $H''$, respectively, then $ {\rm ESC} (G, \overrightarrow{g}, \overrightarrow{\chi}  )$ is an
{\rm ESC} over
$H$, where  $\overrightarrow{\chi} = \overrightarrow{\chi'} \otimes \overrightarrow{\chi''} := \{ \chi _i '\otimes \chi _j'' \}_{i, j}$  and  $\overrightarrow{g} = \overrightarrow{g'} + \overrightarrow{g''}  :=\{ (g_i, g_j) \}_{i, j}$. Conversely,  if  $ {\rm ESC} (G, \overrightarrow{g}, \overrightarrow{\chi}  )$ is an
{\rm ESC} over
$H$ with  $\overrightarrow{\chi} = \overrightarrow{\chi'} \otimes \overrightarrow{\chi''}$  and  $\overrightarrow{g} = \overrightarrow{g'} + \overrightarrow{g''} $, then  $ \rm ESC (H', \overrightarrow{g'}, \overrightarrow{\chi'}  )$  and  $ \rm ESC (H'', \overrightarrow{g''}, \overrightarrow{\chi''}  )$ are {\rm ESC} over $H'$ and $H''$, respectively.

\end{Proposition}

\noindent {\it Proof.} It is clear. \hfill $\Box$

\begin {Lemma} \label{pp2.2.4'} If $H$ is a subgroup of $G$ with  $ \rm ESC (H, \overrightarrow{g}, \overrightarrow{\chi} )$ and $\chi_i'$ is  an   extension on $G$
of $\chi_i$, then  the matrices  of $\rm ESC (G, \overrightarrow{g}, \overrightarrow{\chi'} )$ and  $\rm ESC (H, \overrightarrow{g}, \overrightarrow{\chi} )$ are the same.
\end {Lemma}

We  give a sufficient condition for $(q_{ij}) _{\theta \times \theta}$ to be a {\rm YD}-module over finite commutative group with ${\rm ord } (a_i) = N_i$.
\begin {Theorem} \label{pp2.2.4} Assume that  $G= (a_1) \times (a_2) \times \cdots \times (a_t) $ is a finite commutative group and $(q_{ij}) _{\theta \times \theta}$ is a braided matrix. If
  there exist distinct $1\le \lambda _1 , \lambda_2\cdots, \lambda_\theta \le t$ such that $q_{ij}^{N_{\lambda _i}} = q_{ij}^{N_{\lambda_j}} =1$  for any $1\le i. j\le \theta$ with  $\theta \le t,$ then there exists $ \rm ESC (G, \overrightarrow{g}, \overrightarrow{\chi} )$ such that  $\chi _j (g_i) = q_{ij}$ and  $g_i = a_i$ for  $1\le i, j\le \theta$. Furthermore, $(q_{ij}) _{\theta \times \theta}$ is  a  $G$-{\rm YD}-module.

\end {Theorem}
\noindent {\it Proof.} Let $\lambda_1, \lambda_2, \cdots, \lambda _t$ be a permutation of $1, 2, \cdots, t$.
Let $H '= (a_{\lambda_1}) \times \cdots \times (a_{\lambda_\theta})$,  $H'' = (a_{\lambda _{\theta +1}}) \times \cdots \times (a_{\lambda_t})$. By Proposition \ref {pp2.1.4}, there exists a bicharacter $r$ of  $H'$ such that  $r(a_{\lambda_i} , a_{\lambda_j}) = q_{ij}$. Let  $\chi _j' \in \widehat {H'}$ such that  $\chi _j '(g_i') = q_{ij}$ and   $g_i' = a_{\lambda_i}$ for  $1\le i, j\le \theta$. Consequently
$ \rm ESC (H', \overrightarrow{g'}, \overrightarrow{\chi'}  )$ is an element system with characters of  $H'$. Considering Lemma \ref {pp2.2.4'} we complete the proof. \hfill $\Box$

\vskip.1in
Following the Theorem above, we have
\begin {Corollary} \label{ppp61} If $(q_{ij})_{\theta\times \theta}$ is a braided matrix with finite orders $($i.e.  ${\rm ord }(q_{ij}) < \infty$ for $1\le i, j \le \theta $$ )$, then there exists a finite commutative group $G$ such that $(q_{ij})_{\theta\times \theta}$ is a {\rm YD} -module over $G.$
\end {Corollary}

If there exists a subgroup $H$ of a finite commutative group $G$ such that $ G = (g_1) \times (g_2) \times \cdots \times (g_\theta) \times H $, then  $ \rm ESC (G, \overrightarrow{g}, \overrightarrow{\chi} )$ is called to be full.  $(q_{ij})_{\theta \times \theta}$
is called to be full if it is a braided matrix of full $\rm ESC (G, \overrightarrow{g} , \overrightarrow{\chi} )$.

\begin {Proposition} \label{pp2.2.5} $(q_{ij})_{\theta \times \theta}$ is of a full $G$-{\rm YD} module if and only if there exist a sungroup $H$ and elements $b_i \in G$ for $1\le i \le \theta$ such that $ G = (b_1) \times (b_2) \times \cdots \times (b_\theta) \times H $ and
$q_{ij} ^{{\rm ord } (b_i)} = q_{ij} ^{{\rm ord } (b_j)} = 1$ for $1\le i, j \le \theta$.

\end {Proposition}

\noindent {\it Proof.}  Necessity. It is clear.

Sufficiency. It is clear that there exist $b_i \in H$ for $1+\theta \le i \le t$ such that $G = (b_1) \times \cdots \times (b_t)$. By Theorem \ref {pp2.2.4}, $(q_{ij})_{\theta \times \theta}$ is of a full $G$-{\rm YD} module. \hfill $\Box$

\subsection{ {\rm YD} module of weak elementary groups}

\begin{Proposition} \label{pp2.3.3}
If  $G= (a_1) \times (a_2) \times \cdots \times (a_t) $ is a weak elementary group with $\theta \le t$ and  $N= {\rm
\rm ord } (a_1)$, then  $(q_{ij}) _{\theta  \times \theta}$ is  $G$-{\rm YD} module if and only if ${\rm ord } (q_{ij}) \mid N$, $1\le i, j \le \theta$.

\end{Proposition}
\noindent {\it Proof.} If $(q_{ij}) _{\theta \times \theta}$ is a $G$-{\rm YD} module, then there exists  $\rm ESC (G, \overrightarrow{g}, \overrightarrow{\chi} )$  such that its matrix is  $(q_{ij}) _{\theta \times \theta}$ and   $\chi _j (g_i) = q_{ij} $.  Consequently ${\rm ord } (q_{ij}) \mid N$.

Conversely, if  ${\rm ord } (q_{ij}) \mid N$, then there exists  $ \rm ESC (G, \overrightarrow{g} , \overrightarrow{\chi} )$ such that its matrix is  $(q_{ij}) _{\theta \times \theta}$ with  $\chi _j (a_i) = q_{ij}$ for  $1\le i, j\le \theta$ by Theorem  \ref {pp2.2.4}. \hfill $\Box$

\begin{Remark} It can be omitted that the conditions: $\theta \le t$  in the necessity.
\end{Remark}

\begin{Proposition} \label{pp2.3.4} Assume that   $G=(a_1) \times (a_2) \times \cdots \times (a_t) $ is a weak elementary group with $N= {\rm
\rm ord } (a_1)$, If
 $(q_{ij}) _{\theta \times \theta}$ is of a $G$-{\rm YD} module  and a finite  Cartan type with quantum number  $q$, then  ${\rm ord }(q) \mid N $.

\end{Proposition}
\noindent {\it Proof.} There exists  $ \rm ESC (G, \overrightarrow{g}, \overrightarrow{\chi} )$ such that   $\chi _j (g_i) = q_{ij}$ for  $1\le i, j\le \theta$. It is clear $g_i^N =1$ and $\chi _j (g_i)^N = q_{ij}^N =1$ for  $1\le i, j\le \theta$. Consequently,
${\rm ord }(q) \mid N $  since there exist $i$ such that $q_{ii} =q$. \hfill $\Box$

If $N\in \mathbb N$ and $g^N =1$ for any $g \in G$, then $N$ is called a generalized order of $G.$ For example, if   $G=(a_1) \times (a_2) \times \cdots \times (a_t) $ is a weak elementary group with $N= {\rm
\rm ord } (a_1)$, then $N$ is a generalized order of $G.$

\begin{Proposition} \label{p2.3.4} Assume that $(q_{ij}) _{\theta \times \theta}$ is of a connected $G$-{\rm YD} module with $\theta >1$.

${\rm (i)}$  If $N$ is a generalized order of $ G$,  then  the order of every quantum number can divides $N$ except  last   diagram  in Row 5, Diagram 4  in Row 7 in \cite [Table A.2] {He05}, last   diagram of Row 7,    Diagram 5  of Row 9, Diagram 4, 5, 6  of Row 11 in \cite [Table B] {He06a},
 Diagram 2  in Row 4 in \cite [Table C] {He06a}.

${\rm (ii)}$  If $N$ is a generalized order of $ G$,  then  the square of  order of every quantum number can divides $N$ in   last   diagram  in Row 5, Diagram 4  in Row 7 in \cite [Table A.2] {He05}, last   diagram of Row 7,    Diagram 5  of Row 9, Diagram 4, 5, 6  of Row 11 in \cite [Table B] {He06a},
 Diagram 2  in Row 4 in \cite [Table C] {He06a}.

${\rm (iii)}$  If   $G=(a_1) \times (a_2) \times \cdots \times (a_t) $ is a weak elementary group with $N= {\rm
\rm ord } (a_1)$ and
 $(q_{ij}) _{\theta \times \theta}$ is  a finite  Cartan type with quantum number  $q$, then  ${\rm ord }(q) \mid N $.

\end{Proposition}
\noindent {\it Proof.} {\rm (i)} and {\rm (i)} We check this step by step in \cite [Table A.1, A.2] {He05} and \cite [Table B, C] {He06a}.

{\rm (iii)} It follows {\rm (i)}.  \hfill $\Box$

\begin{Remark} It can be omitted that the conditions: $\theta \le t$  in the necessity of {\rm (ii)}.
\end{Remark}


\subsection{ {\rm YD} module of elementary groups}

\begin {Proposition} \label{pp2.4.2}
 If  $G$ is an elementary $p$-group with $2 \le \theta \le t$ and $q_{ij} =q_{ji}$ for $1\le i, j\le \theta$,  then  $(q_{ij}) _{\theta \times \theta}$ is  of   a connected¡¡ $G$-{\rm YD} module  and an arithmetic root system if and only if
$(q_{ij}) _{\theta \times \theta}$ is of a connected finite Cartan type with quantum number $q$ and ${\rm ord } (q) = p$.

\end {Proposition}
\noindent {\it Proof.} Necessity follows form Proposition \ref {ppp4}.

Sufficiency follows from Theorem   \ref {pp2.2.4}.  \hfill $\Box$

\begin{Remark} It can be omitted that the conditions: $\theta \le t$ and $q_{ij} =q_{ji}$ for $1\le i, j\le \theta$ in the necessity.
\end{Remark}

\subsection {Nichols algebras of  {\rm YD} module   over  $\mathbb W_n$ }

 In this subsection we give the sufficient and necessary condition for Nichols algebra $\mathfrak B(M)$ of  {\rm YD} module  $M$ over  $\mathbb W_n$  with ${\rm supp } (M) \subseteq A$
 to be finite dimensional.

\begin {Lemma}\label {ppp60} Assume that  $M$ is a finite dimensional  $G$- YD module. If $G_1$ is a finite  commutative subgroup of $G$ and $supp M \subseteq G_1$, then $M$ is a $G_1$-YD module with  diagonal braiding.
  \end {Lemma}
\noindent {\it Proof.} It follows from \cite [Lemma 2.3]{ZZC07}. \hfill $\Box$

\begin {Theorem}\label {pp3.8'} Let $G= \mathbb W_n $.
 Assume that  $\sigma _1, \cdots, \sigma _m  \in A$ and
$M = M({\mathcal O}_{\sigma _1} , \rho  ^{(1)})\oplus M({\mathcal
O}_{\sigma _2}, \rho ^{(2)}) \oplus \cdots \oplus M({\mathcal O}_{\sigma _m},
\rho ^{(m)})$ is a
 {\rm YD} module over $kG$ with $ \rho  ^{(i)} \in  \widehat
  {G ^{\sigma_i}}$ for $i=1, 2, \cdots, m$. Then the following conditions are equivalent.

  {\rm (i)}
$\mathfrak B (M)$ is  finite dimensional.

 {\rm (ii)} Every connected component of generalized Dynkin diagram of   $M$ is a finite laced Cartan type with
$\rho ^{(i)}(\sigma_i) = -{\rm id}$ for $1\le i \le m$.

 {\rm (iii)} Every connected component of generalized Dynkin diagram of   $M$ is a finite laced Cartan type.

  {\rm (iv)} Every connected component of generalized Dynkin diagram of   $M$ is a finite  Cartan type.

 {\rm (v)}  $\dim \mathfrak B(M({\mathcal O}_{\sigma _1} , \rho  ^{(1)})\oplus M({\mathcal
O}_{\sigma _2}, \rho ^{(2)}) \oplus \cdots \oplus M({\mathcal O}_{\sigma _j},
\rho ^{(j)}))< \infty$  for $1 \le j \le m.$

  \end {Theorem}
\noindent {\it Proof.} {\rm (i)} $\Rightarrow {\rm (ii)} $. $M$ is an $A$-Yetter Drinfeld module by Lemma \ref {ppp60}. By Proposition \ref {pp2.4.2} and \cite {AZ07},
every connected component of $M$ is a finite laced  Cartan type  with
$\chi^{(i)}(\sigma_i) = -1$ for $1\le i \le m$.

{\rm (ii)} $\Rightarrow {\rm (iii)} \Rightarrow {\rm (iv)}$ It is clear.

{\rm (iv)} $\Rightarrow {\rm (i)} $. Since $G$ is a finite group, the order of  quantum number of every connected component of generalized Dynkin diagram is  finite. By \cite [Lemma 6.4] {WZZ15}, $1< p_{u, u} < \infty$ since $p_{u, u}$ is   quantum number  or square of  quantum number of connected component for any hard super-letter $u.$ Consequently, $\mathfrak B (M)<\infty.$

{\rm (i)} $\Leftrightarrow {\rm (v)} $ It is clear.   \hfill $\Box$

\section {Lyndon basis of Nichols algebras of finite Cartan types  and enveloping algebras of Lie algebras} \label{Lyndon basis}

In this section  it is shown that hard braided Lie Lyndon word,  standard  Lyndon word,  Lyndon-basis-path, hard  Lie Lyndon word and    standard  Lie Lyndon word    are the same with respect to $ \mathfrak B(V)$, Cartan matrix $A_c$ and $U(L^+)$, respectively, where  $V$ and $L$ correspond to the same finite Cartan matrix $A_c$.

\subsection {  Standard basis of {\rm FK} algebra $\mathcal E_n$}

Let $A := \{ x_1, x_2, \cdots, x_\theta \}$  be an alphabet and a basis of $V = {\rm span } (A)$; $A^*:= \{u \mid u \hbox { is a word } \}$. Let $I$ be a graded ideal of $T(V)$ as algebras.

\begin {Definition} \label {pp2} (\cite [Def. 1] {Kh99})
A word $u$ is called a Lyndon word if $| u| =1$ or $| u| \geq2$,  and for each
representation $u=u_1u_2$,  where $u_1$and $u_2$ are nonempty
words,  the inequality $u<u_2u_1$ holds.
\end {Definition}

$u$ is called a standard   word with respect to $T(V)/I$ if $u $ can not be  written as  linear combination of  strictly greater  words in $T(V)/I$.

Let $ {\rm SW} (T(V)/I)=: \{ u \in A^* \mid  u \hbox { is a standard word with respect to } $T(V)/I$ \}$, written as {\rm SW} in short;

 ${\rm L }:= \{ u\in A^* \mid u \hbox { is a Lyndon word}\}. $ Notice that we view
$ {\rm SW} (T(V)/I)$ and ${\rm L }$ are in $T(V)/I$ often for convenience.

$u$ is called  Standard Lyndon ward if $u \in {\rm SW } \cap {\rm L}.$

Lyndon  word $u$ is called a standard  Lie  Lyndon  word with respect to $T(V)/I$ if $[u]^- $ can not be  written as  linear combination of  strictly greater  Lie  Lyndon words in $T(V)/I$.

${\rm  SLLW}(T(V)/I)=: \{ u\in {\rm L}  \mid  u \hbox { is a standard Lie Lyndon word } \}$.

\begin {Lemma}\label {pp10.7''}
{\rm (i)} $ {\rm SW} $ is a basis of $T(V)/I.$

{\rm (ii)} Any factor of a standard  word is a standard  word.

{\rm (iii)} If $u$ is a standard word, then  $u = u_1u_2 \cdots u_r$ with $u_! \ge u_2\ge \cdots \ge u_r $ and $u_i \in {\rm SW } \cap {\rm L}$ for $1\le i \le r.$

\end {Lemma}
\noindent {\it Proof.}  {\rm (i)} If a word $u \notin {\rm SW},$ then $u = \sum\limits _{v \in {\rm SW},  \mid v\mid = \mid u \mid, v >u} a_v v$. It is clear that ${\rm SW}$ is linearly independent in $T(V)/I.$ Consequently, ${\rm SW}$ is a basis of $T(V)/I.$

{\rm (ii)} It is clear.

{\rm (iii)} It follows from Part {\rm (ii)} and \cite [Th. 5.1.5] {Lo83}. \hfill $\Box$

Set $h_u := {\rm inf } \{ m \mid u ^m =0 \} (i.e. = {\rm min }  \{ m \mid u ^m =0 \}$ when there exists $m \in \mathbb N$ such that $u^m=0$), which is called a nil order of $u.$

\begin {Proposition}\label {ppppp7}
  $\{ [u] ^- \mid u \in {\rm SLLW} (T(V)/I)\}$ is a basis of   $\mathfrak L ^- (T(V)/I)$.
\end {Proposition}
\noindent {\it Proof.} Let $J^- := I \cap {\rm Lie} (V).$  By Proposition \ref {p7.10},
$\mathfrak L^-(T(V)/I)  \cong {\rm Lie} (V) / J^-$. By \cite [Th. 2.1] {LR95},
$\{ [u] ^- \mid u \in {\rm SLLW} ({\rm Lie} (V) / J^-)\}$ is a basis of ${\rm Lie} (V) / J^-$. Consequently, $\{ [u] ^- \mid u \in {\rm SLLW} (T(V)/I)\}$ is a basis of   $\mathfrak L ^- (T(V)/I)$.
 \hfill $\Box$
\begin{Proposition} \label{pp2.3.4} Assume that  $I$ is a graded  ideal of $T(V)$ as algebras (e.g. ${\mathcal E}_n = T(V)/ I$). Then

{\rm (i)} $\dim (T(V)/I) \le \prod \limits_{u \in {\rm SW } \cap {\rm L}} h_u$.
 In particular, $\dim ({\mathcal E_n}) \le \prod \limits_{u \in {\rm SW } \cap {\rm L}} h_u$.

{\rm (ii)} The following conditions are equivalent:

{\rm (F1)} $\dim ( T(V)/I )< \infty$;

{\rm (F11)} ${\rm SW} ( T(V)/I )$ is a finite set;

{\rm (F12)} ${\rm SW} ( T(V)/I ) \cap {\rm L}$ is a finite set and $h _u < \infty$ for any $u \in  {\rm SW} ( T(V)/I ) \cap {\rm L}$.

{\rm (F13)}  $\{ [u] ^- \mid u \in {\rm SLLW} (T(V)/I)\}$ is finite and  $h _u < \infty$ for any $u \in  {\rm SLLW} ( T(V)/I ) $.

\end{Proposition}
\noindent {\it Proof.}  {\rm (i)} It follows from Lemma \ref {pp10.7''}.

{\rm (ii)} It follows from Lemma \ref {pp10.7''} that {\rm (F1)} and {\rm (F2)} are equivalent. By Part {\rm (i)},  {\rm (F12)} implies {\rm (F11)}.

Now we show that {\rm (F11)} implies {\rm (F12)}. If it does not hold, then ${\rm SW} ( T(V)/I ) \cap {\rm L}$ is infinite or there exists $u \in {\rm SW} ( T(V)/I ) \cap {\rm L}$ such that $h_u = \infty.$ ${\rm SW} ( T(V)/I ) \cap {\rm L}$ is finite by {\rm (F11)}. Consequently, $u, u^2, u^3, \cdots $ are linearly independent, which contradicts {\rm (F1)}.

  It follows from Proposition \ref {ppppp7} and
Theorem \ref {pp3.8''000} ({\rm F4})   that  {\rm (F1)} and {\rm (F13)} are equivalent. \hfill $\Box$

This also give  an estimate of dimension  for every finite generated algebra $B$ since $B \cong T(V)/I$ as algebras.

\subsection { Lyndon basis of braided Hopf algebras with  diagonal type}

Let $(T(V)/I, C)$ be  a graded braided algebras with diagonal braiding $C$ and canonical basis $ \{ x_1, x_2, \cdots, x_\theta \}$
in this subsection without special announcement.

$[u] $ is called a super-letter or braided   Lie  Lyndon word  if $u$ is a Lyndon word. $[u]^- $ is called a Lie  Lyndon word  if $u$ is a Lyndon word.

  If $(V, C)$ is a braided vector space with diagonal type, then there a commutative group $G$ such that $(V, C)$  become a YD-module over $G$. Consequently, if $(V, C)$ is a finite dimensional diagonal braided vector space, then $\mathfrak B(V)$ is a
 graded  braided Hopf algebras over  commutative group $G$. 

 A Lyndon word  $u$ is said to be a hard braided Lie Lyndon word with respect to $T(V)/I$ if $[u] $ is not a linear
combination   in  $T(V)/I$ of products $[u_1][u_2]\cdot\cdot\cdot[u_i],  i\in
\mathbb N$,  where $u_j$ is a Lyndon word  with $ u <u_j$,  $1\le j \le i$. In this case $[u]$ is called a hard super-letter sometimes (cf. \cite {Kh99}).

 A Lyndon word  $u$ is said to be a hard  Lie Lyndon word with respect to $T(V)/I$ if $[u]^- $ is not a linear
combination in $T(V)/I$ of products $[u_1]^-[u_2]^-\cdot\cdot\cdot[u_i]^-,  i\in
\mathbb N$,  where $u_j$ is a Lyndon word  with $ u <u_j$,  $1\le j \le i$.

  Lyndon  word $u$ is called a standard braided Lie  Lyndon  word with respect to $T(V)/I$ if $[u] $ can not be  written as  linear combination of  strictly greater braided Lie  Lyndon words in $T(V)/I$.

\begin{Remark} The length of every  summand above  is equal since $T(V)/I$ is graded.
\end{Remark}

Let
$ {\rm HBLLW}(T(V)/I) =: \{ u  \in {\rm L}\mid  u \hbox { is a hard braided Lie Lyndon word }\}$,

  ${\rm SBLLW} (T(V)/I)=: \{ u \in {\rm L} \mid  u \hbox { is a standard  braided Lie Lyndon word} \}$.

${\rm HLLW }(T(V)/I) =: \{ u\in {\rm L} \mid  u \hbox { is a hard Lie Lyndon word  }\}$,

Let {\rm RPBW } denote a  restricted {\rm  PBW } basis of $T(V)/I$ (defined in \cite {Kh99}).

\begin {Lemma}\label {ppp1.1}
{\rm (i)}  If  $l \in {\rm L}, $  then   $ [l]^- =  l + \sum \limits _ {w > l, \mid l \mid = \mid w\mid} a_w w$ in $T(V)/I,$  where $ a_w \in k$.

{\rm (ii)} If $l \in  {\rm L}, $ then  $ [l] =  a_l l + \sum \limits _ {w > l, \mid l \mid = \mid w\mid} a_w w$ in $T(V)/I,$ where $a_l, a_w \in k$ with $a_l\not=0.$

{\rm (iii)}   $  {\rm  SW }\cap  {\rm L} \subseteq  {\rm SBLLW}.$

{\rm (iv)}   $ {\rm SW} \cap  {\rm L} \subseteq  {\rm HBLLW}.$

{\rm (v)}   $ {\rm SW} \cap  {\rm L} \subseteq  {\rm SLLW}.$

{\rm (vi)}   $ {\rm SW} \cap  {\rm L } \subseteq  {\rm HLLW}.$

{\rm (vii)} $\dim \mathfrak L^- (T(V)/ I) \ge \mid  {\rm HLLW} (T(V)/I)\mid, \mid  {\rm SLLW} (T(V)/I)\mid, \mid  {\rm SW}(T(V)/I) \cap  {\rm L}\mid.$
\end {Lemma}
\noindent {\it Proof.}  {\rm (ii)} We show this by induction on $\mid l \mid.$  It is clear when $\mid l \mid =1$ since $[l] = l.$ Assume that $l = uv$ is the Shirshow decomposition of $l$.
 If $u' > u$ and  $v'>v$ with $\mid u'\mid = \mid u\mid $
  and $\mid v'\mid = \mid v\mid $, then $u'v' > uv =l$ and $v'u' >vu >l.$
 \begin {eqnarray*} [l] &=& [v] [u] - p _{v, u}[u][v]\\
 &=& (a_v' v + \sum \limits _ {v' > v, \mid v' \mid = \mid v \mid} a_{v'}' v')(a_u 'u + \sum \limits _ {u' > u, \mid u' \mid = \mid u \mid} a_{u'}' u') \\
&& - p _{vu}(a_u 'u + \sum \limits _ {u' > u, \mid u' \mid = \mid u \mid} a_{u'}' u')
 (a_v' v + \sum \limits _ {v' > v, \mid v' \mid = \mid v \mid} a_{v'}' v')\ \ ( \hbox {by inductive  hypothesis}) \\
&=&  a_ll + \sum \limits _ {w > l, \mid l \mid = \mid w\mid} a_ww.
\end {eqnarray*}

{\rm (i)}  The proof is similar to the proof of {\rm (ii)}.

{\rm (iii)} If $l  \notin {\rm  SBLLW}$ with $l \in {\rm L}$, then
 \begin {eqnarray*}  [l] &=&  \sum \limits _ {u>l, \mid u \mid =  \mid l \mid } a_{u}' [u] \\
&=&  \sum \limits _ {u>l, \mid u \mid =  \mid l \mid}  a_u '( a''_u u + \sum\limits _{v > u, \mid v \mid =  \mid u \mid} a''_v v) \ \ \  \hbox {and  }\\
&& [l]  =  a_l''' l + \sum \limits _ {w > l, \mid l \mid = \mid w\mid} a_w '''w\ \ (  \hbox { by } {\rm (ii)}) .
\end {eqnarray*}
Consequently, $ l =  \sum \limits _ {w > l, \mid l \mid = \mid w\mid} a_w w$ and  $l \notin  {\rm SW} \cap  {\rm L}.$

{\rm (iv)} If $l  \notin {\rm  HBLLW}$ with $l \in {\rm L}$, then
 \begin {eqnarray*}  [l] &=& \sum \limits _{r=1} ^m  \sum \limits _ {u_i>l, 1\le  i \le r, \mid u_1u_2\cdots u_r \mid =  \mid l \mid  } a_{u}' [u_1] [u_2] \cdots [u_r] \\
&=&  \sum \limits _{r=1} ^m \sum \limits _ {u_i>l,  1\le  i \le r, \mid u_1u_2\cdots u_r \mid =  \mid l \mid } a_{u}' \prod _{i=1}^r ( a_{u_i}''u_i   + \sum \limits_{w_{ij} > u_i, \mid w_{ij} \mid =  \mid u_i \mid} a'' _{w_{ij}}  w_{ij}) \\
&=&\sum \limits _ {u>l, \mid u \mid =  \mid l \mid}   a_u u
\ \ \  \hbox { and  }\\
&& [l]  =  a_l''' l + \sum \limits _ {w > l, \mid l \mid = \mid w\mid} a_w '''w\ \ (  \hbox { by } {\rm (ii)}) .
\end {eqnarray*}
Consequently, $ l =  \sum \limits _ {w > l, \mid l \mid = \mid w\mid} a_w w$ and  $l \notin  {\rm SW} \cap  {\rm L}.$.

Similarly, we have {\rm (v)} and  {\rm (vi)}.

{\rm (vii)}  Obviously, both $\{ [u]^- \mid  u \in  {\rm HLLW } (T(V)/I)  \}$ and $\{ [u]^- \mid  u \in  {\rm SLLW } (T(V)/I)  \}$  are two linearly independent sets in ${\mathfrak L}^- (T(V)/I)$, which implies that $ \{ { [u]^- \mid  u \in  {\rm SW } (T(V)/I)  } \cap {\rm L} \}$  is linearly independent sets in ${\mathfrak L}^- (T(V)/I)$.
\hfill $\Box$

\begin {Lemma}\label {pp9.1}  $ {\rm HBLLW} (T(V)/ I) \subseteq {\rm SBLLW} (T(V)/ I)$
and $ {\rm HLLW} (T(V)/ I) \subseteq {\rm SLLW} (T(V)/ I)$.
\end {Lemma}
\noindent {\it Proof.}  If $u\notin {\rm SBLLW}, $ then  $[u] = \sum \limits_{v > u, v\in{\rm  L}, \mid u \mid =  \mid v \mid } a_v [v]$. Consequently, $u \notin {\rm HBLLW}.$  \hfill $\Box$

\begin {Lemma}\label {pp10.7'} Assume that   $T(V)/I$ is  a graded  braided Hopf algebras with diagonal braiding.
If both $B$ and $B_1$ are generators of RPBW basis   of $T(V)/I$  with $B \subseteq B_1$, then $B=B_1.$

\end {Lemma}
\noindent {\it Proof.}  If there exists $y \in B_1$ and $y \notin B, $  then $y = \sum \limits _{v_1 \ge v_2 \ge \cdots \ge v_r, v_i \in B} a_v v_1v_2 \cdots v_r$, which contradicts to that $B_1$ is a generator of  {\rm RPBW } basis   of $T(V)/I$.
\hfill $\Box$

\begin {Theorem}\label {pp10.7'''}  {\rm (i)} If  $T(V)/I$ is  a graded  braided Hopf algebras with diagonal braiding, then
${\rm SW } (T(V)/ I)\cap {\rm L} = {\rm HBLLW} (T(V)/I).$

Furthermore, if $\Delta ( \mathfrak B(V)  )$ is not  an arithmetic root system, then $\dim \mathfrak L^-( \mathfrak B(V)) = \infty.$

{\rm (ii)}
If $V$ is of  a   Cartan type  with $1 < {\rm ord } q_{ii} < \infty$ for $1\le i \le n$,  then the following conditions   are equivalent:

{\rm (a)} $\dim \mathfrak B(V)< \infty$;

{\rm  (b)}  $\dim \mathfrak L^-(V)<\infty$;

 {\rm (c)}  $\Delta ( \mathfrak B(V)  )$ is  an arithmetic root system;

 {\rm (d)} $V$ is of a finite   Cartan type.

\end {Theorem}
\noindent {\it Proof.} {\rm (i)} ${\rm SW }\cap {\rm L} \subseteq {\rm HBLLW}$ by Lemma \ref {ppp1.1} {\rm (iv)}. If there exists $y \in  {\rm HBLLW}$ and $y \notin {\rm SW}\cap {\rm L}, $  then $y = \sum \limits_{v_1 \ge v_2 \ge \cdots \ge v_r, v_i \in {\rm SW }\cap {\rm L},
\mid v_1v_2\cdots v_r \mid = \mid y \mid, v_1v_2\cdots v_r>y} a_v v_1v_2 \cdots v_r$ by Lemma \ref {pp10.7''} {\rm (iii)}, which contradicts to that ${\rm HBLLW}$ is  generator of RPBW basis   of $T(V)/I$.

{\rm (ii)}
${\rm (a)}
\Rightarrow {\rm (d)}. $  It is clear. $ {\rm (d)}
\Rightarrow {\rm (a)}$.
It follows from Theorem \ref {pp3.8''}. ${\rm (a)}
\Rightarrow {\rm (b)}$.  It is clear. ${\rm (b)}
\Rightarrow {\rm (c)}$.  It follows from {\rm (i)}.  ${\rm (c)}
\Rightarrow {\rm (d)}$.  It is clear.
\hfill $\Box$

\subsection { Lyndon basis of Nichols algebras with finite Cartan types}

Let $V$ be a braided vector space with finite Cartan type, which corresponds to Cartan matrix $A_c$ and  a  positive root system $\Phi^+$.
\begin {Lemma}\label {pp10.7''''}
For any $\alpha \in \Phi ^+$, there exists unique  $l \in {\rm SW} \cap {\rm L}$ such that $deg (l) = \alpha.$

\end {Lemma}
\noindent {\it Proof.} It follows from \cite [Th. 4] {He06} and Theorem \ref {pp10.7'''}. \hfill $\Box$

\begin{Proposition} \label {ppp104}(\cite[Proposition 2.9]{LR95})
If $l$ is a standard Lyndon ward, then
  $l$ is of the form $l=l_1
\ldots l_k a$, where
\begin{itemize}
    \item $l_i$ is a standard  Lyndon word, for each $i =1, \ldots , k$;
    \item $l_i$ is a beginning of $l_{i-1}$, for each $i>1$;
    \item $a$ is a letter.
\end{itemize}

\end{Proposition}

In what follows, we describe a set of  hard super- letters  for each connected  finite Cartan
 type.

In Figure 1 the trees give the  standard Lyndon words. Every tree in Figure 1 is called a Lyndon-basis-tree.
The root of each tree is the leftmost vertex. Every path from the root to a  vertex which is a white vertex in stand of a black vertex is called a Lyndon-basis-path. Let ${\rm P}(A_c) := \{ u \mid u \hbox { is a Lyndon-basis-path }\}$, which is written as ${\rm P}(V)$ sometime.

\begin{Theorem} \label{ppp101} If $V$ is a connected finite Cartan type with Cartan matrix $A_c$, then
$$ {\rm  HBLLW}(\mathfrak B(V)) = {\rm SW} (\mathfrak B(V)) \cap {\rm L} = {\rm P} (V) \hbox { \ and \  } \dim \mathfrak L^- (V) \ge \mid \Phi ^+\mid .$$ Furthermore, $u$ is a hard braided Lie Lyndon word if and only if $u$ is a   standard Lyndon word.

\end {Theorem}
\noindent {\it Proof.} It follows from the proofs in \cite [Section 4.2]{AA08} that  ${\rm SW } (\mathfrak B(V))\cap {\rm L} \subseteq  {\rm P}(V)$.\emph{}

 By \cite [Th. 4]{He06}, $\mid  {\rm HBLLW} \mid = \mid \Phi^+\mid. $ By\cite [Th. 3.4]{LR95}, $\mid {\rm P}(V)\mid = \mid \Phi^+\mid. $   By Theorem \ref {pp10.7'''}, $ {\rm SW }\cap {\rm L }= {\rm HBLLW}. $ Consequently, considering Lemma \ref {ppp1.1},  we complete the proof.
\hfill $\Box$

\begin{Remark}
In \cite [before Lemma 3.9]{AA07} Nicol¨¢s Andruskiewitsch and  Iv¨¢n Ezequiel Angiono said   ``It is
easy to see that a braiding of Cartan type is standard, see the first part of
the proof of [H1, Th. 1]". Therefore, a braiding of finite Cartan type is standard and       ${\rm SW}(\mathfrak B(V)) \cap {\rm L} = {\rm P}(V)$ holds by \cite [subsection 4.2] {AA08}.
\end{Remark}

For example, for $C_n,$  all of  standard Lyndon word  are    $x_i x_{i+1} \cdots x_j$, $1\le i \le j \le n-1; $
$x_i x_{i+1} \cdots x_{n-1} x_i \cdots x_j $, $1\le i \le n-1, i\le j\le n; $
$x_i x_{i+1} \cdots x_{n-1} x_nx_{n-1} \cdots x_j $, $1\le i \le n-1, i+1\le j\le n. $
It will be proved that all of  hard  super- letters  are    $[x_i x_{i+1} \cdots x_j] $ , $1\le i \le j  \le n-1; $
$[x_i x_{i+1} \cdots x_{n-1} x_i \cdots x_j ]$, $1\le i \le n-1, i\le j\le n; $
$[x_i x_{i+1} \cdots x_{n-1} x_nx_{n-1} \cdots x_j] $, $1\le i \le n-1, i+1\le j\le n. $
By \cite {Bo89} all of  positive roots   are    $e_i + e_{i+1} +\cdots +e_j $, $1\le i \le j \le n-1; $
$e_i+ e_{i+1} +\cdots+ e_{n-1} +e_i +\cdots +e_j $, $1\le i \le n-1, i\le j\le n; $
$e_i +e_{i+1} +\cdots+ e_{n-1}+ e_n+e_{n-1} +\cdots+ e_j $ , $1\le i \le n-1, i+1\le j\le n. $

\subsection {Lyndon basis of enveloping algebras of  Lie algebras}

Let ${\rm Lie }  (A)$ be the Lie algebra generated by $A := \{ x_1, x_2, \cdots, x_\theta\}$ in $T(V)$ with $V := {\rm span} (A)$, which is a basis of $V.$  Obviously, ${\rm Lie }  (A)={\rm Lie }  (V)$.  Let $J$ be a graded Lie ideal of ${\rm Lie }  (V)$ and $I$ the ideal of $T(V)$ by $J.$
 Assume that  $L : = {\rm Lie }  (A)/J $ is a Lie algebra. By \cite [Section 2]{LR95}, $U(L) = T(V) / I$ is the  enveloping algebra of $L$ with Lie operation $[x, y]^- = xy -yx.$

\begin{Proposition} \label{ppp111'''}  If  $J$ is a graded Lie ideal of ${\rm Lie }  (V)$ and $I$ is an ideal of $T(V)$ generated by $J$, then
$${\rm SLLW} (T(V)/I) = {\rm HLLW} (T(V)/I) = {\rm SLLW} ({\rm Lie}(V)/J) ={\rm SW} (T(V)/I) \cap {\rm L},$$ where ${\rm SLLW} ({\rm Lie}(V)/J):= {\rm SLLW} ( T(V)/J)$.

Furthermore, if $(T(V)/I, C)$ be  a braided graded algebras, then  the sets above is equal to
$ {\rm HBLLW} (T(V)/I).$
.
\end {Proposition}

\noindent {\it Proof.} We show this by following steps. {\rm (i)} If $u \notin {\rm SLLW} ({\rm Lie}(V)/J)$, then
\begin {eqnarray*} & [u] ^- \equiv &  \sum \limits _ {w > u, \mid u \mid = \mid w\mid, w\in {\rm L}} a_w [w]^- \ \  ( \hbox {mod } J) \hbox { and }\\
& [u] ^- \equiv& \sum \limits _ {w > u, \mid u \mid = \mid w\mid, w\in {\rm L}} a_w [w]^- \ \  ( \hbox {mod } I),
\end {eqnarray*} where $ a_w \in k$.
Consequently, $u \notin {\rm HLLW} (T(V)/I)$ and $u \notin {\rm SLLW} (T(V)/I)$, which implies ${\rm HLLW} (T(V)/I) \subseteq  {\rm SLLW} ({\rm Lie}(V)/J) $ and ${\rm SLLW} (T(V)/I) \subseteq  {\rm SLLW} ({\rm Lie }  (V)/J) $.

{\rm (ii)}  It is clear that $U(L)$ is a pointed Hopf algebra generated primitive elements. Thus ${\rm HLLW}(T(V)/ I)$  is   generators  of {\rm PBW} basis of $U(L)$ by \cite [Th. 2, Cor. 1] {Kh99}. By \cite [Th. 2.7]{LR95}, ${\rm SLLW} ( {\rm Lie }  (A)/ J)$  is generators  of {\rm PBW} basis of $U(L).$

{\rm (iii) }  Considering Part {(i)} {\rm (ii)} and Lemma \ref {pp10.7'}, we have ${ \rm HLLW}(T(V)/ I) = {\rm SLLW} ({\rm Lie }  (V)/ J)$.

{\rm (iv) } By \cite [Cor. 2.8]{LR95}, $
 {\rm  SW } (T(V)/I) \cap L = {\rm SLLW} ({\rm Lie }  (V)/ J)$.

{\rm (v) } Considering Lemma \ref {ppp1.1}, we complete the proof.  \hfill $\Box$

\subsection { Lyndon basis of  enveloping  algebras of  simple Lie algebras}

Let $L$ be a simple Lie algebra with Cartan matrix $A_c = (a_{ij})_{n \times n}$ and following Serre relations:

(S1) $[h_i,    h_j ] ^- = 0,   $ $1\le i,    j \le n. $

(S2) $[x_i,    y_i]^- = h_i,   $  $[x_i,    y_j]^-=0$,     $i \not= j $,    $1\le
i,    j \le n. $

(S3) $[h_i,    x_j]^- = a_{ij} x_j$,    $[h_i,    y_j] ^-= -a_{ij} y_j$,   $1\le
i,    j \le n. $

$(S_4 ^+ )$  $({\rm ad } x_i)^{-a_{ij} +1} x_j =0$,    $i \not= j. $

$(S_4 ^- )$  $({\rm
 ad } y_i)^{-a_{ij} +1} y_j =0$,    $i \not= j. $

Let $ X:= \{ x_1, x_2, \cdots, x_n\}$, $Y:= \{ y_1, y_2, \cdots, y_n\}$, $\eta := \{ h_1, h_2, \cdots, h_n\}$, $L^+ := {\rm Lie} (X)/J^+ $, $L^- := {\rm Lie} (Y)/J^- $,
 where  $J^+$ is a  Lie ideal of   ${\rm Lie }  (X)$ generated by $S_{4}^+$.

\begin{Theorem} \label{ppp111}  Assume that  $L$ is a simple Lie algebra corresponding with Cartan matrix $A_c$. then
$${\rm SLLW} (U(L^+)) ={ \rm HLLW}(U(L^+)) = {\rm SLLW} ( {\rm Lie }  (X) / J^+) = {\rm SW}(U(L^+))\cap {\rm L} = {\rm P}(A_c).$$
Furthermore, if  $V$ is a connected finite Cartan type with Cartan matrix $A_c$, then
the sets above are equal to $$ {\rm  HBLLW}(\mathfrak B(V)) = {\rm SW} (\mathfrak B(V)) \cap {\rm L} = {\rm P} (V).$$

\end {Theorem}
\noindent {\it Proof.}  By \cite [Th. 3.4]{LR95}, $SW(U(L^+))\cap L ={\rm P}(A_c).$
The others follows from Proposition \ref {ppp111'''} and  Theorem \ref {ppp101}. \hfill $\Box$

\appendix

\section*{Appendix}

\section {Proof of Theorem \ref {1.5}.}

 \noindent {\it Proof.} {\rm (i)}  There exists subrack $X \subseteq \mathcal{O}_{\sigma}^G$ such that
$R \cup S= X$ is a subrack  decomposition of $X$, and there exist $a\in R$, $b\in S$ such that ${\rm sq} (a, b) := a\rhd (b \rhd (a \rhd b)) \not= b$. Let $H$ be the subgroup generated by $X$ in  $G$. It is clear ${\mathcal O}_{\sigma _1}^H \not= {\mathcal O}_{\sigma _2} ^H$ for any $\sigma _1 \in R, \sigma _2 \in S. $ In fact, if ${\mathcal O}_{\sigma _1}^H = {\mathcal O}_{\sigma _2} ^H$, then there exists $t\in H$, such that $\sigma _2 = t \sigma _1t^{-1}$. $t = \prod _{i=1} ^n t_i$ whit $t_i \in X$ or $t_i^{-1}\in X.$ Therefore $\sigma_2 \in R$, which is a contradiction.

{\rm (ii)} By {\rm (i)}, we have ${\mathcal O}_{a}^H \not= {\mathcal O}_{b} ^H$  and they are not square commutative.  By \cite [Th. 8.6] {HS08}, $\dim \mathfrak B (M({\mathcal O}_{a}^H, \lambda  )+ M({\mathcal O}_{b} ^H, \xi )) = \infty,  $ for any
$\lambda \in \widehat {H^a}$, $\xi  \in \widehat {H^b}$. Consequently,  {\rm dim} $\mathfrak B(\mathcal{O}_{\sigma}^G,  \rho) = \infty $ for any $\rho \in \widehat { G^\sigma}$ by  \cite [Lemma 3.2{ \rm ii}]{AFGV08}. \hfill $\Box$

\section {Lyndon-basis-tree}

\vskip 1cm
$$\hbox {Figure 1}$$\\

 $A_n:$\begin{picture}(100,     15) \put(27,     1){\makebox(0,     0)[t]{$\circ $}}
\put(60,     1){\makebox(0,     0)[t]{$\circ $}}  \put(126,     1){\makebox(0,
0)[t]{$\circ $}} \put(159,     1){\makebox(0,     0)[t]{$\circ $}}
\put(28,     -1){\line(1,     0){30}}
 \put(104,    -1){\line(1,     0){20}} \put(127,     -1){\line(1,     0){30}}
\put(22,    15){i} \put(58,     15){i+1}  \put(124,
15){n-1} \put(157,     15){n}
\end{picture}\\ \\

 $B_n:$\begin{picture}(100,     15) \put(27,     1){\makebox(0,     0)[t]{$\circ $}}
\put(60,     1){\makebox(0,     0)[t]{$\circ $}}  \put(126,     1){\makebox(0,
0)[t]{$\circ $}} \put(159,     1){\makebox(0,     0)[t]{$\circ $}}
\put(192,     1){\makebox(0,     0)[t]{$\circ $}}\put(225,     1){\makebox(0,     0)[t]{$\circ $}}\put(290,     1){\makebox(0,     0)[t]{$\circ $}}\put(323,     1){\makebox(0,     0)[t]{$\circ $}}
\put(28,     -1){\line(1,     0){30}}
 \put(104,    -1){\line(1,     0){20}} \put(127,     -1){\line(1,     0){30}}
 \put(160,     -1){\line(1,     0){30}} \put(193,     -1){\line(1,     0){30}}
 \put(291,    -1){\line(1,     0){30}}
\put(22,    15){i} \put(56,     15){i+1}  \put(122, 15){n-1} \put(157,     15){n}
 \put(190,     15){n} \put(221,     15){n-1} \put(283,     15){i+2} \put(318,     15){i+1}
\end{picture}\\ \\

$C_n:$\begin{picture}(100,     15) \put(27,     1){\makebox(0,     0)[t]{$\circ $}}
\put(60,     1){\makebox(0,     0)[t]{$\circ $}}  \put(126,     1){\makebox(0,
0)[t]{$\circ $}} \put(159,     1){\makebox(0,     0)[t]{$\circ $}}\put(159,    -22){\makebox(0,     0)[t]{$\bullet $}}
\put(192,     1){\makebox(0,     0)[t]{$\circ $}}\put(192,    -22){\makebox(0,     0)[t]{$\bullet $}}
\put(253,     1){\makebox(0,     0)[t]{$\circ $}}\put(253,    -22){\makebox(0,     0)[t]{$\bullet $}}
\put(284,     1){\makebox(0,     0)[t]{$\circ $}}\put(284,    -22){\makebox(0,     0)[t]{$\bullet $}}
\put(316,     -22){\makebox(0,     0)[t]{$\circ $}}
\put(28,     -1){\line(1,     0){30}}\put(127,     -1){\line(1,     0){30}}
 \put(104,    -1){\line(1,     0){20}} \put(127,     -1){\line(4,     -3){30}}
 \put(160,     -1){\line(1,     0){30}} \put(193,     -1){\line(1,     0){30}}
 \put(254,    -1){\line(1,     0){29}} \put(160,     -23){\line(1,     0){30}}\put(193,     -23){\line(1,     0){30}}
 \put(254,    -23){\line(1,     0){29}} \put(286,    -23){\line(1,     0){29}}
\put(22,    15){i} \put(56,     15){i+1}  \put(122, 15){n-1} \put(157,     15){n} \put(190,     15){n-1}  \put(250,     15){i+2} \put(280,     15){i+1}
\put(157,     -15){i} \put(190,     -15){i+1}  \put(250,     -15){n-2} \put(280,     -15){n-1} \put(312,     -15){n}
\end{picture}\\ \\

$D_n:$\begin{picture}(100,     15) \put(27,     1){\makebox(0,     0)[t]{$\circ $}}
\put(60,     1){\makebox(0,     0)[t]{$\circ $}}  \put(126,     1){\makebox(0, 0)[t]{$\circ $}} \put(159,     1){\makebox(0,     0)[t]{$\circ $}}\put(159,    -24){\makebox(0,     0)[t]{$\circ $}}
\put(192,     1){\makebox(0,     0)[t]{$\circ $}}
\put(253,     1){\makebox(0,     0)[t]{$\circ $}}
\put(284,     1){\makebox(0,     0)[t]{$\circ $}}
\put(28,     -1){\line(1,     0){30}}\put(127,     -1){\line(1,     0){30}}
 \put(104,    -1){\line(1,     0){20}} \put(127,     -1){\line(4,     -3){30}}
 \put(160,     -1){\line(1,     0){30}} \put(193,     -1){\line(1,     0){30}}
 \put(254,    -1){\line(1,     0){29}}
\put(22,    5){i} \put(56,     5){i+1}  \put(122, 5){n-2} \put(157,     5){n} \put(190,     5){n-1}
\put(250,    5){i+2} \put(280,     5){i+1} \put(157,     -15){n-1}
\end{picture}\\ \\

$F_4:$ \begin{picture}(50,     15) \put(27,     1){\makebox(0,     0)[t]{$\circ $}}
\put(22,    10){4}
\end{picture} ,\!\!\!\!\!\!\!\!\!\!\!\!\!\!
\begin{picture}(100,     15) \put(83,     1){\makebox(0,     0)[t]{$\circ $}}
\put(50,     1){\makebox(0,     0)[t]{$\circ $}}
\put(51,     -1){\line(1,     0){30}}
\put(50,    10){3} \put(81,     10){4}
\end{picture},
\begin{picture}(190,     15) \put(83,     1){\makebox(0,     0)[t]{$\circ $}}
\put(50,     1){\makebox(0,     0)[t]{$\circ $}}\put(116,     1){\makebox(0,     0)[t]{$\circ $}}
\put(149,     1){\makebox(0,     0)[t]{$\circ $}}\put(182,     1){\makebox(0,     0)[t]{$\circ $}}
\put(116,    -20){\makebox(0,     0)[t]{$\circ $}}
\put(51,     -1){\line(1,     0){30}}\put(84,     -1){\line(1,     0){30}}\put(117,     -1){\line(1,     0){30}}\put(150,     -1){\line(1,     0){30}}\put(84,     -1){\line(4,     -3){30}}
\put(50,    10){2} \put(81,     10){3}\put(116,     10){4}\put(149,     10){3}\put(182,     10){4}\put(120,     -18){3}
\end{picture},\\ \\ \\

\begin{picture}(100,     15) \put(83,     1){\makebox(0,     0)[t]{$\circ $}}
\put(50,     1){\makebox(0,     0)[t]{$\circ $}}\put(116,     1){\makebox(0,     0)[t]{$\circ $}}
\put(149,     1){\makebox(0,     0)[t]{$\circ $}}\put(182,     1){\makebox(0,     0)[t]{$\circ $}}
\put(215,     1){\makebox(0,     0)[t]{$\circ $}}\put(248,     1){\makebox(0,     0)[t]{$\circ $}}
\put(281,     1){\makebox(0,     0)[t]{$\circ $}}\put(314,     1){\makebox(0,     0)[t]{$\circ $}}
\put(347,    1){\makebox(0,     0)[t]{$\circ $}}
\put(215,     24){\makebox(0,     0)[t]{$\bullet $}}\put(248,   24  ){\makebox(0,     0)[t]{$\bullet $}}
\put(281,     24){\makebox(0,     0)[t]{$\bullet $}}\put(314,    24){\makebox(0,     0)[t]{$\bullet $}}
\put(347,    24){\makebox(0,     0)[t]{$\bullet $}}\put(380,    24){\makebox(0,     0)[t]{$\circ $}}
\put(149,    -20){\makebox(0,     0)[t]{$\circ $}}\put(182,    -20){\makebox(0,     0)[t]{$\circ $}}
\put(215,    -20){\makebox(0,     0)[t]{$\circ $}}\put(248,    -20){\makebox(0,     0)[t]{$\circ $}}
\put(51,     -1){\line(1,     0){30}}\put(84,     -1){\line(1,     0){30}}\put(117,     -1){\line(1,     0){30}}
\put(150,     -1){\line(1,     0){30}}\put(150,     -1){\line(1,     0){30}}\put(183,     -1){\line(1,     0){30}}
\put(216,     -1){\line(1,     0){30}}\put(249,     -1){\line(1,     0){30}}\put(282,     -1){\line(1,     0){30}}
\put(315,     -1){\line(1,     0){30}}\put(183,     -1){\line(4,     3){31}}\put(216,     21){\line(1,     0){30}}
\put(249,    21){\line(1,     0){30}}\put(282,    21){\line(1,     0){30}}\put(315,     21){\line(1,     0){30}}
\put(348,     21){\line(1,     0){30}}
\put(117,     -1){\line(4,     -3){30}}\put(183,     -1){\line(4,     -3){30}}
\put(150,     -22){\line(1,     0){30}}\put(216,     -22){\line(1,     0){30}}
\put(50,    6){1} \put(81,     6){2}\put(116,     6){3}\put(149,     6){4}\put(182,     6){3}
\put(215,    6){4} \put(248,     6){2}\put(281,     6){3}\put(314,     6){3}\put(347,     6){2}
\put(149,     -18){3}\put(182,     -18){2}\put(215,     -18){2}\put(248,     -18){3}
\put(215,     28){1}\put(248,     28){2}\put(281,     28){3}\put(314,     28){4}\put(347,     28){3}
\put(380,     28){2}
\end{picture}\\ \\

$G_2:$ \begin{picture}(155,     15) \put(27,     1){\makebox(0,     0)[t]{$\circ $}}
\put(50,     1){\makebox(0,     0)[t]{$\bullet $}}\put(83,     1){\makebox(0,     0)[t]{$\bullet $}}
\put(116,     1){\makebox(0,     0)[t]{$\circ $}}\put(58,    -20){\makebox(0,     0)[t]{$\circ $}}
\put(83,    -20){\makebox(0,     0)[t]{$\circ $}}\put(116,     -20){\makebox(0,     0)[t]{$\bullet $}}
\put(149,     -20){\makebox(0,     0)[t]{$\circ $}}
\put(28,     -1){\line(1,     0){30}}\put(51,     -1){\line(1,     0){30}}\put(84,     -1){\line(1,     0){30}}
\put(84,     -22){\line(1,     0){30}}\put(28,     -1){\line(4,     -3){30}}\put(51,     -1){\line(4,     -3){30}}
\put(117,     -22){\line(1,     0){30}}
\put(27,    6){1} \put(50,     6){1}\put(83,     6){1}\put(116,     6){2}\put(48,     -15){2}\put(81,     -15){2}
\put(114,     -15){1}\put(147,     -15){2}
\end{picture} ,
\begin{picture}(100,     15) \put(50,     1){\makebox(0,     0)[t]{$\circ $}}
\put(50,    10){2}
\end{picture} . \\ \\

$E_8:$ \begin{picture}(40,     15) \put(27,     1){\makebox(0,     0)[t]{$\circ $}}
\put(22,    10){8}
\end{picture}, \!\!\!\!\!\!\!\!\!\!\!\!
\begin{picture}(100,     15) \put(83,     1){\makebox(0,     0)[t]{$\circ $}}
\put(50,     1){\makebox(0,     0)[t]{$\circ $}}
\put(51,     -1){\line(1,     0){30}}
\put(50,    10){7} \put(81,     10){8}
\end{picture},
\begin{picture}(123,     15) \put(83,     1){\makebox(0,     0)[t]{$\circ $}}
\put(50,     1){\makebox(0,     0)[t]{$\circ $}}\put(116,     1){\makebox(0,     0)[t]{$\circ $}}
\put(51,     -1){\line(1,     0){30}}\put(84,     -1){\line(1,     0){30}}
\put(50,    10){6} \put(81,     10){7}\put(114,     10){8}
\end{picture},\!\!\!\!\!\!\!
\begin{picture}(158,     15) \put(83,     1){\makebox(0,     0)[t]{$\circ $}}\put(50,     1){\makebox(0,     0)[t]{$\circ $}}\put(116,     1){\makebox(0,     0)[t]{$\circ $}}\put(149,     1){\makebox(0,     0)[t]{$\circ $}}
\put(51,     -1){\line(1,     0){30}}\put(84,     -1){\line(1,     0){30}}\put(117,     -1){\line(1,     0){30}}
\put(50,    10){5} \put(81,     10){6}\put(114,     10){7}\put(148,     10){8}
\end{picture},\\

\begin{picture}(195,     15) \put(83,     1){\makebox(0,     0)[t]{$\circ $}}\put(50,     1){\makebox(0,     0)[t]{$\circ $}}\put(116,     1){\makebox(0,     0)[t]{$\circ $}}\put(149,     1){\makebox(0,     0)[t]{$\circ $}}
\put(182,     1){\makebox(0,     0)[t]{$\circ $}}
\put(51,     -1){\line(1,     0){30}}\put(84,     -1){\line(1,     0){30}}\put(117,     -1){\line(1,     0){30}}
\put(150,     -1){\line(1,     0){30}}
\put(50,    10){4} \put(81,     10){5}\put(114,     10){6}\put(148,     10){7}\put(181,     10){8}
\end{picture},
\begin{picture}(197,     15) \put(83,     15){\makebox(0,     0)[t]{$\circ $}}\put(50,     15){\makebox(0,     0)[t]{$\circ $}}\put(116,     15){\makebox(0,     0)[t]{$\circ $}}\put(149,     15){\makebox(0,     0)[t]{$\circ $}}
\put(182,     15){\makebox(0,     0)[t]{$\circ $}}\put(22,     15){\makebox(0,     0)[t]{$\circ $}}
\put(51,     13){\line(1,     0){30}}\put(84,     13){\line(1,     0){30}}\put(117,     13){\line(1,     0){30}}
\put(150,     13){\line(1,     0){30}}\put(23,     13){\line(1,     0){26}}
\put(50,    19){4} \put(81,     19){5}\put(114,     19){6}\put(148,     19){7}\put(181,     19){8}\put(22,     19){3}
\end{picture}, \\ \\ \\

\begin{picture}(100,     15) \put(22,     1){\makebox(0,     0)[t]{$\circ $}}\put(50,     -11){\makebox(0,     0)[t]{$\circ $}}\put(78,     -22){\makebox(0,     0)[t]{$\circ $}}\put(78,     1){\makebox(0,     0)[t]{$\circ $}}\put(106,     -10){\makebox(0,     0)[t]{$\circ $}}\put(106,     14){\makebox(0,     0)[t]{$\circ $}}
\put(134,     27){\makebox(0,     0)[t]{$\circ $}}\put(134,     1){\makebox(0,     0)[t]{$\circ $}}
\put(134,     -22){\makebox(0,     0)[t]{$\circ $}}
\put(162,     40){\makebox(0,     0)[t]{$\circ $}}\put(162,     14){\makebox(0,     0)[t]{$\circ $}}
\put(162,     -12){\makebox(0,     0)[t]{$\circ $}}
\put(190,     27){\makebox(0,     0)[t]{$\circ $}}\put(190,     1){\makebox(0,     0)[t]{$\circ $}}
\put(190,     -25){\makebox(0,     0)[t]{$\circ $}}
\put(218,     14){\makebox(0,     0)[t]{$\circ $}}\put(218,     -12){\makebox(0,     0)[t]{$\circ $}}
\put(246,     2){\makebox(0,     0)[t]{$\circ $}}\put(246,     -25){\makebox(0,     0)[t]{$\circ $}}
\put(274,     -12){\makebox(0,     0)[t]{$\circ $}}\put(302,     -25){\makebox(0,     0)[t]{$\circ $}}
\put(23,     -1){\line(2,     -1){25}}\put(51,     -13){\line(2,     -1){25}}\put(51,     -13){\line(2,     1){25}}
\put(79,     0){\line(2,     1){25}}\put(79,     0){\line(2,     -1){25}}
\put(107,     11){\line(2,     1){25}}\put(107,     11){\line(2,     -1){25}}\put(107,     -13){\line(2,     -1){25}}
\put(135,     24){\line(2,     1){25}}\put(135,     24){\line(2,     -1){25}}\put(135,     -1){\line(2,     -1){25}}
\put(163,     37){\line(2,     -1){25}}\put(163,     11){\line(2,     -1){25}}\put(163,     -15){\line(2,     -1){25}}
\put(191,     24){\line(2,     -1){25}}\put(191,     -2){\line(2,     -1){25}}
\put(219,     11){\line(2,     -1){25}}\put(219,     -15){\line(2,     -1){25}}
\put(248,     -2){\line(2,     -1){25}}\put(275,     -15){\line(2,     -1){25}}
\put(22,    6){2} \put(50,     -8){4} \put(78,     -18){3}\put(78,    2){5}\put(106,    15){6}\put(106,    -9){3}
\put(134,    28){7}\put(134,    2){3}\put(134,    -19){4}\put(162,    42){8} \put(162,    16){3}\put(162,    -10){4}
\put(190,    29){3} \put(190,    3){4}\put(190,    -23){5}
\put(218,    17){4}\put(218,    -10){5}\put(246,    4){5}\put(246,    -23){6}
 \put(274,    -10){6}\put(302,    -23){7}
\end{picture}\\ \\ \\

\begin{picture}(100, 15) \put(22, -103){\makebox(0, 0)[t]{$\circ $}}\put(42,  -93){\makebox(0, 0)[t]{$\circ $}}
\put(62,    -101){\makebox(0,    0)[t]{$\circ $}}\put(82,    -90){\makebox(0,    0)[t]{$\circ $}}
\put(82,    -110){\makebox(0,    0)[t]{$\circ $}}\put(102,    -80){\makebox(0,    0)[t]{$\circ $}}
\put(102,    -100){\makebox(0,    0)[t]{$\circ $}}\put(122,    -110){\makebox(0,    0)[t]{$\circ $}}
\put(122,    -70){\makebox(0,    0)[t]{$\circ $}}\put(122,    -90){\makebox(0,    0)[t]{$\circ $}}
\put(142,    -60){\makebox(0,    0)[t]{$\circ $}}\put(142,    -80){\makebox(0,    0)[t]{$\circ $}}
\put(142,    -100){\makebox(0,    0)[t]{$\circ $}}\put(142,    -120){\makebox(0,    0)[t]{$\circ $}}
\put(162,    -50){\makebox(0,    0)[t]{$\circ $}}\put(162,    -80){\makebox(0,    0)[t]{$\circ $}}
\put(162,    -90){\makebox(0,    0)[t]{$\circ $}}\put(162,    -110){\makebox(0,    0)[t]{$\circ $}}
\put(182,    -50){\makebox(0,    0)[t]{$\circ $}}\put(182,    -70){\makebox(0,    0)[t]{$\circ $}}
\put(182,    -90){\makebox(0,    0)[t]{$\circ $}}\put(182,    -100){\makebox(0,    0)[t]{$\circ $}}
\put(182,    -31){\makebox(0,    0)[t]{$\circ $}}\put(202,    -45){\makebox(0,    0)[t]{$\circ $}}
\put(202,    -60){\makebox(0,    0)[t]{$\circ $}}\put(202,    -78){\makebox(0,    0)[t]{$\circ $}}
\put(202,    -110){\makebox(0,    0)[t]{$\circ $}}\put(200,    -20){\makebox(0,    0)[t]{$\circ $}}\put(200,    -34){\makebox(0,    0)[t]{$\circ $}}
\put(216,    -38){\makebox(0,    0)[t]{$\circ $}}\put(232,    -46){\makebox(0,    0)[t]{$\circ $}}
\put(232,    -24){\makebox(0,    0)[t]{$\circ $}}\put(248,    -30){\makebox(0,    0)[t]{$\circ $}}
\put(248,    -15){\makebox(0,    0)[t]{$\circ $}}\put(264,    24){\makebox(0,    0)[t]{$\circ $}}\put(264,    6){\makebox(0,    0)[t]{$\circ $}}
\put(264,    -24){\makebox(0,    0)[t]{$\circ $}}\put(280,    15){\makebox(0,    0)[t]{$\circ $}}\put(280,    -3){\makebox(0,    0)[t]{$\circ $}}
\put(280,    -33){\makebox(0,    0)[t]{$\circ $}}\put(296,    6){\makebox(0,    0)[t]{$\circ $}}\put(296,    -12){\makebox(0,    0)[t]{$\circ $}}
\put(296,    -42){\makebox(0,    0)[t]{$\circ $}}\put(312,    15){\makebox(0,    0)[t]{$\circ $}}\put(312,    -3){\makebox(0,    0)[t]{$\circ $}}
\put(328,    6){\makebox(0,    0)[t]{$\circ $}}\put(344,    -3){\makebox(0,    0)[t]{$\circ $}}
\put(360,    -12){\makebox(0,    0)[t]{$\circ $}}\put(232,    -75){\makebox(0,    0)[t]{$\circ $}}\put(232,    -93){\makebox(0,    0)[t]{$\circ $}}
\put(248,    -66){\makebox(0,    0)[t]{$\circ $}}\put(248,   -84){\makebox(0,    0)[t]{$\circ $}}
\put(264,    -67){\makebox(0,    0)[t]{$\circ $}}\put(280,   -76){\makebox(0,    0)[t]{$\circ $}}
\put(296,   -85){\makebox(0,    0)[t]{$\circ $}}
\put(232,    -140){\makebox(0,    0)[t]{$\bullet $}}\put(248,    -150){\makebox(0,    0)[t]{$\bullet $}}
\put(264,    -160){\makebox(0,    0)[t]{$\bullet $}}\put(280,    -170){\makebox(0,    0)[t]{$\bullet $}}
\put(296,    -180){\makebox(0,    0)[t]{$\bullet $}}\put(312,    -190){\makebox(0,    0)[t]{$\circ $}}
\put(328,    -200){\makebox(0,    0)[t]{$\circ $}}
\put(248,    -130){\makebox(0,    0)[t]{$\bullet $}}\put(264,    -140){\makebox(0,    0)[t]{$\bullet $}}
\put(280,    -150){\makebox(0,    0)[t]{$\bullet $}}\put(296,    -160){\makebox(0,    0)[t]{$\bullet $}}
\put(312,    -170){\makebox(0,    0)[t]{$\bullet $}}\put(328,    -180){\makebox(0,    0)[t]{$\bullet $}}
\put(344,    -190){\makebox(0,    0)[t]{$\circ $}}
\put(264,    -120){\makebox(0,    0)[t]{$\bullet $}}\put(280,    -130){\makebox(0,    0)[t]{$\bullet $}}
\put(296,    -140){\makebox(0,    0)[t]{$\bullet $}}\put(312,    -150){\makebox(0,    0)[t]{$\bullet $}}
\put(328,    -160){\makebox(0,    0)[t]{$\bullet $}}\put(344,    -170){\makebox(0,    0)[t]{$\bullet $}}
\put(360,    -180){\makebox(0,    0)[t]{$\circ $}}
\put(280,    -110){\makebox(0,    0)[t]{$\bullet $}}\put(296,    -120){\makebox(0,    0)[t]{$\bullet $}}
\put(312,    -130){\makebox(0,    0)[t]{$\bullet $}}\put(328,    -140){\makebox(0,    0)[t]{$\bullet $}}
\put(344,    -150){\makebox(0,    0)[t]{$\bullet $}}\put(360,    -160){\makebox(0,    0)[t]{$\bullet $}}
\put(376,    -170){\makebox(0,    0)[t]{$\circ $}}
\put(225,    -145){\makebox(0,    0)[t]{$\nearrow $}}\put(240,    -154){\makebox(0,    0)[t]{$\nearrow $}}
\put(255,    -163){\makebox(0,    0)[t]{$\nearrow $}}\put(271,    -172){\makebox(0,    0)[t]{$\nearrow $}}
\put(288,    -184){\makebox(0,    0)[t]{$\nearrow $}}\put(304,    -195){\makebox(0,    0)[t]{$\nearrow $}}
\put(319,    -204){\makebox(0,    0)[t]{$\nearrow $}}
\put(296,    -100){\makebox(0,    0)[t]{$\bullet $}}\put(312,    -110){\makebox(0,    0)[t]{$\bullet $}}
\put(328,    -120){\makebox(0,    0)[t]{$\bullet $}}\put(344,    -130){\makebox(0,    0)[t]{$\bullet $}}
\put(360,    -140){\makebox(0,    0)[t]{$\bullet $}}\put(376,    -150){\makebox(0,    0)[t]{$\bullet $}}
\put(392,    -160){\makebox(0,    0)[t]{$\circ $}}\put(392,    -180){\makebox(0,    0)[t]{$\circ $}}
\put(312,    -90){\makebox(0,    0)[t]{$\bullet $}}\put(328,    -100){\makebox(0,    0)[t]{$\bullet $}}
\put(344,    -110){\makebox(0,    0)[t]{$\bullet $}}\put(360,    -120){\makebox(0,    0)[t]{$\bullet $}}
\put(376,    -130){\makebox(0,    0)[t]{$\bullet $}}\put(392,    -140){\makebox(0,    0)[t]{$\bullet $}}
\put(408,    -150){\makebox(0,    0)[t]{$\circ $}}
\put(408,    -170){\makebox(0,    0)[t]{$\circ $}}\put(408,    -190){\makebox(0,    0)[t]{$\circ $}}
\put(408,    -210){\makebox(0,    0)[t]{$\circ $}}\put(408,    -230){\makebox(0,    0)[t]{$\circ $}}
\put(408,    -250){\makebox(0,    0)[t]{$\circ $}}
\put(328,    -80){\makebox(0,    0)[t]{$\bullet $}}\put(344,    -90){\makebox(0,    0)[t]{$\bullet $}}
\put(360,    -100){\makebox(0,    0)[t]{$\bullet $}}\put(376,    -110){\makebox(0,    0)[t]{$\bullet $}}
\put(392,    -120){\makebox(0,    0)[t]{$\bullet $}}\put(408,    -130){\makebox(0,    0)[t]{$\bullet $}}
\put(424,    -140){\makebox(0,    0)[t]{$\circ $}}
\put(424,    -160){\makebox(0,    0)[t]{$\bullet $}}\put(424,    -180){\makebox(0,    0)[t]{$\bullet $}}
\put(424,    -200){\makebox(0,    0)[t]{$\bullet $}}\put(424,    -220){\makebox(0,    0)[t]{$\bullet $}}
\put(424,    -240){\makebox(0,    0)[t]{$\circ $}}
\put(435,    -150){\makebox(0,    0)[t]{$\swarrow $}}\put(435,    -170){\makebox(0,    0)[t]{$\swarrow $}}
\put(435,    -190){\makebox(0,    0)[t]{$\swarrow $}}\put(435,    -210){\makebox(0,    0)[t]{$\swarrow $}}
\put(435,    -230){\makebox(0,    0)[t]{$\swarrow $}}
\put(218,    -67){\makebox(0,    0)[t]{$\circ $}}\put(218,    -89){\makebox(0,    0)[t]{$\circ $}}
\put(218,    -118){\makebox(0,    0)[t]{$\circ $}}
\put(200,    -2){\makebox(0,    0)[t]{$\circ $}}\put(216,    -20){\makebox(0,    0)[t]{$\circ $}}
\put(216,    -2){\makebox(0,    0)[t]{$\circ $}}\put(232,    -2){\makebox(0,    0)[t]{$\circ $}}
\put(248,    15){\makebox(0,    0)[t]{$\circ $}}\put(248,    -2){\makebox(0,    0)[t]{$\circ $}}
\put(264,    24){\makebox(0,    0)[t]{$\circ $}}\put(264,    6){\makebox(0,    0)[t]{$\circ $}}
\put(24,    -105){\line(2,    1){16}}\put(44,    -94){\line(2,    -1){16}}\put(64,    -103){\line(2,    1){16}}
\put(64,    -103){\line(2,    -1){16}}\put(84,    -93){\line(2,    1){16}}\put(84,    -93){\line(2,    -1){16}}
\put(104,    -103){\line(2,    -1){16}}\put(104,    -83){\line(2,    -1){16}}\put(104,    -83){\line(2,    1){16}}
\put(124,    -73){\line(2,    1){16}}\put(124,    -73){\line(2,    -1){16}}\put(124,    -94){\line(2,    -1){16}}\put(124,    -113){\line(2,    -1){16}}
\put(144,    -63){\line(2,    1){16}}\put(144,    -83){\line(1,    0){16}}\put(144,    -103){\line(2,    1){16}}\put(144,    -103){\line(2,    -1){16}}
\put(164,    -52){\line(1,    0){16}}\put(164,    -83){\line(2,    1){16}}\put(164,    -83){\line(2,    -1){16}}\put(164,    -93){\line(2,    -1){16}}
\put(183,    -50){\line(0,    1){14}}\put(184,    -52){\line(4,    1){16}}\put(184,    -73){\line(2,    1){16}}
\put(184,    -73){\line(2,    -1){16}}\put(184,    -103){\line(2,    -1){16}}
\put(201,    -37){\line(4,    -1){14}}\put(217,    -42){\line(2,    -1){14}}
\put(217,    -22){\line(4,    -1){14}}\put(233,    -25){\line(2,    -1){14}}\put(233,    -25){\line(2,    1){14}}
\put(249,    12){\line(2,    1){14}}\put(249,    12){\line(2,    -1){14}}\put(249,    -18){\line(2,    -1){14}}
\put(265,    22){\line(2,    -1){14}}\put(265,    4){\line(2,    -1){14}}\put(265,    -26){\line(2,    -1){14}}
\put(281,    13){\line(2,    -1){14}}\put(281,    -5){\line(2,    -1){14}}\put(281,    -35){\line(2,    -1){14}}
\put(297,    4){\line(2,    1){14}}\put(297,    4){\line(2,    -1){14}}
\put(313,    13){\line(2,    -1){14}}\put(329,    4){\line(2,    -1){14}}\put(345,    -5){\line(2,    -1){14}}
\put(217,    -71){\line(3,    -1){14}}\put(218,    -90){\line(3,    -1){13}}\put(233,    -77){\line(2,    1){14}}\put(233,    -77){\line(2,    -1){14}}
\put(220,    -121){\line(2,    -3){14}}\put(233,    -142){\line(4,    -3){14}}\put(249,    -152){\line(4,    -3){14}}
\put(265,    -162){\line(4,    -3){14}}\put(281,    -172){\line(4,    -3){13}}\put(297,    -182){\line(4,    -3){14}}
\put(313,    -192){\line(4,    -3){13}}
\put(233,    -97){\line(2,    -5){14}}\put(249,    -132){\line(4,    -3){14}}\put(265,    -142){\line(4,    -3){14}}
\put(281,    -152){\line(4,    -3){14}}\put(297,    -162){\line(4,    -3){13}}\put(313,    -172){\line(4,    -3){14}}
\put(329,    -182){\line(4,    -3){13}}
\put(249,    -87){\line(2,    -5){14}}\put(265,    -122){\line(4,    -3){14}}\put(281,    -132){\line(4,    -3){14}}
\put(297,    -142){\line(4,    -3){14}}\put(313,    -152){\line(4,    -3){13}}\put(329,    -162){\line(4,    -3){14}}
\put(345,    -172){\line(4,    -3){13}}
\put(265,    -70){\line(1,    -3){14}}\put(281,    -112){\line(4,    -3){14}}\put(297,    -122){\line(4,    -3){14}}
\put(313,    -132){\line(4,    -3){14}}\put(329,    -142){\line(4,    -3){13}}\put(345,    -152){\line(4,    -3){14}}
\put(361,    -162){\line(4,    -3){14}}
\put(281,    -80){\line(2,    -3){14}}\put(297,    -102){\line(4,    -3){14}}\put(313,    -112){\line(4,    -3){14}}
\put(329,    -122){\line(4,    -3){14}}\put(345,    -132){\line(4,    -3){13}}\put(361,    -142){\line(4,    -3){14}}
\put(377,    -152){\line(4,    -3){13}}\put(392,    -164){\line(0,    -1){16}}
\put(297,    -88){\line(5,   -1){14}}\put(313,    -92){\line(4,    -3){14}}\put(329,    -102){\line(4,    -3){14}}
\put(345,    -112){\line(4,    -3){14}}\put(361,    -122){\line(4,    -3){13}}\put(377,    -132){\line(4,    -3){14}}
\put(393,    -142){\line(4,    -3){13}}
\put(408,    -154){\line(0,   -1){16}}\put(408,    -174){\line(0,    -1){16}}\put(408,    -194){\line(0,    -1){16}}
\put(408,    -214){\line(0,    -1){16}}\put(408,    -234){\line(0,    -1){16}}
\put(297,    -45){\line(4,   -5){30}}\put(329,    -82){\line(4,    -3){14}}\put(345,    -92){\line(4,    -3){14}}
\put(361,    -102){\line(4,    -3){14}}\put(377,    -112){\line(4,    -3){13}}\put(393,    -122){\line(4,    -3){14}}
\put(409,    -132){\line(4,    -3){13}}
\put(424,    -144){\line(0,   -1){16}}\put(424,    -164){\line(0,    -1){16}}\put(424,    -184){\line(0,    -1){16}}
\put(424,    -204){\line(0,    -1){16}}\put(424,    -224){\line(0,    -1){16}}
\put(184,    -32){\line(4,    3){14}}\put(184,    -32){\line(4,    -1){14}}\put(203,    -64){\line(2,    -1){14}}
\put(203,    -83){\line(2,    -1){14}}\put(203,    -112){\line(2,    -1){14}}
\put(201,    -22){\line(1,    0){14}}\put(200,    -20){\line(0,    1){14}}
\put(201,    -4){\line(1,    0){14}}\put(217,    -4){\line(1,    0){14}}
\put(233,    -4){\line(1,    0){14}}\put(233,    -4){\line(3,   4){14}}
\put(249,    -70){\line(6,    1){14}}\put(265,    -70){\line(2,    -1){14}}\put(281,    -79){\line(2,   -1){14}}
\put(22,   -100){1} \put(42,    -90){3}\put(62,    -99){4}\put(82,    -87){5}\put(82,    -110){2}\put(102,    -78){6}\put(102,    -97){2}\put(122,    -68){7}\put(122,    -90){2}\put(122,    -108){4}\put(140,    -58){8}\put(140,    -78){2}\put(140,    -99){4}\put(140,    -118){3}\put(160,    -48){2}\put(160,    -78){4}\put(160,    -104){5}\put(160,    -123){3}\put(176,    -48){4}\put(180,    -68){5}\put(180,    -88){3}\put(179,    -114){3}
\put(180,    -28){5}\put(195,    -45){3}\put(197,    -58){6}\put(197,    -78){3}\put(197,    -108){4}
\put(216,    -40){4}\put(232,    -50){2}\put(226,    -22){4}\put(238,    -40){2}\put(238,    -17){5}
\put(264,    24){6}\put(258,    -6){2}\put(264,    -22){2}\put(280,    17){2}\put(272,    -8){4}\put(272,    -39){4}
\put(296,    7){4}\put(293,    -11){3}\put(293,    -41){3}\put(312,    18){5}\put(315,    -4){3}
\put(328,    8){3}\put(344,    -1){4}\put(360,    -10){2}\put(232,    -73){4}\put(230,    -91){2}\put(248,    -64){5}\put(245,    -83){2}
\put(440,    -148){2}\put(440,    -168){4}\put(440,    -188){5}\put(440,    -208){6}\put(440,    -228){7}

\put(213,   -167){1 }\put(225,    -172){3}\put(240,    -180){4}\put(257,    -190){5}\put(275,    -200){6}
\put(290,    -210){7}\put(305,    -220){8}
\put(193,    -21){6}\put(195,    -33){3}\put(212,    -64){3}\put(212,    -84){4}\put(212,    -112){2}
\put(198,    1){7}\put(215,    -18){3}\put(215,    1){3}\put(230,    1){4}
\put(242,    14){5}\put(245,    -2){2}\put(264,    -64){2}\put(280,    -73){4}\put(296,    -81){3}
\end{picture}\\ \\ \\ \\ \\ \\ \\ \\ \\ \\ \\ \\ \\ \\ \\

\section {Generalized Dynkin diagrams of finite Cardan type}
The list below is the generalied Dynkin diagrams of finite Cartan type with quantum number $q.$\\

$$ \hbox {Figure 2}$$\\

 \noindent
{\rm (i)} $B_n$:
\ \ \ \ \ \ \ \ \  \ \ \
$\begin{picture}(100,     15)
\put(27,     1){\makebox(0,    0)[t]{$\bullet$}}
\put(60,     1){\makebox(0,     0)[t]{$\bullet$}}
\put(93,    1){\makebox(0,    0)[t]{$\bullet$}}
\put(159,     1){\makebox(0,     0)[t]{$\bullet$}}
\put(192,    1){\makebox(0,     0)[t]{$\bullet$}}
\put(225,    1){\makebox(0,    0)[t]{$\bullet$}}
\put(28,     -1){\line(1,     0){33}}
\put(61,     -1){\line(1,     0){30}}
\put(130,    -1){\makebox(0,    0)[t]{$\cdots\cdots\cdots\cdots$}}
\put(160,    -1){\line(1,     0){30}}
\put(193,     -1){\line(1,     0){30}}
\put(22,    -15){1}
\put(58,     -15){2}
\put(91,     -15){3}
\put(157,     -15){n-2}
\put(191,     -15){n-1}
\put(224,     -15){n}
\put(22,    10){$q^{2}$}
\put(58,     10){$q^{2}$}
\put(91,     10){$q^{2}$}
\put(157,     10){$q^{2}$}
\put(191,     10){$q^{2}$}
\put(224,     10){$q$}
\put(40,     5){$q^{-2}$}
\put(73,     5){$q^{-2}$}
\put(172,    5){$q^{-2}$}
\put(205,     5){$q^{-2}$}
\end{picture}  $  \hfill ($q^2 \not=1$).\\

 \noindent  {\rm (ii)} $C_n$:
\ \ \ \ \ \ \ \ \  \ \ \
$\begin{picture}(100,     15)
\put(27,     1){\makebox(0,    0)[t]{$\bullet$}}
\put(60,     1){\makebox(0,     0)[t]{$\bullet$}}
\put(93,    1){\makebox(0,    0)[t]{$\bullet$}}
\put(159,     1){\makebox(0,     0)[t]{$\bullet$}}
\put(192,    1){\makebox(0,     0)[t]{$\bullet$}}
\put(225,    1){\makebox(0,    0)[t]{$\bullet$}}
\put(28,     -1){\line(1,     0){33}}
\put(61,     -1){\line(1,     0){30}}
\put(130,    -1){\makebox(0,    0)[t]{$\cdots\cdots\cdots\cdots$}}
\put(160,    -1){\line(1,     0){30}}
\put(193,     -1){\line(1,     0){30}}
\put(22,    -15){1}
\put(58,     -15){2}
\put(91,     -15){3}
\put(157,     -15){n-2}
\put(191,     -15){n-1}
\put(224,     -15){n}
\put(22,    10){$q$}
\put(58,     10){$q$}
\put(91,     10){$q$}
\put(157,     10){$q$}
\put(191,     10){$q$}
\put(224,     10){$q^2$}
\put(40,     5){$q^{-1}$}
\put(73,     5){$q^{-1}$}
\put(172,    5){$q^{-1}$}
\put(205,     5){$q^{-2}$}
\end{picture}$  \hfill ($q^2 \not=1$).\\

 \noindent {\rm (iii)} $F_4$:
 \ \ \ \ \ \ \ \ \  \ \ \
$\begin{picture}(100,     15)
\put(27,     1){\makebox(0,    0)[t]{$\bullet$}}
\put(60,     1){\makebox(0,     0)[t]{$\bullet$}}
\put(93,    1){\makebox(0,    0)[t]{$\bullet$}}
\put(126,     1){\makebox(0,   0)[t]{$\bullet$}}
\put(28,     -1){\line(1,     0){33}}
\put(61,     -1){\line(1,     0){30}}
\put(94,    -1){\line(1,     0){30}}
\put(22,    -15){1}
\put(58,     -15){2}
\put(91,     -15){3}
\put(124,    -15){4}
\put(22,    10){$q^2$}
\put(58,     10){$q^2$}
\put(91,     10){$q$}
\put(124,     10){$q$}
\put(40,     5){$q^{-2}$}
\put(73,     5){$q^{-2}$}
\put(106,    5){$q^{-1}$}
\end{picture}$  \hfill ($q^2 \not=1$).\\

 \noindent {\rm (iv)} $G_2$:
\ \ \ \ \ \ \ \ \  \ \ \
$\begin{picture}(100,     15)
\put(27,     1){\makebox(0,    0)[t]{$\bullet$}}
\put(60,     1){\makebox(0,     0)[t]{$\bullet$}}
\put(28,     -1){\line(1,     0){33}}
\put(22,    -15){1}
\put(58,     -15){2}
\put(22,    10){$q$}
\put(58,     10){$q^3$}
\put(40,     5){$q^{-3}$}
\end{picture}$    \hfill ($q^3 \not=1, q\not= -1$).\\ \\

Let  $q_{ij} q_{ji} =q^{-1}$ and  $q_{ii} =q \not= 1$  for $1\le i, j \le n$ with $i\not= j$ in the following diagrams.

\noindent {\rm (v)} $A_n (n\ge 1)$ $
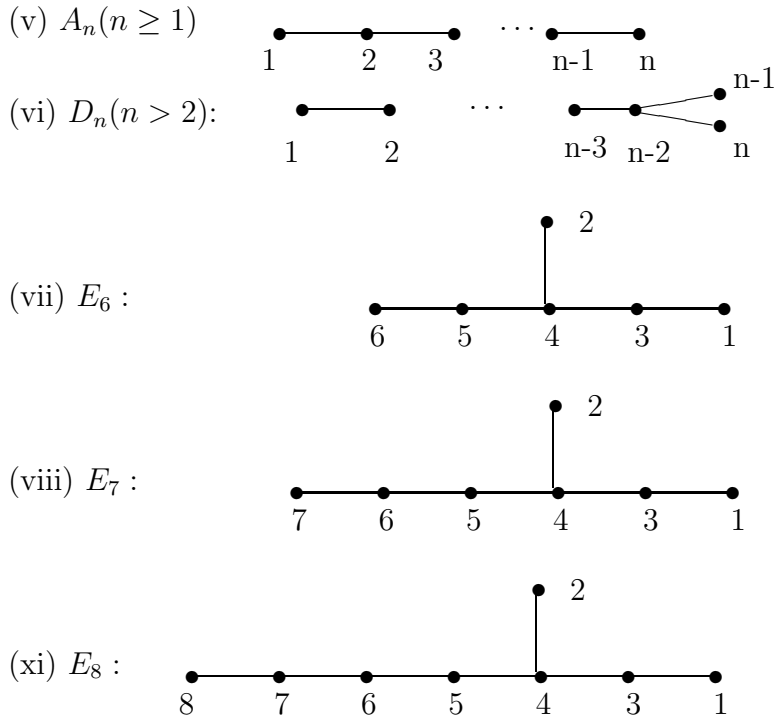
\begin{picture}(100,    15) \put(27,   1){\makebox(0,
0)[t]{$\bullet$}} \put(60,    1){\makebox(0,   0)[t]{$\bullet$}}
\put(93,    1){\makebox(0,    0)[t]{$\bullet$}} \put(28,   -1){\line(1,
0){33}} \put(61,    -1){\line(1,    0){30}} \put(20,    -15){1} \put(58,
-15){2} \put(83,    -15){3} \put(110,   0){$\dots$}\put(130,
1){\makebox(0,   0)[t]{$\bullet$}}\put(163,   1){\makebox(0,
0)[t]{$\bullet$}}\put(129,   -1){\line(1,   0){33}} \put(130,   -15){n-1}
\put(163,  -15){n}
\end{picture}$ \\

\noindent {\rm (vi)}
$D_n (n>2)$: $\begin{picture}(200,    15)  \put(130,    7){\makebox(0,
0)[t]{$\bullet$}} \put(153,    7){\makebox(0,    0)[t]{$\bullet$}}
 \put(185,    1){\makebox(0,    0)[t]{$\bullet$}}
\put(185,   13){\makebox(0,    0)[t]{$\bullet$}}  \put(126,    -13){n-3}
\put(150,    -15){n-2} \put(190,    -13){n} \put(190,    13){n-1}
\put(152,   5){\line(6,    1){30}} \put(182,    -1){\line(-6,    1){30}}

 \put(128,    5){\line(1,    0){25}}

 \put(27,   7){\makebox(0,   0)[t]{$\bullet$}} \put(60,    7){\makebox(0,
0)[t]{$\bullet$}} \put(28,   5){\line(1,   0){33}}  \put(20,    -15){1}
\put(58,   -15){2}  \put(90,   5){$\dots$}

\end{picture}$ \\ \\ \\

\noindent {\rm (vii)} $E_6:$\begin{picture}(100,    15) \put(158,   34){\makebox(0,
0)[t]{$\bullet$}} \put(93,   1){\makebox(0,   0)[t]{$\bullet$}}
\put(126,    1){\makebox(0,   0)[t]{$\bullet$}} \put(159,
1){\makebox(0,   0)[t]{$\bullet$}} \put(192,   1){\makebox(0,
0)[t]{$\bullet$}} \put(225,   1){\makebox(0,   0)[t]{$\bullet$}}
\put(157,   1){\line(0,   1){30}}  \put(94,   -1){\line(1,    0){30}}
\put(127,   -1){\line(1,   0){30}} \put(160,   -1){\line(1,    0){30}}
\put(193,   -1){\line(1,   0){30}}   \put(91,   -15){6} \put(124,
-15){5} \put(157,   -15){4} \put(191,    -15){3}\put(170,   28){2}
\put(224,   -15){1}

\end{picture}\\ \\ \\

\noindent {\rm (viii)} $E_7:$\begin{picture}(100,    15) \put(158,   34){\makebox(0,
0)[t]{$\bullet$}} \put(60,    1){\makebox(0,    0)[t]{$\bullet$}}
\put(93,   1){\makebox(0,   0)[t]{$\bullet$}} \put(126,  1){\makebox(0,
0)[t]{$\bullet$}} \put(159,    1){\makebox(0,   0)[t]{$\bullet$}}
\put(192,   1){\makebox(0,    0)[t]{$\bullet$}} \put(225,
1){\makebox(0,   0)[t]{$\bullet$}}  \put(157,   1){\line(0,   1){30}}
\put(61,   -1){\line(1,    0){30}} \put(94,   -1){\line(1,   0){30}}
\put(127,   -1){\line(1,    0){30}} \put(160,   -1){\line(1,   0){30}}
\put(193,   -1){\line(1,    0){30}}  \put(58,    -15){7} \put(91,
-15){6} \put(124,   -15){5} \put(157,    -15){4} \put(191,    -15){3}
\put(224,   -15){1} \put(170,   28){2}
\end{picture}\\ \\ \\

\noindent {\rm (xi)} $E_8:$\begin{picture}(100,    15) \put(158,   34){\makebox(0,
0)[t]{$\bullet$}} \put(27,    1){\makebox(0,    0)[t]{$\bullet$}}
\put(60,    1){\makebox(0,    0)[t]{$\bullet$}} \put(93,
1){\makebox(0,   0)[t]{$\bullet$}} \put(126,    1){\makebox(0,
0)[t]{$\bullet$}} \put(159,    1){\makebox(0,    0)[t]{$\bullet$}}
\put(192,   1){\makebox(0,    0)[t]{$\bullet$}} \put(225,
1){\makebox(0,   0)[t]{$\bullet$}} \put(28,    -1){\line(1,    0){33}}
\put(157,   1){\line(0,    1){30}} \put(61,    -1){\line(1,    0){30}}
\put(94,   -1){\line(1,    0){30}} \put(127,    -1){\line(1,    0){30}}
\put(160,   -1){\line(1,    0){30}} \put(193,    -1){\line(1,    0){30}}
\put(22,   -15){8} \put(58,    -15){7} \put(91,    -15){6} \put(124,
-15){5} \put(157,    -15){4} \put(191,    -15){3} \put(224,    -15){1}
\put(170,   28){2}
\end{picture}\\ \\

$A_n$,  $D_n$, $E_6$, $E_7$, $E_8$  are called the finite laced Cartan types.

\section {Other}

\vskip.1in
It is possible that $ \dim \mathfrak B(V)) = \infty$ and $ GK ( \mathfrak B(V)) < \infty$. For example, let $V= kx_1$ with $q_{11}=1$. By simple computation, $\mathfrak B(V) = k[x]$ with $\dim (\mathfrak B(V)) = \infty$ and $ GK ( \mathfrak B(V)) < \infty$ by \cite  [Pro. 1.15]{MR87}.

\begin {Lemma}\label {pp7.1} If $ A \rtimes \mathbb{S}_n $ is a subgroup of  $\mathbb Z_2 ^n \rtimes \mathbb{S}_n$, then  $A= 0$, or $ A=\mathbb Z_2^n,$ or $ A = \{ 0, (1, 1, \cdots, 1)\}$, or
$ A= \{a \in \mathbb Z_2 ^n \mid  \sum \limits  _{i=1}^n a_i \hbox {  is  even } \}$.
\end {Lemma}

\noindent {\it Proof.} It is clear that $ A \rtimes \mathbb{S}_n $ is a subgroup of  $\mathbb Z_2 ^n \rtimes \mathbb{S}_n$ when $A$ is a set in this Lemma. Conversely, if   $A$ is not 1th set, 2th set, 3th set in this Lemma,  then there exists  $ 0 \not= a \in A$ such that $a\not= (1, 1, \cdots, 1).$

 If $a_{i_1}, a_{i_2}, \cdots , a_{i_r}$ are $1$ and other  components of $a$ are $0$ with $0< r <n$. Let $\tau \in \mathbb S_n$ sending $j$ to $i_j$ for $1\le j \le r$.  We  have  $b= :\tau \cdot a  = (  \stackrel {r} {\overbrace{1, 1, \cdots, 1}}, 0, \cdots, 0)\in A$. Similarlty, $c:= ( 0,  \stackrel {r} {\overbrace{1, 1, \cdots, 1}}, 0, \cdots, 0)\in A$. Thus  $b-c = ( 1,  \stackrel {r-1} {\overbrace{0, 0, \cdots, 0}}, 1, 0, \cdots, 0)\in A$ and  $(  1, 1, 0,  \cdots, 0)\in A$. We keep on doing this, we have
$(  1, 1, 1, 1, 0   \cdots, 0), \cdots,  (  \stackrel {2k} {\overbrace{1, 1, \cdots, 1}}, 0, \cdots, 0)\in A$. This implies $ A= \{a \in \mathbb Z_2 ^n \mid  \sum \limits  _{i=1}^n a_i \hbox {  is  even } \}$. \hfill $\Box$

\begin {Lemma}\label {pp7.1''}  If $(V, \alpha, \delta)$ is $H$- {\rm YD} module, then
the braiding $C$ of $V$ can be lifted to  $\mathfrak B(V) = T(V)/I.$

Here $t_{ij}$ and $g_{ij}$ are in Section \ref  {fk conjecture}. \end {Lemma}

\noindent {\it Proof.}  By \cite [Section 2.1]{AS02}, $I$ is a {\rm YD} module over $H.$

\begin {eqnarray*} C((u+I) \otimes
 (v +I)) &=& \sum ({u+I}) _{(-1)} \cdot (v+I) \otimes  \otimes ({u+I}) _{(0)} \\
&=&\sum ({u} _{(-1)} \cdot v  +I)  \otimes (u_{(0)}
+I ) \\
&=& \sum _{i=1} ^m (u_i +I) \otimes (v_i + I),
 \end {eqnarray*}
 where $C(u\otimes v) =\sum _{i=1} ^m u_i  \otimes v_i. $ \hfill $\Box$

\begin {Lemma}\label {6.1} By decomposition (\ref {pe4.5.1}) of $G^{}(12)$, there exists  $\zeta_{st}(t_{ij}) \in (\mathbb S_n) ^{(12)}$ such that $g_{st}t_{uv}=\zeta_{st}(t_{uv})g_{s't'}$  for any $1\le u, v, s, t\le n.$ Then
for any $2 \le  j, k, j_1, k_1 \le n,$
$\zeta_{12}(t_{12}) = (12)$ since ${\rm id}(12)=(12){\rm id}$;

$\zeta_{1j}(t_{12}) = (12)$ since $(2j)(12)=(12)(1j)$;

$\zeta_{2j}(t_{12}) = (kj)$ since $(1j)(12)=(12)(2j)$;

 $\zeta_{kj}(t_{12}) = (1j)$ since $(1k)(2j)(12)=(kj)(1k)(2j)$;

 $\zeta_{12}(t_{1j}) = {\rm id}$ since ${\rm id}(1j)={\rm id} (1j)$;

 $\zeta_{1j}(t_{1j}) = (12)$ since $(2j)(1j)=(12)(2j)$;

$\zeta_{1j_1}(t_{1j}) ={\rm id}$ since $(2j_{1})(1j)={\rm id}(1j)(2j_{1}), \ \ j<j_{1}$;

$\zeta_{1j_1}(t_{1j}) = (jj_{1})(12)$ since $(2j_{1})(1j)=(jj_{1})(12)(1j_{1})(2j),\ \ j>j_{1}$;

$\zeta_{2j}(t_{1j}) = {\rm id}$ since $(1j)(1j)={\rm id}\ {\rm id}$;

$\zeta_{2j_1}(t_{1j}) = (jj_{1})$ since $(1j_{1})(1j)=(jj_{1})(1j_{1})$, $j \not= j_1$;

$\zeta_{kj}(t_{1j}) = (12)(kj)$ since $(1k)(2j)(1j)=(12)(kj)(2k)$, $j \not= k$;

$\zeta_{kj_1}(t_{1j}) = {\rm id}$ since $(1k)(2j_{1})(1j)={\rm id}(2j_{1}),\ \ \ j=k$;

$\zeta_{kj_1}(t_{1j}) = (kj)$ since $(1k)(2j_{1})(1j)=(kj)(1k)(2j_{1}),\ \ j\neq j_{1}$,   $j \not= k$;

$\zeta_{12}(t_{2j}) = {\rm id}$ since ${\rm id}(2j)={\rm id}(2j)$;

$\zeta_{1j}(t_{2j}) = {\rm id}$ since $(2j)(2j)={\rm id}\ {\rm id}$;

$\zeta_{1j_1}(t_{2j}) = (jj_{1})$ since $(2j_{1})(2j)=(jj_{1})(2j_{1})$, $j \not= j_1$;

$\zeta_{2j}(t_{2j}) = (12)$ since $(1j)(2j)=(12)(1j)$;

$\zeta_{2j_1}(t_{2j}) = (jj_{1})(12)$ since $(1j_{1})(2j)=(jj_{1})(12)(1j)(2j_{1}),\ \ j<j_{1}$;

$\zeta_{2j_1}(t_{2j}) = {\rm id}$ since $(1j_{1})(2j)={\rm id}(1j_{1})(2j),\ \ j>j_{1}$;

$\zeta_{kj}(t_{2j}) = {\rm id}$ since $(1k)(2j)(2j)={\rm id}(1k)$;

$\zeta_{kj_1}(t_{2j}) = (12)(kj_{1})$ since $(1k)(2j_{1})(2j)=(12)(kj_{1})(1j_{1}),\ \ j=k$,  $k < j_1$;

$\zeta_{kj_1}(t_{2j}) = (j_{1}j)$ since $(1k)(2j_{1})(2j)=(j_{1}j)(2j_{1})(1k),\ \ j\neq j_{1}, ~ j\neq k$,  $k < j_1$;

$\zeta_{12}(t_{kj}) = (kj)$ since ${\rm id}(kj)={\rm id}(kj)$;

$\zeta_{1j}(t_{kj}) = (kj)$ since $(2j)(kj)=(kj)(2k)$;

$\zeta_{1k}(t_{kj}) = (kj)$ since $(2k)(kj)=(kj)(2j)$;

$\zeta_{1j_1}(t_{kj}) = (kj)$ since $(2j_{1})(kj)=(kj)(2j_{1}),\ \ k\neq j_{1}, ~ j\neq j_{1}$;

$\zeta_{kj}(t_{kj}) = (kj)$ since $(1j)(kj)=(kj)(1k)$;

$\zeta_{2k}(t_{kj}) = (kj)$ since $(1k)(kj)=(kj)(1j)$;

$\zeta_{2k}(t_{kj}) = (kj)$ since $(1j_{1})(kj)=(kj)(1j_{1}),\ \ k\neq j_{1},~  j\neq j_{1}$;

$\zeta_{2j_1}(t_{kj}) = (12)$ since $(1k)(2j)(kj)=(12)(1k)(2j)$;

$\zeta_{kj}(t_{kj}) = (kj)$ since $(1k_{1})(2j_{1})(kj)=(kj)(1k)(2j_{1}),\ \ k_{1}=j$, $k < j_1$;

$\zeta_{k_1j_1}(t_{kj}) = (kj)$ since $(1k_{1})(2j_{1})(kj)=(kj)(1k_{1})(2k),\ \ k_{1}<k, ~ j_{1}=j$;

$\zeta_{k_1j_1}(t_{kj}) =(12) (jkk_1)$ since $(1k_{1})(2j_{1})(kj)=(12)(jkk_{1})(1k)(2k_{1}),\ \ k_{1}>k, ~ j_{1}=j$;

$\zeta_{k_1j_1}(t_{kj}) = (kj)$ since $(1k_{1})(2j_{1})(kj)=(kj)(1j)(2j_{1}),\ \ j_{1}>j, ~ k_{1}=k$;

$\zeta_{k_1j_1}(t_{kj}) = (12)(kjj_{1})$ since $(1k_{1})(2j_{1})(kj)=(12)(kjj_{1})(1j_{1})(2j),\ \ j_{1}<j, ~ k_{1}=k$;

$\zeta_{k_1j_1}(t_{kj}) = (kj)$ since $(1k_{1})(2j_{1})(kj)=(kj)(1k_{1})(2j),\ \ k_{1}\neq j, ~ j_{1}=k$,  $k_1 < j$;

$\zeta_{k_1j_1}(t_{kj}) = (kj)$ since $(1k_{1})(2j_{1})(kj)=(kj)(1k_{1})(2j_{1}),\ \ k_{1}\neq k, ~ k_{1}\neq j, ~ j_{1}\neq j, ~ j_{1}\neq k$.
`

\end {Lemma}

\section {Relation between bi-one arrow Nichols algebras and $\mathfrak{B}({\mathcal O}_s, \rho)$}

In this section it is shown that
bi-one arrow Nichols algebras and
$\mathfrak{B}({\mathcal O}_s, \rho)$ introduced in \cite {Gr00, AZ07,  AFZ09} are the
same  up to isomorphisms.

For any ${\rm RSR} (G, r, \overrightarrow \rho, u)$, we can
construct an arrow Nichols algebra
 $\mathfrak{B} (kQ_1^1, ad (G, r, \overrightarrow{\rho},
$ $u))$ ( see \cite [Pro. 2.4] {ZCZ08}), written as $\mathfrak{B} (G,
r, \overrightarrow{\rho}, $ $ u)$ in short.
Let us recall the  precise description of arrow {\rm YD} module.
For an ${\rm RSR}(G, r, \overrightarrow \rho, u)$ and a $kG$-Hopf
bimodule $(kQ_1^c, G, r, \overrightarrow {\rho}, u)$ with the module
operations $\alpha^-$ and $\alpha^+$, define a new left $kG$-action
on $kQ_1$ by
$$g\rhd x:=g\cdot x\cdot g^{-1},\ g\in G, x\in kQ_1,$$
where $g\cdot x=\alpha^-(g\otimes x)$ and $x\cdot
g=\alpha^+(x\otimes g)$ for any $g\in G$ and $x\in kQ_1$. With this
left $kG$-action and the original left (arrow) $kG$-coaction
$\delta^-$, $kQ_1$ is a  Yetter-Drinfeld   $kG$-module. Let  $Q_1^1:=\{a\in Q_1 \mid
s(a)=1\}$, the set of all arrows with starting vertex $1$. It is
clear that $kQ_1^1$ is a  Yetter-Drinfeld   $kG$-submodule of $kQ_1$, denoted by
$(kQ_1^1, ad(G, r, \overrightarrow {\rho}, u))$, called the arrow {\rm YD} module.

\begin {Lemma} \label {1.1'} For any $s\in G$ and  $\rho \in  \widehat
{G^s}$, there exists a bi-one arrow Nichols algebra $\mathfrak{B} (G,
r, \overrightarrow{\rho}, u)$ such that
$$\mathfrak{B}({\mathcal O}_s, \rho)
\cong \mathfrak{B} (G, r, \overrightarrow{\rho}, u)$$ as graded
braided Hopf algebras in $^{kG}_{kG} \! {\mathcal YD}$.

\end {Lemma}
{\bf Proof.}  Assume that $V$ is the representation space of $\rho$
 with $\rho (g) (v)= g\cdot v$ for any $g\in G, v\in V$. Let $C = {\mathcal O_s}$,
$r = r_C C$, $r_C ={\rm deg } \rho$, $u(C)
= s$, $I_C(r, u)= \{1 \}$ and  $ (v)\rho _C^{(1)}(h)= \rho (h^{-1})(v)$
for any $h\in G$, $v\in V$. We get a bi-one arrow Nichols algebra $
\mathfrak{B}( G, r, \overrightarrow{\rho}, u)$.

We now  only need to show that $M ({\mathcal O}_s, \rho) \cong
(kQ_1^1, ad( G, r, \overrightarrow{\rho}, u)) $ in $^{kG}_{kG} \!
{\mathcal YD}$. We recall the notation in \cite [Proposition
1.2]{ZCZ08}. Assume  $J_C(1) = \{1, 2, \cdots, n\}$ and $X_C^{(1)} =
V$ with basis $\{x_C^{(1, j)} \mid j =1, 2, \cdots, n\}$ without
loss of  generality. Let $v_j $ denote $x_C^{(1,j)}$ for convenience. In
fact,
 the  left and right coset decompositions of $G^s$ in $G$ are
\begin {eqnarray} \label {e1.1.1}
G =\bigcup_{i=1}^m g_{i} G^s    \ \ \hbox {and } \ \    G
&=&\bigcup_{i=1} ^m G^sg_{i} ^{-1} \ \ ,
\end {eqnarray} respectively.

 Let $\psi$ be a map from $M({\mathcal O}_s,
\rho)$ to $(kQ_1^1, {\rm ad } (G, r, \overrightarrow{\rho}, u))$ by
sending $ g_i v_j$ to $a _{t_i, 1}^{(1, j)}$ for any $1\le i\le m,
1\le j \le n$. Since the dimension is $mn$, $\psi$ is a
bijective. See
\begin {eqnarray*}
\delta ^- (\psi (g_i v_j)) &=& \delta ^- (a _{t_i, 1}^{(1,j)})\\
&=& t_i \otimes a_{t_i, 1}^{(1, j)} = (id \otimes \psi) \delta ^-
(g_iv_j). \end {eqnarray*} Thus $\psi$ is a $kG$-comodule
homomorphism. For any $h\in G$, assume  $hg_i = g_{i'} \gamma $ with
$\gamma \in G^s$. Thus $g_i ^{-1} h^{-1} = \gamma ^{-1}
g_{i'}^{-1}$, i.e. $\zeta _i (h^{-1})= \gamma ^{-1}$, where
$\zeta_i$ was defined  in \cite [(0.3)]{ZZC07}. Since $\gamma \cdot
x^{(1,j)} \in V$, there exist $k_{C, h^{-1}}^{(1, j, p)}\in k $,
$1\le p \le n$, such that $\gamma \cdot x^{(1,j)} =\sum _{p=1}^n
k_{C, h^{-1}} ^{(1, j, p)} x^{(1, p)}$. Therefore
\begin {eqnarray}
\label {e1.11}x^{(1,j)}\cdot \zeta _i (h^{-1}) &=& \gamma \cdot
x^{(1,j)} \ \ ( \hbox {by definition of } \rho _C^{(1)} ) \nonumber
\\ &=&\sum _{p=1}^n k_{C, h^{-1}} ^{(1, j, p)} x^{(1,
p)}.\end {eqnarray}

See
\begin {eqnarray*}
 \psi (h \cdot g_iv_j) &=& \psi ( g _{i'}(\gamma v_j)) \\
&=& \psi (g_{i'} (\sum _{p=1}^n k_{C, h^{-1}} ^{(1, j, p)} v_p) )\\
 &= &\sum _{p=1}^n k_{C, h^{-1}} ^{(1, j, p)} a _{t_{i'}, 1} ^{(1,
p)}
\end {eqnarray*}
and
\begin {eqnarray*}
h \rhd (\psi (g_i v_j)) &=& h \rhd (a _{t_i, 1}^{(1,j)})\\
&=& a_{ht_i, h}^{(1, j)} \cdot h^{-1} \\
&=& \sum _{p =1}^n k _{C, h^{-1}} ^{(1, j, p)} a_{t_{i'}, 1}^{(1,
p)} \ \ (\hbox {by \cite [Pro.1.2]{ZCZ08} and }  (\ref {e1.11})).
\end {eqnarray*} Therefore
$\psi $ is a $kG$-module homomorphism. $\Box$

 Therefore  we write  Hopf bimodule $(kQ_1^1, ad( G, r, \overrightarrow{\rho}, u)) $ in the proof above as $M ({\mathcal O}_s, {\rm ad }(\rho))$ and  Nichols algebra $\mathfrak{B} (G,
r, \overrightarrow{\rho}, u)$ in the lemma above as
$\mathfrak{B}({\mathcal O}_s,  {\rm ad } (\rho))$  in short.

\begin {Remark} \label {1.2'} The representation $\rho$ in   $\mathfrak{B}({\mathcal O}_s, \rho)$   introduced in \cite {Gr00, AZ07} and
$\rho _C^{(i)}$ in {\rm RSR} are different.  $\rho (g)$ acts on its representation
space from  the left and $\rho _C^{(i)} (g)$
acts on its representation space from  the right.
\end {Remark}

Otherwise, when $\rho = \chi$ is a one dimensional representation, then $(kQ_1^1, ad( G, r, \overrightarrow{\rho}, u)) $ is PM (see \cite [Def. 1.1] {ZZC07}). Thus the formulae  are available in   \cite [Lemma 1.9] {ZZC07}. That is, $g \cdot a _{t} = a _{g t_i, g}$, $  a _{ t_i} \cdot g = \chi (\zeta _{i}(g)) a_{t_ig, g}$.

\section*{Acknowledgement}
Y.Z.Z. was supported by the Australian Research Council through Discovery Project DP140101492.

\end {document}